\documentclass[12pt]{amsart}
\usepackage{mathtools}
\usepackage{etoolbox}
\usepackage{calligra}
\usepackage{amsthm}
\usepackage{amsmath}
\usepackage{amssymb}
\usepackage{amscd}
\usepackage{bbm}
\usepackage[all]{xy}
\usepackage{mathrsfs}
\usepackage[top=72pt,bottom=96pt,left=72pt,right=72pt]{geometry}
\newtheorem{Lemma}{Lemma}[section]
\def\YM{\mathrm{YM}}
\def\qS{\mathscr{S}}
\def\dirac{\mathrm{D}}
\newtheorem{Remark}[Lemma]{Remark}

\newtheorem{Theorem}[Lemma]{Theorem}

\newtheorem{Definition}[Lemma]{Definition}
\newtheorem{Proposition}[Lemma]{Proposition}
\numberwithin{equation}{section}
\def\F{\mathfrak{F}}
\def\f{\mathfrak{f}}

\def\qS{\mathscr{S}}
\def\C{\mathbb{C}}
\def\torus{\mathbb{T}}
\def\E{\mathbb{E}}
\def\F{\mathfrak{F}}
\def\qGG{\mathfrak{qGG}}
\def\qgs{q\mathcal{GS}}
\def\Hor{\mathrm{Hor}}
\def\dvol{\mathrm{dvol}}
\def\id{\mathrm{id}}
\def\Ker{\mathrm{Ker}}
\def\Im{\mathrm{Im}}
\def\YM{\mathrm{YM}}

\def\H{\mathrm{Hor}}
\def\Mor{\textsc{Mor}}
\def\N{\mathbb{N}}
\def\Obj{\textsc{Obj}}

\def\R{\mathbb{R}}
\def\Z{\mathbb{Z}}
\def\C{\mathbb{C}}
\def\C{\mathbb{C}}

\def\Rep{\mathbf{Rep}}
\def\triv{\mathrm{triv}}

\def\l{\mathrm{L}}

\def\ad{\mathrm{ad}}
\def\Ad{\mathrm{Ad}}
\def\c{\mathrm{c}}
\def\r{\mathrm{R}}
\def\u{\mathbbm{u}}
\def\U{\mathcal{U}}
\def\D{\mathcal{D}}
\def\G{\mathcal{G}}

\def\T{\mathcal{T}}

\def\Vert{\mathrm{Ver}}
\def\V{\mathrm{V}}

\begin{document}
\date{\today}
\title{Yang--Mills--Connes Theory and Quantum Principal $SU(N)$--Bundles over Toric Non--Commutative Manifolds}
\author{Gustavo Amilcar Salda\~na Moncada}
\address{Gustavo Amilcar Salda\~na Moncada\\
CIMAT, Unidad Guanajuato}
\email{gamilcar@ciencias.unam.mx}
\begin{abstract}
This paper has two main objectives. The first one is to show that the Connes' formulation of Dirac theory can be applied in the framework of quantum principal bundles for any $n$--dimensional spectral triple, any quantum group, any quantum principal connection and any finite--dimensional corepresentation of the quantum group. The second objective is to demonstrate that, for quantum principal $SU(N)$--bundles over toric non--commutative manifolds, one can define a Yang–Mills functional that measures the squared norm of the curvature of a quantum principal connection, in contrast to the Yang–Mills functional proposed by Connes, which measures the squared norm of the curvature of a quantum linear connection compatible with the metric. An illustrative example based on the non--commutative $n$--torus and $U(1)$ is presented, highlighting the differences and similarities between the two functionals.

 \begin{center}
  \parbox{300pt}{\textit{MSC 2010:}\ 46L87, 58B99.}
  \\[5pt]
  \parbox{300pt}{\textit{Keywords:}\ Spectral Triple, quantum principal bundles, quantum principal connections.}
 \end{center}
\end{abstract}
\maketitle

\section{Introduction}

In reference \cite{con},  A. Connes presented a formulation of spinorial matter and Yang--Mills theory within the framework of non--commutative geometry. This formulation is based on the concepts of {\it spectral triples}, {\it quantum vector bundles}, and {\it quantum linear connections}, and it has been extensively studied ever since. In particular,  Connes defined the non--commutative Yang--Mills functional as the functional that measures the squared norm of the curvature of any quantum linear connection compatible with the metric.

On the other hand, in differential geometry,  the Yang--Mills functional is defined as the functional that measures the squared norm of the curvature of a principal connection. Following this idea, in references \cite{sald2,saldym} we developed another non--commutative Yang--Mills functional, which now measures the squared norm of the curvature of any quantum principal connection. However, a natural question arises from all this: are both non--commutative Yang--Mills functionals equivalent? If not, do they share any critical point?

This paper was developed as an attempt to provide a preliminary answer to these questions, while formulating a {\it physical} Yang--Mills theory  based on quantum principal bundles and incorporating the Dirac theory formulated by  Connes

Since one of our goals is to develop a {\it physical} Yang--Mills theory, we will restrict ourselves to quantum principal $\G$--bundles in which $\G$ is the canonical quantum group associated with the Lie group $SU(N)$, for $N$ $\in$ $\N$. This assumption is not trivial, it was made because, as far as the author knows, there are no {\it physical reasons} to believe that the structure group of the bundle modeling fundamental interactions in nature must be a purely non--commutative geometrical object. Nevertheless, it is worth mentioning that, from a {\it completely mathematical} point of view, the theory can in principle be applied for any quantum group.

To accomplish our purposes, this paper is divided into six sections. Following this introduction, Section 2 presents the main aspects of the theory of quantum groups, bicovariant first--order differential calculus, and the universal differential envelope $\ast$--calculus. In addition, in Section 2.1 we focus on showing how these concepts correctly generalize the notions of matrix compact Lie groups and the algebra of $\C$--valued differential forms on them. Finally, in Section 2.2 we introduce a $1$--dimensional non--standard first--order differential calculus on $U(1)$.

In the third section, we present the fundamentals of the theory of quantum principal bundles and quantum principal connections developed by  M. Durdevich in references \cite{micho1,micho2,micho3,stheve}, among other papers.  In Section 3.2, as an example, we provide a brief exposition of toric non--commutative manifolds and the corresponding quantum principal bundles over them, following references \cite{lan2,lan4}. 

In the fourth section, we present the theory of associated quantum vector bundles, the canonical Hermitian structure defined on them and the gauge quantum linear connection induced by a quantum principal connection. This section is based on \cite{sald2}, and the reader is encouraged to consult that reference for further details. It is worth mentioning that the structures commented above are mutually compatible; in other words, every gauge quantum linear connection is compatible with the canonical Hermitian structure. This compatibility provides a natural framework for incorporating Connes' formulation of Dirac theory into the setting of quantum principal bundles.

It is worth noticing that, up to this point, the paper is  entirely bibliographical (except for Section 2.3, Proposition \ref{teodio1} and Proposition \ref{grassmann}). We have deliberately taken the time to provide the reader with a clear and coherent context for the theory required in what follows, specially because Durdevich's formulation of quantum principal bundles and the theory of \cite{sald2} are not widely known. In this sense, the paper is reasonably self--contained.

The fifth section is the most important part of the paper, as it fulfills the main objectives of our work. Section 5.1 is about some general concepts needed. In Section 5.2 we define the gauge Dirac operator associated with an {\it arbitrary quantum principal bundle} and we present some of its properties related to the quantum gauge group. This is one of the purposes of this paper. Section 5.3 constitutes the core of the paper: here, we introduce the two non--commutative Yang--Mills functionals and we present the corresponding non--commutative geometrical Yang--Mills equations for quantum principal $SU(N)$--bundles over non--commutative toric manifolds. This is the second purpose of this paper.

In Section 6 we present an illustrative example that demonstrates the full application of our theory. Specifically, we consider quantum principal bundles over the non-commutative $n$--torus and explicitly compute the gauge Dirac operator associated with a specific quantum principal connection. In addition, we explicitly determine the space of critical points (i.e., the solutions of the corresponding Yang--Mills equation) for both non--commutative Yang--Mills functionals. Finally, the last section is about some concluding comments and Appendix~A provides a brief summary of the similarities and differences between Durdevich's formulation of quantum principal bundles and the Brzezi\'nski–Majid formulation presented in \cite{libro,br1,br2}, as well as the reasons why we have chosen to work with Durdevich's theory.

Connes' formulation of Yang--Mills theory is developed under the convention that the differential in all graded differential $\ast$--algebras is anti--Hermitian, i.e., $d(a^\ast) = - (da)^\ast$ \cite{con}. On the other hand, the theory of quantum principal bundles is constructed under the convention that the differential is Hermitian, i.e., $d(a^\ast) = (da)^\ast$. In this sense, in Sections 2, 3, and 4, we adapt the theory to the anti--Hermitian convention. The proofs of all results presented in these sections are analogous to those under the Hermitian convention, differing only by the inclusion of a factor $i = \sqrt{-1}$ in the differential and in the quantum germs map. Moreover, throughout this work we use Sweedler's notation.

In the whole text, we will use the words {\it classical} or {\it classical case} to refer to the differential geometry's framework. In the same way, we will use the word {\it quantum} to refer to non--commutative geometry's framework. Furthermore, all {\it classical} principal $H$--bundles $$\pi_X: Y\longrightarrow X $$ are considered under the assumptions that $X$, $Y$ are compact smooth manifolds and $H$ is a compact matrix Lie group.

\section{Compact Quantum Groups and Quantum Differential Forms}

Let us start with basics of Woronowicz's theory of quantum groups \cite{woro1, woro2}. However, we are going to use a more general graded differential $\ast$--algebra on quantum groups that the one presented in \cite{woro2}. This space is called the {\it universal differential envelope $\ast$--calculus} \cite{micho1,stheve}.

\subsection{Basic Topics}
A compact matrix quantum group (or for this paper, a quantum group) will be denoted by $\G$, and its dense $\ast$--Hopf (sub)algebra will be denoted by
\begin{equation}
    \label{2.f0}
    G^\infty:=(G,\cdot,\mathbbm{1},\Delta,\epsilon,S,\ast),
\end{equation}
where $\Delta$ is the coproduct, $\epsilon$ is the counit and $S$ is the antipode \cite{woro1}. The space $G^\infty$ shall be treated as the algebra of all {\it polynomial functions} defined on $\G$.  In the same way, a (right) $\G$--corepresentation on a $\C$--vector space $V$ is a linear map $$\delta^V: V \longrightarrow V \otimes G$$ such that
\begin{equation} 
\label{2.f1}
(\id_V\otimes \epsilon )\circ \delta^V\cong \id_V, \qquad (\id_V\otimes \Delta )\circ \delta^V=(\delta^V\otimes \id_G ) \circ \delta^V.
\end{equation}

\noindent We say that the corepresentation is {\it finite--dimensional} if $\dim_{\C}(V)< |\N|$. The map $\delta^V$ is often referred to as \emph{(right) coaction} of $\G$ on $V$.

Given two $\G$--corepresentations $\delta^V$, $\delta^W$,  a {\it corepresentation morphism} is a linear map 
\begin{equation}
\label{2.f3}
T:V\longrightarrow W \quad \mbox{ such that }\quad (T\otimes \id_G)\circ \delta^V=\delta^W \circ T.
\end{equation}
We will denote by 
\begin{equation}
\label{2.f4}
\Mor(\delta^V,\delta^W),\qquad \Obj(\Rep_{\G}).
\end{equation}
the set of all corepresentation morphisms between two corepresentations $\delta^V$, $\delta^W$, and  the set of all finite--dimensional $\G$--corepresentations, respectively. 

In reference \cite{woro1}, Woronowicz proved that for every finite--dimensional $\G$--corepresentation $\delta^V$, there always exists an inner product $\langle -|-\rangle$ on $V$ that makes $\delta^V$ unitary. Henceforth, we will assume that every finite--dimensional $\G$--corepresentation is unitary.  A proof of the following theorem can be found in reference \cite{woro1}.

\begin{Theorem}
\label{rep}
Let $\T$ be a complete set of mutually non--equivalent irreducible (necessarily finite--dimensional) $\G$--corepresentations with $\delta^\C_\triv$ $\in$ $\T$ (the trivial corepresentation on $\C$). For any $\delta^V$ $\in$ $\T$ that coacts on $V$,
\begin{equation}
\label{2.f8}
\delta^V(e_j)=\sum^{n_{V}}_{i=1} e_i\otimes g^{V}_{ij},
\end{equation}
where $\{ e_i\}^{n_{V}}_{i=1}$ is an orthonormal basis of $V$ and $\{g^{V}_{ij}\}^{n_{V}}_{i,j=1}$ $\subseteq$ $G$. Then $\{g^{V}_{ij}\}_{\delta^V,i,j}$ is a linear basis of $G$, where the index $\delta^V$ runs on $\T$ and $i$, $j$ run from $1$ to $n_{V}=\mathrm{dim}_\C(V)$.
\end{Theorem}

For every $\delta^V$ $\in$ $\T$, the set $\{g^{V}_{ij}\}^{n_{V}}_{ij=1}$ satisfies
\begin{equation}
\label{2.f9}
\begin{aligned}
    \Delta(g^{V}_{ij})= &\sum^{n_{V}}_{k=1} g^{V}_{ik}\otimes g^{V}_{kj},\quad  \epsilon(g^{V}_{ij})=\delta_{ij},\quad S(g^{V}_{ij})=g^{V\,\ast}_{ji}
\end{aligned}
\end{equation}
with $\delta_{ij}$ being the Kronecker delta, among other properties (\cite{woro1}).

In the rest of the paper, we will consider only bicovariant  first--order differential $\ast$--calculus (abbreviated ``$\ast$--FODC") over $\G$ (\cite{stheve}). This theory is standard in the context of quantum groups, and the reader is referred to \cite{stheve} for further details. In particular, in Sections 6, 11 of \cite{stheve} it is proven that every bicovariant $\ast$--FODC over $\G$ 
$$(\Gamma,d),\qquad d: G\longrightarrow \Gamma$$ (is isomorphic to one that) satisfies
 \begin{equation}
    \label{2.f9.4}
    \Gamma=G\otimes \mathfrak{qg}^\#, \qquad \mathfrak{qg}^\#:= {\Ker(\epsilon)\over \mathcal{R}},
\end{equation}
where $\mathcal{R}$ $\subseteq$ $\Ker(\epsilon)$ is a right $G$--ideal such that $S(\mathcal{R})^\ast\subseteq \mathcal{R}$ and $\Ad(\mathcal{R})\subseteq \mathcal{R}\otimes G$, with
\begin{equation}
    \label{2.f9.5}
    \Ad:G\longrightarrow G\otimes G, \qquad g \longmapsto g^{(2)}\otimes S(g^{(1)})\,g^{(3)}
\end{equation}
the right adjoint coaction of $G$. Recall that for this paper, we will work with the convention $$d(g^\ast)=-(dg)^\ast$$ for all $g$ $\in$ $G$. We will refer to $\mathfrak{qg}^\#$ as the {\it quantum dual Lie algebra} and the dimension of the $\ast$--FODC $(\Gamma,d)$ is defined as the dimension of $\C$--vector space $\mathfrak{qg}^\#$. This definition comes from the fact that, a linear basis of $\mathfrak{qg}^\#$ is a $G$--bimodule basis of $\Gamma$ (\cite{stheve}).

The space $\mathfrak{qg}^\#$ allows to consider the quantum germs map 
\begin{equation}
\label{2.f14}
\pi:G \longrightarrow \mathfrak{qg}^\#,\qquad 
g \longmapsto S(g^{(1)})\,dg^{(2)}=-(dS(g^{(1)}))g^{(2)}.
\end{equation}
The map $\pi$ has several useful properties, for example, the restriction map $\pi|_{\Ker(\epsilon)}$ is surjective and for all $g$ $\in$ $G$ we have (\cite{stheve})
\begin{equation}
\label{properties}
\ker(\pi)=\mathcal{R}\oplus \C\mathbbm{1},  \qquad   dg=g^{(1)}\pi(g^{(2)}), \qquad  \pi(g)^\ast=\pi(S(g)^\ast).
\end{equation}
 Furthermore, there is $\G$--corepresentation on $\mathfrak{qg}^\#$ given by 
\begin{equation}
\label{2.f15}
\ad: \mathfrak{qg}^\#\longrightarrow \mathfrak{qg}^\#\otimes G, \qquad \ad\circ \pi= (\pi\otimes \id_G)\circ \Ad.
\end{equation}
On the other hand, there is a right $G$--module structure on $\mathfrak{qg}^\#$ given by 
\begin{equation}
\label{2.f16}
\theta \diamondsuit  g:=\pi(hg-\epsilon(h)g))=S(g^{(1)})\,\theta\, g^{(2)}
\end{equation}
for every $\theta=\pi(h)$ $\in$ $\mathfrak{qg}^\#$; and it satisfies  $(\theta \diamondsuit g)^\ast=-\theta^\ast\diamondsuit S(g)^\ast$ (\cite{stheve}).

Let $(\Gamma,d)$ be a bicovariant $\ast$--FODC over $\G$. Consider
\begin{equation}
    \label{2.f12.1}
    \begin{aligned}
        \mathfrak{qg}^{\#\wedge}=\otimes^\bullet& \mathfrak{qg}^\#/A^\wedge, \qquad \otimes^\bullet\mathfrak{qg}^\#:=\bigoplus_k (\otimes^k\mathfrak{qg}^\#)\\
\otimes^k\mathfrak{qg}^\#:=&\underbrace{\mathfrak{qg}^\#\otimes\cdots\otimes \mathfrak{qg}^\#}_{k\; times}, \qquad \otimes^0\mathfrak{qg}^\#:=\C\,\mathbbm{1},
    \end{aligned}
\end{equation}
where $A^\wedge$ is the two--side ideal of $\otimes^\bullet\mathfrak{qg}^\#$ generated by elements 
\begin{equation}
    \label{2.f12.11}
    \pi(g^{(1)})\otimes \pi(g^{(2)})\qquad \mbox{ for all }\qquad g \,\in\, \mathcal{R}.
\end{equation}
According to \cite{micho1}, the space 
\begin{equation}
    \label{2.f12.3}
    \Gamma^{\wedge}= G\otimes \mathfrak{qg}^{\#\wedge}.
\end{equation}
has a structure of a graded differential $\ast$--algebra generated by its degree $0$ elements $\Gamma^{\wedge\, 0}=G$, where $$d:\Gamma^\wedge\longrightarrow \Gamma^\wedge $$ is the natural extension of $d:G\longrightarrow \Gamma=\Gamma^{\wedge\,1}$ in such a way that $\mathfrak{qg}^{\#\wedge}$ is a graded differential $\ast$--subalgebra and the following equation holds
\begin{equation}
    \label{2.f12.2}
    d\pi(g)=-\pi(g^{(1)})\pi(g^{(2)})
\end{equation}
for all $g$ $\in$ $G$.
The graded differential $\ast$--algebra $(\Gamma^\wedge,d,\ast)$ is called the universal differential envelope $\ast$--calculus. For more detail, see reference \cite{micho1}. It is worth mentioning that the right $G$--module structure of $\mathfrak{qg}^\#$ (see equation (\ref{2.f16})) can be extended to $\mathfrak{qg}^{\#\wedge}$ by means of 
\begin{equation}
    \label{2.f12.4}
    \mathbbm{1}\diamondsuit g=\epsilon(g),\quad (\theta_1\theta_2)\diamondsuit g=(\theta_1\diamondsuit g^{(1)})(\theta_2\diamondsuit g^{(2)}).
\end{equation}
for all $g$ $\in$ $H$, $\theta_1$, $\theta_2$ $\in$ $\mathfrak{qg}^\#$. 

Moreover, the $\ast$--Hopf algebra structure $(\Delta,\epsilon,S)$ of $G$  can be extended to a graded differential $\ast$--Hopf algebra structure 
\begin{equation}
\label{2.f3.8}
\Gamma^{\wedge\,\infty}:=(\Gamma^\wedge,\cdot,\mathbbm{1},\Delta,\epsilon, S,d,\ast).
\end{equation}
In other words, we can define maps
 \begin{equation}
\label{2.f3.5}
\Delta: \Gamma^\wedge \longrightarrow \Gamma^\wedge \otimes  \Gamma^\wedge, \qquad \epsilon: \Gamma^\wedge\longrightarrow \C,\qquad S:\Gamma^\wedge\longrightarrow \Gamma^\wedge.
\end{equation}
that are graded differential $\ast$--algebra morphisms and  $\Gamma^{\wedge\,\infty}$ satisfies the definition of $\ast$--Hopf algebra. In particular, we have 
\begin{equation}
    \label{coproduc.1}
    \Delta(\theta)=\mathbbm{1}\otimes \theta+\ad(\theta)
\end{equation}
for all $\theta$ $\in$ $\mathfrak{qg}^\#$. The reader is encourage to consult reference \cite{micho1} for details.

\subsection{Example: Compact Matrix Lie Groups with the Standard Differential Calculus}

    Let $H$ be a compact matrix Lie group ($H\subseteq M_k(\C))$), and let $\G$ be its associated quantum group \cite{woro1}. If 
    \begin{equation}
        \label{2.2.f1}
        G^\infty=(G,\cdot,\mathbbm{1},\Delta,\epsilon,S,\ast,)
    \end{equation}
     is its dense $\ast$--Hopf algebra, then $G$ is the $\ast$--algebra generated by the smooth $\C$--valued functions
    \begin{equation}
    \label{2.2.f2}
    \begin{aligned}
    w_{ij}: H \longrightarrow \C, \qquad A=(a_{ij}) \longmapsto a_{ij}.
    \end{aligned}
\end{equation}
    The Hopf algebra structure $(\Delta,\epsilon,S)$ of $G$ is defined as the pull--back of the group structure of $H$ \cite{woro1}. In particular, we have $\Delta(w_{ij})= \displaystyle \sum^n_{k=1} w_{ik}\otimes w_{kj}.$  

If we consider $$\mathcal{R}:=\Ker^2(\epsilon):=\{\displaystyle \sum^n_{i=1}a_i\,b_i \mid a_i,\, b_i \,\in\, \Ker(\epsilon) \, \mbox{ for some }\, n\,\in\, \N  \},$$ then, the $\ast$--FODC 
induced by equation \ref{2.f9.4} is bicovariant and it is given by (\cite{woro2,war,appendix})
\begin{equation}
    \label{ec.2.102}
    (\Gamma=G\otimes \mathfrak{h}^\#_\C,d),
\end{equation}
where
\begin{equation}
    \label{ec.2.103}
\mathfrak{h}^\#_\C={\Ker(\epsilon)\over \mathcal{R}}={\Ker(\epsilon)\over \Ker^2(\epsilon)}
\end{equation}
 is the complexification of the dual space of the Lie algebra $\mathfrak{h}$ of $H$, and $d=id'$ with $d'$ the de Rham differential of $H$. 
 
 As is proven in \cite{appendix}, the quantum germs map 
\begin{equation*}
    \begin{aligned}
\pi: G \longrightarrow \mathfrak{h}^\#_\C, \qquad   
          g \longmapsto S(g^{(1)})\,dg^{(2)}
    \end{aligned}
\end{equation*}
is the differential of $g$ et $e$, i.e., $\pi(g)=(dg)_e$ and the right $G$--module structure on $\mathfrak{h}^\#_\C$  given in equation (\ref{2.f16}) can be simplified to (\cite{appendix})
 \begin{equation}
     \label{2.2.f4}
     \pi(h)\diamondsuit g=\epsilon(g)\pi(h).
 \end{equation}
 
Now, in light of \cite{appendix}, the universal differential envelope $\ast$--calculus of $(\Gamma,d)$ is given by
\begin{equation}
    \label{2.2.f5}
    (\Gamma^\wedge=G\otimes \bigwedge \mathfrak{h}^\#_\C,d,\ast),
\end{equation}
where $\bigwedge \mathfrak{h}^\#_\C$ is the exterior algebra of $\mathfrak{h}^\#_\C$ and $d=id'$ with $d'$ the de Rham differential. Hence $(\Gamma^\wedge,d,\ast)$ is a graded differential $\ast$--subalgebra of the graded differential $\ast$--algebra 
\begin{equation}
    \label{2.2.f5.1}
    (\Omega^\bullet_\C(H)= C^\infty_\C(H)\otimes \bigwedge \mathfrak{h}^\#_\C,d,\ast)
\end{equation}
of $\C$--valued differential forms of $H$ and it suffices to consider convergent sequences $\{ g_i\}^\infty_{i=1}$ $\subseteq$ $G$ in $(\Gamma^\wedge,d,\ast)$ to obtain $(\Omega^\bullet_\C(H),d,\ast).$ 

In accordance with \cite{appendix}, the graded differential $\ast$--Hopf algebra structure $$ \Gamma^{\wedge\,\infty}=(\Gamma^{\wedge},\cdot,\mathbbm{1},\Delta,\epsilon,S,d,\ast)$$ coincides with the pull--back of group structure of $TH_\C=H\times \mathfrak{h}_\C$ on $\C$--valued differential forms of $H$. Furthermore, if we consider the {\it classical} adjoint action of $H$ on $H$ 
$$\Ad^{\mathrm{class}}:H\times H\longrightarrow H,\qquad \Ad^{\mathrm{class}}(R,S):=\Ad^{\mathrm{class}}_S(R):=S^{-1}RS;$$ and the {\it classical} adjoint action of $H$ on $\mathfrak{h}_\C$
\begin{equation}
    \label{adjointclassical}
    \ad^{\mathrm{class}}:\mathfrak{h}_\C\times H\longrightarrow \mathfrak{h}_\C,\qquad \ad^{\mathrm{class}}(v,S):=\ad^{\mathrm{class}}_S(v)=(d\,\Ad^{\mathrm{class}}_S)_e(v),
\end{equation}
we have
\begin{equation}
    \label{adstructure1}
        \Ad(g)=g\circ \Ad^{\mathrm{class}},\qquad  \ad(\theta)=\theta\circ \ad^{\mathrm{class}}
\end{equation}
for all $g$ $\in$ $G$ and all $\theta$ $\in$ $\mathfrak{h}^\#_\C$. 

Finally, if consider the {\it classical} adjoint Lie algebra representation of  $\mathfrak{h}_\C$
\begin{equation}
    \label{adstructure2}
    \mathrm{c}: \mathfrak{h}_\C\times \mathfrak{h}_\C\longrightarrow \mathfrak{h}_\C, \qquad (u,v)\longmapsto [v,u]_\C
\end{equation}
with $[-,-]_\C$ the complexified Lie bracket, then
\begin{equation}
    \label{adstructure4}
    \mathrm{c}^T(\theta)=i\, \theta\circ \mathrm{c}
\end{equation}
for all $\theta$ $\in$ $\mathfrak{h}^\#_\C$, where 
\begin{equation}
    \label{adstructure3}
    \mathrm{c}^T:=(\id_{\mathfrak{h}^\#_\C}\otimes \pi)\circ \ad: \mathfrak{h}^\#_\C\longrightarrow \mathfrak{h}^\#_\C\otimes \mathfrak{h}^\#_\C.
\end{equation}
The reader is encouraged to consult  reference \cite{appendix} for more details. It is worth mentioning that reference \cite{appendix} was developed under the Hermitian convention of the differential; so equation (\ref{adstructure4}) in that paper is written as $c^T(\theta)=\theta\circ c$.

The example developed in this subsection explicitly shows that the universal differential envelope $\ast$--calculus $(\Gamma^\wedge,d,\ast)$ is a proper generalization of the algebra $(\Omega^\bullet_\C(H),d,\ast)$ of $\C$--valued differential forms of $H$ in non--commutative geometry.  

In this way, for a given quantum group $\G$ and a bicovariant $\ast$--FODC $(\Gamma,d)$ over a quantum group $\G$, the triplet $$(\Gamma^\wedge,d,\ast)$$ will be interpreted as the $\ast$--algebra of {\it quantum differential forms} of $\G$. In this sense,  the quantum dual Lie algebra 
$$\mathfrak{qg}^\#={\Ker(\epsilon)\over \mathcal{R}}$$ plays the role of the {\it dualization} $\mathfrak{h}^\#_\C$ of $\mathfrak{h}_\C$  and the  $G$--corepresentation $\ad$ plays the role of the {\it dualization} of the right adjoint action $\ad^{\mathrm{class}}$ of $H$ on $\mathfrak{h}_\C$. Moreover, by equation (\ref{adstructure4}) the linear map
\begin{equation}
    \label{trans1}
    \mathrm{c}^T:= (\id_{\mathfrak{qg}^\#}\otimes \pi)\circ \ad : \mathfrak{qg}^\#\longrightarrow \mathfrak{qg}^\#\otimes \mathfrak{qg}^\# 
\end{equation}
will be interpreted as the {\it quantum Lie bracket}.

\subsection{Example: The Group $U(1)$ with Non--Standard Differential Calculus}

Let us consider the Lie group $U(1)$ and its associated quantum group $\mathcal{U}(1)$. As usual, the dense $\ast$--Hopf algebra will be denoted by $G^\infty=(G,\cdot,\mathbbm{1},\Delta,\epsilon,S,\ast,)$ where 
\begin{equation}
    \label{2.3.f1}
    G=\C[z,z^{-1}]
\end{equation}
is the polynomial Laurent algebra with $z^{-1}=z^\ast$; and  
\begin{equation}
    \label{2.3.f2}
    \Delta(z^n)=z^n\otimes z^n,\qquad \epsilon(z^n)=1,\qquad S(z^n)=z^{n\ast}.
\end{equation}
for every $n$ $\in$ $\Z$.

The next step is to find the universal differential envelope $\ast$--calculus of a $\ast$--FODC over $\mathcal{U}(1)$ different to the one showed in Section $2.2$.

\begin{Proposition}
    \label{2.p.2}
    Let 
    \begin{equation}
    \label{2.3.f11}
    \mathcal{R}=( z^2-\mathbbm{1})\,G\;\;\subseteq\;\; \Ker(\epsilon).
\end{equation}
 Then 
\begin{enumerate}
    \item It induces a $1$--dimensional bicovariant $\ast$--FODC $(\Gamma,d)$.
    \item For $n \geq 2$, $\Gamma^{\wedge\,n}=\otimes^n_G \Gamma$ 
\end{enumerate}
\end{Proposition}
\begin{proof}
    Let $a=z^2-\mathbbm{1}$.
    \begin{enumerate}
        \item  Clearly $\mathcal{R}$ is a right $G$--ideal. Notice that $a$ is not invertible in $G$. Furthermore, $S(a)^\ast=a$ and $\Ad(a)=a\otimes \mathbbm{1}$; which implies that $S(\mathcal{R})^\ast\subseteq \mathcal{R}$ and $\Ad(\mathcal{R})= \mathcal{R}\otimes G$. Therefore, the induced $\ast$--FODC $(\Gamma,d)$ by equation (\ref{2.f9.4}) is bicovariant. 

        On the other hand, since 
        \begin{equation}
            \label{2.3.12}
                        \pi(z)=i\,[z-\mathbbm{1}]_{\mathcal{R}}\not=0,\qquad  \qquad \pi(z^{2n+1})=\pi(z), \qquad  \qquad \pi(z^{2n})=0
        \end{equation}
        for all $n$ $\in$ $\Z$, we obtain 
\begin{equation}
        \label{2.3.f13}
            \mathfrak{qu}(1)^\#:=\mathfrak{qg}^\#=\mathrm{span}_\C\{ \pi(z)\}
\end{equation}
        and hence,  $(\Gamma,d)$ is $1$--dimensional.
        \item Let $ag$ $\in$ $\mathcal{R}$, with $g$ $\in$ $G$. Then, by equation (\ref{2.f16}) we have
        \begin{eqnarray*}
            \pi((ag)^{(1)})\otimes \pi((ag)^{(2)})
            &=&
            \pi(z^2g^{(1)})\otimes \pi(z^2g^{(2)})-\pi(g^{(1)})\otimes \pi(g^{(2)})
            \\
            &=&
            (\pi(z^2)\diamondsuit g^{(1)}+\epsilon(z^2)\pi(g^{(1)}))\otimes (\pi(z^2)\diamondsuit g^{(2)}+\epsilon(z^2)\pi(g^{(2)}))
            \\
            &-&
            \pi(g^{(1)})\otimes \pi(g^{(2)})
            \\
            &=&
             \pi(g^{(1)})\otimes\pi(g^{(2)})-
            \pi(g^{(1)})\otimes \pi(g^{(2)})\,=\,0.
        \end{eqnarray*}
  Thus (see equation (\ref{2.f12.1})) $$\mathfrak{qu}(1)^{\#\wedge \,n}=\otimes^n \mathfrak{qu}(1)^\#/A^{\wedge\,n}=\otimes^n\, \mathfrak{qu}(1)^\#$$ and therefore $\Gamma^{\wedge\,n}=\otimes^n_G \Gamma$ (see equations (\ref{2.f9.4}), (\ref{2.f12.3})).
    \end{enumerate}
\end{proof}

Since $\Ad(z)=z\otimes \mathbbm{1}$,  we get (see equation (\ref{2.f15}))
\begin{equation}
    \label{2.3.f16}
    \ad(\pi(z))=\pi(z)\otimes \mathbbm{1}.
\end{equation}
This implies that the $\ad$ corepresentation is the trivial corepresentation on the corresponding quantum dual Lie algebra $\mathfrak{qu}(1)^\#$. In addition, it is clear that the inner product $\langle-|-\rangle$  defined by $\langle \pi(z) \mid \pi(z)\rangle=1$  makes $\ad$ unitary. Furthermore, by equation (\ref{2.3.f16}) and the fact that $\pi(\mathbbm{1})=0$, it follows that the quantum Lie bracket is identically zero (see equation \ref{trans1}), i.e., 
\begin{equation}
    \label{2.3.16.1}
    \mathrm{c}^T=0.
\end{equation}

Finally, set 
\begin{equation}
    \label{2.3.f17}
    \beta_{U(1)}:=\{\vartheta:=\pi(z)\}
\end{equation}
and in accordance with equations (\ref{properties}), (\ref{2.f12.2}) we obtain
\begin{equation}
    \label{2.3.f18}
    \vartheta^\ast=\vartheta,\qquad d\vartheta=-\vartheta\,\vartheta.
\end{equation}

In summary, taking into account Section $2.2$,  we have presented two different bicovariant finite--dimensional $\ast$--FODC's for the Lie group $U(1)$ characterized by their corresponding quantum dual Lie algebras:
$$\mathfrak{u}(1)^\#_\C,\qquad\qquad \mathfrak{qu}(1)^\#.$$

\section{Quantum Principal Bundles}

This work is develop in M. Durdevich's framework of quantum principal bundles and associated quantum vector bundles. There are mathematical reasons to consider this formulation instead of the Brzezi\'nski–Majid formulation shown in \cite{libro}. The \emph{principal reason concerns all aspects related to curvature}, and the interested reader can find a discussion of this point in Appendix~A.  For further details on the theory of this section, see \cite{micho1, micho2, micho3,stheve}, remembering that for this paper $d(a^\ast)=-(da^\ast)$.

\subsection{Basic Topics}

Let $(B,\cdot,\mathbbm{1},\ast)$ be a quantum space and let $\G$ be a quantum group. A {\it quantum principal $\G$--bundle} over $B$ (abbreviated ``qpb") is a quantum structure formally represented by the triple 
\begin{equation}
\label{2.f17}
\zeta=(P,B,\Delta_P),\qquad \Delta_P:P \longrightarrow P\otimes G
\end{equation}
where $(P,\cdot,\mathbbm{1},\ast)$ is called the {\it quantum total space},  $(B,\cdot,\mathbbm{1},\ast)$ is a $\ast$--subalgebra, which receives the name {\it quantum base space} and 
\begin{enumerate}
\item $\Delta_P$ is a $\ast$--algebra morphism and a $\G$--corepresentation.
\item $\Delta_P(x)=x\otimes \mathbbm{1}$ if and only if $x$ $\in$ $B$.
\item The linear map $\beta:P\otimes P\longrightarrow P\otimes G$ given by $$\beta(x\otimes y):=x\cdot \Delta_P(y):=(x\otimes \mathbbm{1})\cdot \Delta_P(y) $$ is surjective. 
\end{enumerate}

The previous definition is the \emph{standard} one for quantum principal bundles in Durdevich's formulation, with the only difference that, in general, there is no need to work with a quantum group; a Hopf algebra is sufficient \cite{micho2,stheve}.

Given a qpb $\zeta$, a {\it differential calculus} on it is:
 \begin{enumerate}
 \item A graded differential $\ast$--algebra $(\Omega^\bullet(P),d,\ast)$ generated by $\Omega^0(P)=P$ ({\it quantum differential forms of $P$}).
 \item  A bicovariant $\ast$--FODC $(\Gamma,d)$ over $\G$.
 \item The map $\Delta_P$ is extendible to a graded differential $\ast$--algebra morphism $$\Delta_{\Omega^\bullet(P)}:\Omega^\bullet(P)\longrightarrow \Omega^\bullet(P)\otimes \Gamma^{\wedge},$$ where $(\Gamma^\wedge,d,\ast)$ is the universal differential envelope $\ast$--calculus. Here, we have considered that $\otimes$ is the tensor product of graded differential $\ast$--algebras.
 \end{enumerate}

Notice that if $\Delta_{\Omega^\bullet(P)}$ exists, then it is unique because all our graded differential $\ast$--algebras are generated by their degree $0$ elements. Furthermore,  $\Delta_{\Omega^\bullet(P)}$ is a graded differential $\Gamma^\wedge$--corepresentation on $\Omega^\bullet(P)$ (\cite{micho2}). In this way, the space of horizontal forms is defined as 
\begin{equation}
\label{2.f18}
\Hor^\bullet P\,:=\{\varphi \in \Omega^\bullet(P)\mid \Delta_{\Omega^\bullet(P)}(\varphi)\, \in \, \Omega^\bullet(P)\otimes G \},
\end{equation}
and it is a graded  $\ast$--subalgebra of $\Omega^\bullet(P)$ (\cite{stheve}). Since $\Delta_{\Omega^\bullet(P)}(\Hor^\bullet P)\subseteq \Hor^\bullet P\otimes G,$ the map 
\begin{equation}
\label{2.f19}
\Delta_\Hor:=\Delta_{\Omega^\bullet(P)}|_{\Hor^\bullet P}: \Hor^\bullet P \longrightarrow \Hor^\bullet P\otimes G
\end{equation}
is a $\G$--corepresentation on $\Hor^\bullet P$. Also, one can define the space of {\it base} forms ({\it quantum differential forms of $B$}) as 
\begin{equation}
\label{2.f20}
\Omega^\bullet(B):=\{\mu \in \Omega^\bullet(P)\mid \Delta_{\Omega^\bullet(P)}(\mu)=\mu\otimes \mathbbm{1}\}.
\end{equation}
The space of base forms is a graded differential $\ast$--subalgebra of $(\Omega^\bullet(P),d,\ast)$ and in general, it is not generated by $\Omega^0(B)=B$. An explicit example of this fact can be found in reference \cite{appendix}. Finally, the space of vertical forms is defined as (see equation (\ref{2.f12.1}))
\begin{equation}
    \label{vert}
    \Vert^\bullet \,P:=P\otimes \mathfrak{qg}^{\#\wedge}. 
\end{equation}
As the reader can verify in Lemma 3.1 of reference \cite{micho2}, the space $\Vert^\bullet \,P$ can be equipped with a graded differential $\ast$--algebra structure.  The reader can find a proof of the following proposition in  Lemma 3.7 of \cite{micho2}.

\begin{Proposition}
    \label{seq}
    The following sequence of $\ast$--$P$--bimodules
\begin{equation}
\label{3.f1.4}
0\longrightarrow  \Hor^1 P \lhook\joinrel\relbar\joinrel\rightarrow  \Omega^1(P) \xlongrightarrow{\pi_\V} \Vert^1 P \longrightarrow 0
\end{equation}
is exact, where $$\pi_V: \Omega^1(P)\longrightarrow \Vert^1 P,\qquad a\,db\longmapsto a^{(0)}b^{(0)}\otimes S(a^{(1)}b^{(1)})\,a^{(2)}db^{(2)}$$ is a surjective graded differential $\ast$--algebra morphism. 
\end{Proposition}

By {\it dualizing} the notion of principal connections in differential geometry (\cite{nodg}), we have 
\begin{Definition}
    \label{qpc's}
    Let $\zeta$ be a qpb with a differential calculus. A quantum principal connection (abbreviated ``qpc") on $\zeta$ is a linear map
    $$\omega:\mathfrak{qg}^\#\longrightarrow \Omega^{1}(P)$$ such that
    \begin{equation}
    \label{qpc1}
    \Delta_{\Omega^\bullet(P)}(\omega(\theta))=(\omega\otimes \id_G)\ad(\theta)+\mathbbm{1}\otimes\theta, \qquad \omega(\theta^\ast)=\omega(\theta)^\ast,
\end{equation}
for all $\theta$ $\in$ $\mathfrak{qg}^\#$, where $\ad$ is the $\G$--corepresentation given in equation (\ref{2.f15}). Equivalently (\cite{micho2}), a quantum principal connection is a linear map $$\omega:\mathfrak{qg}^\#\longrightarrow \Omega^1(P)$$ such that
\begin{equation}
    \label{qpc3}
    \omega(\theta)\,\in\, \Mor(\ad,\Delta_{\Omega^\bullet(P)}),\quad  (\pi_\V\circ \omega)(\theta)=\mathbbm{1}\otimes \theta,\quad  \omega(\theta^\ast)=\omega(\theta)^\ast
\end{equation}
for all $\theta$ $\in$ $\mathfrak{qg}^\#$, where $\Mor(\ad,\Delta_{\Omega^\bullet(P)})=\{\psi:\mathfrak{qg}^\# \longrightarrow \Omega^\bullet(P) \mid \psi \mbox{ is linear and}$ $(\psi\otimes \id_{G})\circ \ad=(\id_{\Omega^\bullet(P)}\otimes \rho_0)\circ \Delta_{\Omega^\bullet(P)} \circ \psi\}$, with $\rho_0:\Gamma^\wedge\longrightarrow G$ the canonical projection map. 
\end{Definition}
In light of Theorem~12.8 of reference \cite{stheve}, every qpc $\omega$ can be regarded as a left $P$--module monomorphism   from $\mathrm{Ver}^1\,P$ to $\Omega^1(P)$ such that sequence~(\ref{3.f1.4}) splits, in complete analogy with the \emph{classical} case.

In analogy with the {\it classical} case, it can be proved that the set 
\begin{equation}
\label{2.f24.1}
\mathfrak{qpc}(\zeta):=\{\omega:\mathfrak{qg}^\#\longrightarrow \Omega^1(P)\mid \omega \mbox{ is a qpc on }\zeta \}
\end{equation}
is not empty, and it is an affine space modeled by the $\R$--vector space of {\it connection displacements} (\cite{micho2,stheve}) 
\begin{equation}
\label{2.f24.2}
\begin{aligned}
    \overrightarrow{\mathfrak{qpc}(\zeta)}:=\{\lambda:\mathfrak{qg}^\# \longrightarrow \Hor^1 P &\mid \lambda \mbox{ is linear and } \\
    &(\lambda\otimes \id_G)\circ \ad=\Delta_\Hor \circ \lambda,\;\; \lambda\circ \ast=\ast\circ \lambda \}.
\end{aligned}
\end{equation}

A qpc is called {\it regular} if for all $\varphi$ $\in$ $\Hor^{k}P$ and $\theta$ $\in$ $\mathfrak{qg}^\#$, we have 
\begin{equation}
\label{2.f25}
\omega(\theta)\,\varphi=(-1)^{k}\varphi^{(0)}\omega(\theta\diamondsuit\varphi^{(1)}), 
\end{equation}
where $\Delta_\Hor(\varphi)=\varphi^{(0)}\otimes\varphi^{(1)}.$ A qpc $\omega$ is called {\it multiplicative} if 
\begin{equation}
\label{2.f26}
\omega(\pi(g^{(1)}))\omega(\pi(g^{(2)}))=0
\end{equation}
for all $g$ $\in$ $\mathcal{R}$. 

It should be noted that every qpc that comes from the {\it dualization} of a {\it classical} principal connection is regular and multiplicative (\cite{micho1}). In fact, let  $$\zeta=(P,B,\Delta_P)$$  be the  {\it dualization} of a {\it classical} principal $H$--bundle $$\pi_X: Y\longrightarrow X,$$ with $H$ a (compact matrix) Lie group, $Y$ the total (compact) space, $X$ the base (compact) space and $\pi_X$ the bundle--projection (\cite{micho1}). Let us take the differential calculus on $\zeta$ given by considering the space of $\C$--valued differential forms of $Y$ and the bicovariant $\ast$--FODC showed in Section $2.2$. Let $\omega$ be a qpc given by the pull--back of a principal connection of $\pi_X: Y\longrightarrow X$.  By equation (\ref{2.2.f4}) and the fact that the space of $\C$--valued differential forms of $Y$ is graded--commutative, it follows that $\omega$ is regular. Furthermore, a simple calculation using equations (\ref{2.f16}), (\ref{2.2.f4}) proves that $$\omega(\pi(g^{(1)}))\omega(\pi(g^{(2)}))=\omega(\pi(b))\omega(\pi(a))+\omega(\pi(a))\omega(\pi(b))=0$$ if $g=a\,b$ $\in$ $\Ker^2(\epsilon)$. By linearity we conclude that $\omega$ is also multiplicative.

By {\it dualizing} the notion of the covariant derivative of a principal connection in differential geometry (\cite{nodg}),  for a given qpc $\omega$, we define its {\it covariant derivative} as the first--order linear map (\cite{micho3})
\begin{equation}
\label{2.f30}
D^{\omega}: \Hor^\bullet P \longrightarrow \Hor^\bullet P
\end{equation}
such that for every $\varphi$ $\in$ $\Hor^k P$ $$
D^{\omega}(\varphi)=  d\varphi-(-1)^{k}\varphi^{(0)}\omega(\pi(\varphi^{(1)})).$$ 
On the other hand, the first--order linear map 
\begin{equation}
\label{2.f30.1}
\widehat{D}^{\omega}:=-\ast\circ D^{\omega}\circ\ast: \Hor^\bullet P \longrightarrow \Hor^\bullet P
\end{equation}
is called the {\it dual covariant derivative} of $\omega$. Explicitly, we have (\cite{micho3}) 
\begin{equation}
    \label{2.f30.5.1}
    \widehat{D}^{\omega}(\varphi)=d\varphi +\omega (\pi(S^{-1}(\varphi^{(1)})))\varphi^{(0)}
\end{equation}
for every $\varphi$ $\in$ $\Hor^k P$.  According to \cite{micho3}, only for regular qpc's can we guaranty that $D^\omega=\widehat{D}^{\omega}$,
which is the situation for qpc's arising from {\it classical} principal connections (\cite{micho1}). In other words, $D^\omega$ and $\widehat{D}^\omega$ are two different horizontal operators that generalize the covariant derivative of a principal connection in differential geometry. For more details, see reference \cite{micho3}.

According to references \cite{micho2,stheve} the covariant derivatives satisfy $$D^\omega,\;\widehat{D}^\omega\;\in \Mor(\Delta_\Hor,\Delta_\Hor)$$ and hence, for every finite--dimensional $\G$--corepresentation, we can consider 
\begin{equation}
\label{new1}
        D^\omega:\Mor(\delta^V,\Delta_\H)\longrightarrow \Mor(\delta^V,\Delta_\H),\qquad \tau  \longmapsto D^\omega(\tau),
    \end{equation}
\begin{equation}
\label{new2}
        \widehat{D}^\omega:\Mor(\ad,\Delta_\H)\longrightarrow \Mor(\ad,\Delta_\H),\qquad \tau  \longmapsto \widehat{D}^\omega(\tau)
    \end{equation}
given by  $D^\omega(\tau)(v)=D^\omega(\tau(v))$  and   $\widehat{D}^\omega(\tau)(v)=\widehat{D}^\omega(\tau(v))=(D^\omega(\tau(v)^\ast))^\ast$ 
for all $v$ $\in$ $V$.

In Durdevich's formulation, to define the curvature of a qpc, it is necessary the following auxiliary map.

\begin{Definition}
    \label{embbededdifferential}
    We define an {\it embedded differential} as a linear map $$\Theta:\mathfrak{qg}^\#\longrightarrow \mathfrak{qg}^\#\otimes \mathfrak{qg}^\# $$ such that
\begin{enumerate}
    \item  $\Theta$ $\in$ $\Mor(\ad,\ad^{\otimes 2})$, where $\ad^{\otimes 2}$ is the tensor product $\G$--corepresentation of $\ad$ with itself (\cite{woro1}). 
    \item If $\Theta(\theta)=\displaystyle\sum^m_{i,j=1}\theta_i\otimes \theta'_j$, then $d\theta=\displaystyle\sum^m_{i,j=1}\theta_i \theta'_j$ and $\Theta(\theta^\ast)=\displaystyle\sum^m_{i,j=1}\theta'^\ast_j \otimes \theta^\ast_i$.
\end{enumerate}
\end{Definition}
\noindent In light of \cite{micho1}, embedded differentials always exist. Furthermore, the Maurer--Cartan equation shows that in the dualization of the \emph{classical} case, $\displaystyle -{1\over 2}c^T$ is an embedded differential (see equation (\ref{trans1})).

 We define
\begin{equation}
\label{trans}
[ T_1,T_2]:=m_X \circ (T_1\otimes T_2)\circ \mathrm{c}^T:\mathfrak{qg}^\# \longrightarrow X,
\end{equation}
\begin{equation}
    \label{embed}
    \langle T_1,T_2\rangle:=m_X\circ(T_1\otimes T_2)\circ\Theta: \mathfrak{qg}^\#\longrightarrow X
\end{equation}
for every map $T_1,\,\,T_2:\mathfrak{qg}^\# \longrightarrow X$, where $X$ any $\ast$--algebra and  $m_X$ is the product map of $X$.  

Fix $\Theta$ an embedded differential. Hence, according to \cite{micho2,stheve}, we define the curvature of a quantum principal connection $\omega$ as the linear map
\begin{equation}
\label{2.f28}
R^{\omega}:=d\omega-\langle \omega,\omega\rangle : \mathfrak{qg}^\#\longrightarrow  \Omega^2(P)
\end{equation}
and in Theorem 12.11 of reference \cite{stheve} is proven that
\begin{equation}
\label{2.f29}
R^{\omega}\,\in\, \Mor(\ad,\Delta_\Hor).
\end{equation}
for every qpc $\omega$. Notice that the map $\Theta$ is necessary in order for $R^\omega$ to have domain $\mathfrak{qg}^\#$ and for equation~(\ref{2.f29}) to hold. In particular, this allows one to interpret $R^{\omega}$ as the \emph{quantum field strength} in physical terms. We encourage the reader to check references \cite{micho1,micho2,saldym} for a more detail explanation about the definition of the curvature in Durdevich's formulation of qpb's.

 According to references \cite{micho2, stheve}, if $\omega$ is multiplicative, then $R^\omega$ does not depend on the choice of $\Theta$. This is the main reason to study multiplicative qpc's.  In particular, this is the case of every qpc that comes from the {\it dualization} of a {\it classical} principal connection, where the {\it standard} embedded differential used is $\displaystyle -{1\over 2}\mathrm{c}^T$ (\cite{woro1}).  

 The presence of $\Theta$ in the definition of $R^\omega$ in the {\it non--commutative geometrical} case can be interpreted as a {\it quantum phenomenon} in which there can be several non--equivalent ways to embed $\mathfrak{qg}^\#$ into $\mathfrak{qg}^\#\otimes \mathfrak{qg}^\#$ and produce horizontal quadratic expressions with $\omega$.

It is easy to check that the following map is well--defined
\begin{equation}
    \label{ad.f3}
        \wedge: \Mor(\ad,\Delta_{\Omega^\bullet(P)})\longrightarrow \Mor(\ad,\Delta_{\Omega^\bullet(P)}),\qquad 
        \psi  \longmapsto \widehat{\psi}:=\ast \circ \psi \circ \ast.
\end{equation}
In addition,  in light of Proposition 12.13 of reference \cite{stheve}, we get
\begin{equation}
    \label{inersperado}
    \widehat{\langle \psi,\phi\rangle}=\ast\circ \langle\psi,\phi\rangle\circ \ast=(-1)^{kl}\,\langle\widehat{\phi},\widehat{\psi}\rangle,
\end{equation}
where $\widehat{\psi}=\ast\circ \psi\circ \ast,$ $ \widehat{\phi}=\ast\circ \phi\circ \ast$  with   $\psi,\; \phi\; \in\; \Mor(\ad,\Delta_{\Omega^\bullet(P)})$   such that  $\Im(\psi)\subseteq \Omega^k(P),\;  \Im(\phi)\subseteq \Omega^l(P).$  It is worth mentioning that equation~(\ref{inersperado}) differs by a sign from the result of Proposition~12.13 in \cite{stheve}. This discrepancy arises because, in our case, we work with the anti--Hermitian convention for the differential.

In particular, this implies that
\begin{equation}
\label{curast}
    R^{\omega}=- \widehat{R}^\omega.
\end{equation}
Moreover, for every $\tau$ $\in$ $\Mor(\ad,\Delta_\Hor)$,  we obtain
   $\widehat{D}^\omega(\tau)=(\wedge \circ D^\omega\circ \wedge)(\tau).$

\begin{Definition}
    \label{Soperator}
    Let $\omega$ be a qpc. We define the operator 
    $$S^{\omega}: \Mor(\ad,\Delta_\H) \longrightarrow \Mor(\ad,\Delta_\H)$$
   given by
    \begin{equation*}
S^{\omega}(\tau):=\langle \omega,\tau\rangle-(-1)^k\langle\tau,\omega\rangle-(-1)^k[\tau,\omega]
\end{equation*}
for every $\tau$ $\in$ $\Mor(\ad,\Delta_\H)$ with $\Im(\tau)$ $\subseteq$ $\Hor^k P$. Similarly, we define the dual $S^\omega$ operator as $$\widehat{S}^\omega :=\wedge\circ S^\omega \circ \wedge.$$ 
\end{Definition}
\noindent The reader is encouraged to consult this reference for more details on the operator $S^{\omega}$. For example, $S^{\omega}=0$ when $\omega$ is regular. Using Definition \ref{Soperator} and Proposition 4.7 of reference \cite{micho2}, we obtain

\begin{Definition}
    \label{twisted}
    We define the twisted covariant derivative of a qpc $\omega$ as the operator $$DS^\omega:=D^\omega-S^\omega: \Mor(\ad,\Delta_\H) \longrightarrow \Mor(\ad,\Delta_\H).$$ Explicitly, for every $\tau$ $\in$ $\Mor(\ad,\Delta_\H)$ with $\Im(\tau)$ $\subseteq$ $\Hor^k P$ we have 
    \begin{equation}
        \label{twcoder}
        DS^\omega(\tau)=d\tau-\,\langle \omega,\tau\rangle+(-1)^k\,\langle\tau,\omega\rangle.
    \end{equation}
    In the same way, we define the dual twisted covariant derivative as the operator $$\widehat{DS}^\omega=\wedge \circ DS^\omega\circ \wedge.$$
\end{Definition}

\noindent For example, for every $\lambda$ $\in$ $\overrightarrow{\mathfrak{qpc}(\zeta)}$, we get (see equation (\ref{2.f24.2}))
    \begin{equation}
        \label{twcoder1}
        \widehat{DS}^\omega(\lambda)=-DS^\omega(\lambda)\qquad \mbox{ and }\qquad \widehat{DS}^\omega(R^\omega)=DS^\omega(R^\omega).
    \end{equation}

\begin{Remark}
\label{rema}
From this point onward until the end of the paper, we shall restrict our attention exclusively to qpb's for which the quantum base space $(B,\cdot,\mathbbm{1},\ast)$ is stable under functional calculus \cite{con}. According to Appendix B of reference \cite{micho3}, in this case, for every $\delta^V$ $\in$ $\T$ there exists a set  

$$\{T^\l_k \}^{d_{V}}_{k=1} \subseteq \Mor(\delta^V,\Delta_P)$$ for some $d_{V}$ $\in$ $\N$ such that
\begin{equation}
    \label{generators}
\sum^{d_{V}}_{k=1}x^{V\,\ast}_{ki}x^{V}_{kj}=\delta_{ij}\mathbbm{1},
\end{equation}
with $x^{V}_{ki}:=T^\l_k(e_i)$. Here, $\T$  is a complete set of mutually non--equivalent irreducible (necessarily finite--dimensional) $\G$--corepresentations with $\delta^\C_\triv$ $\in$ $\T$ (the trivial corepresentation on $\C$), and $\{e_i\}^{n_{V}}_{i=1}$ is the orthonormal basis of $V$ given in Theorem \ref{rep}.
\end{Remark}

It is worth mentioning that in Proposition 2.7 of reference \cite{sald2}, we show the specific form of the operators $\{T^\l_k\}_{k=1}^{d_V}$ in differential geometry. Chapter 6 of reference \cite{con} also works with $\ast$--algebras stable under holomorphic functional calculus because in this situation, in  A. Connes' words, {\it all possible notions of positivity will coincide}.

\subsection{Example: Quantum Principal Bundles over Toric Non--Commutative Manifolds}

This subsection will be entirely based on references \cite{lan4,lan2}. Consider the {\it deformation functor} presented in \cite{lan4,lan2} $$L_\Xi: {\mathcal V}\longrightarrow {\mathcal V}_{\Xi}.$$ Here, ${\mathcal V}$ is the category of nuclear Fr\'echet spaces whose topology is determined by a countable family of semi--norms equipped with a smooth torus--action. On the other hand, the category ${\mathcal V}_{\Xi}$ is essentially the same as ${\mathcal V}$, but with a different braided monoidal structure $\widetilde{\otimes}_{\Xi}$ with respect to the canonical one $\widetilde{\otimes}$ of ${\mathcal V}$, where $\widetilde{\otimes}$ is the corresponding completion of the algebraic tensor product $\otimes$. The new braided monoidal structure $\widetilde{\otimes}_{\Xi}$ is defined by using the tours--action and some fixed real antisymmetric $n\times n$ matrix $\Xi$. For more details, see \cite{lan4,lan2}.

Moreover, in accordance with \cite{lan4,lan2}, by taking a $\ast$--algebra $\mathcal{A}$ in the category $\mathcal{V}$ (in particular, the product is commutative, continuous and torus--equivariant), its {\it image} under $L_\Xi$ is a $\ast$--algebra $\mathcal{A}_\Xi$ in the category $\mathcal{V}_\Xi$. In particular, its product is non--commutative, continuous and tours--equivariant \cite{lan4,lan2}.

Let $$\pi_X:Y\longrightarrow X $$
be a {\it classical} principal $H$--bundle. Assume that $X$ and $Y$ are closed, oriented, spin Riemannian manifolds. Furthermore, suppose that $X$ is an $n$--dimensional manifold equipped with an isometric smooth action of an $\widetilde{n}$--torus ($\widetilde{n}\geq 2$) such that this action can be lifted to an isometrical action on $Y$ via a covering torus. Also, assume that this lifted action commutes with the right action of $H$ on $Y$ that provides the structure of a principal $H$--bundle. 

Let $\Omega^\bullet_\C(X)$, $\Omega^\bullet_\C(Y)$, $\Omega^\bullet_\C(H)$ be the spaces of $\C$--valued differential forms of $X$, $Y$ and $H$, respectively. Consider now the graded differential $\ast$--algebras 
\begin{equation}
    \label{ec.012}
(\Omega^\bullet(X_\Xi)=L_\Xi(\Omega^\bullet_\C(X)),d:=,\ast),\quad (\Omega^\bullet(Y_\Xi)=L_\Xi(\Omega^\bullet_\C(Y)),d:=,\ast),
\end{equation}
where $d=id'$ with $d'$ the de Rham differential, and we have extended the deformation functor $L_\Xi$ to the level of differential forms by imposing that it commutes with $d'$.

The spaces $X_\Xi:=\Omega^0(X_\Xi)=L_\Xi(C^\infty_\C(X))$, $Y_\Xi:=\Omega^0(Y_\Xi)=L_\Xi(C^\infty_\C(Y))$ receive the name of {\it toric non--commutative manifolds}.  The  product of $X_\Xi$ can be defined as follows: every element $x$ $\in$ $C^\infty_\C(X)$ is the sum of a unique rapidly convergent functions with respect to the spectral decomposition of the action of the $\widetilde{n}$--torus $$x=\displaystyle \sum_{r\in \Z^{\widetilde{n}}}x_r;$$ therefore, given $x$, $x'$ $\in$ $X_\Xi$, its product is given by 
\begin{equation}
    \label{3.2.0}
    x\cdot x':=\sum_{r,\,r'\,\in\, \Z^{\widetilde{n}}} \mathrm{e}^{2\pi i \,r\cdot \Xi\cdot r'}\,x_r\,x'_{r'},
\end{equation}
where $$r\cdot \Xi\cdot r'=(r_1,\cdots r_{\widetilde{n}})\cdot \Xi \cdot \left(\begin{array}{lcr}
r'_1\\
\;\vdots\\
r'_{\widetilde{n}}\\
\end{array}\right)=\sum^{\widetilde{n}}_{k,j=1} r_k\,\Xi_{kj}\,r'_j.$$  Since $\Xi$ is antisymmetric, it follows that $r\cdot \Xi\cdot r'=-r'\cdot \Xi\cdot r$ and hence
\begin{equation}
    \label{ec.3.29}
    x'\cdot x=\sum_{r,r'\in \Z^{\widetilde{n}}} \mathrm{e}^{-2\pi i \,r\cdot \Xi\cdot r'}\,x'_{r'}\,x_r.
\end{equation}
In particular, by defining $$\hat{x}'=\displaystyle \sum_{r\in \Z^{\widetilde{n}}} \mathrm{e}^{4\pi i \,r\cdot \Xi\cdot r'}\,x'_r, \qquad \widetilde{x}'=\displaystyle \sum_{r\in \Z^{\widetilde{n}}} \mathrm{e}^{-4\pi i \,r\cdot \Xi\cdot r'}\,x'_r $$ we get
\begin{equation}
    \label{ec.3.30}
   x\cdot x'= \hat{x}'\cdot x, \qquad x'\cdot x=x\cdot \widetilde{x}'.
\end{equation}
The last equation provides a commutation--relation for the element $x$ acting on the left and on the right. This new product extends to the level of differential forms $\Omega^\bullet(X_\Xi)$ and of course, equation (\ref{ec.3.30}) extends to $\Omega^\bullet(X_\Xi)$: for elements $\mu$,  $\in$ $\Omega^k(X_\Xi)$, $\mu'$ $\in$ $\Omega^l(X_\Xi)$, we have
\begin{equation}
    \label{ec.3.30.1}
    \mu\cdot \mu'=(-1)^{kl}\hat{\mu}'\cdot \mu,\qquad \mu'\cdot \mu=(-1)^{kl}\mu\cdot \widetilde{\mu}', 
\end{equation}
for some $\hat{\mu}'$, $\widetilde{\mu}'$ $\in$ $\Omega^l(X_\Xi)$. The same product and relations hold  for $Y_\Xi$ and $\Omega^\bullet(Y_\Xi)$.

In accordance with \cite{lan4,lan2},  the right $H$--action $$r:Y\times H\longrightarrow Y$$ defines by its pull--back $r^\#$ on differential forms, a graded differential $\ast$--algebra morphism 
\begin{equation}
    \label{3.2.3}
    \Delta_{\Omega^\bullet(Y_\Xi)}: \Omega^\bullet(Y_\Xi)\longrightarrow \Omega^\bullet(Y_\Xi) \,\widetilde{\otimes}\, \Omega^\bullet_\C(H) 
\end{equation}
in a way that
\begin{equation}
    \label{3.2.1}
    \zeta_\Xi=(Y_\Xi,X_\Xi,\Delta_{Y_\Xi}) \qquad \mbox{ where }\qquad \Delta_{Y_\Xi}:=\Delta_{\Omega^\bullet(Y_\Xi)}|_{Y_\Xi}
\end{equation}
is a qpb over $X_\Xi$ with a differential calculus \cite{lan4,lan2}. It is worth mentioning that the previous tensor product $\widetilde{\otimes}$ is the corresponding completion of the algebraic tensor product. Although this lies beyond the purely
geometric--algebraic approach showed in the previous subsection, according to references \cite{lan4,lan2}, {\it mutatis mutandis}, everything still works.

Since the space of horizontal differential forms of the bundle $\pi_X$ is given by (\cite{nodg}) $$C^\infty_\C(Y)\,\Omega^\bullet_\C(X),$$ by applying the deformation functor we get (\cite{lan4})
\begin{equation}
    \label{ec.3.31}
    \Hor^\bullet \,Y_\Xi=Y_\Xi\cdot\Omega^\bullet(X_\Xi).
\end{equation}

Notice that we can consider $\Omega^\bullet(X_\Xi)$ as a graded differential $\ast$--subalgebra of $\Omega^\bullet(Y_\Xi)$ since the pull--back $\pi^\#_X$ on differential forms of the bundle projection $\pi_X$ induces a $\ast$--algebra monomorphism (\cite{lan4,lan2})    
\begin{equation}
    \label{3.2.2}
    \iota: \Omega^\bullet(X_\Xi)\longrightarrow \Omega^\bullet(Y_\Xi).
\end{equation}
In addition, it is clear that $X_\Xi$ is stable under functional calculus \cite{lan4,lan2}.

It is worth mentioning that we have decided to use $d=i\,d'$ in order to preserve the anti--Hermitian convention of the differential and of course, explicit examples of this type of qpb's are presented in \cite{con,lan4}.

According to Section 4 of reference \cite{micho1}, the following isomorphism of graded $\ast$--algebras holds (here, the tensor product is the tensor product of graded vector spaces) $$\Omega^\bullet(Y)\cong \Hor^\bullet\,Y\otimes \mathfrak{h}^{\#\wedge}_\C,$$ where the graded $\ast$--algebra structure of $\Hor^\bullet\,Y\otimes \mathfrak{h}^{\#\wedge}_\C$ is given by
\begin{equation}
    \label{2.fver2}
    \begin{aligned}
        (\psi\otimes\theta)(\phi\otimes \vartheta)&:=(-1)^{kl}\psi \phi^{(0)}\otimes (\theta\diamondsuit \phi^{(1)})\vartheta,\\
(\psi\otimes\theta)^{\ast}&:=\psi^{(0)\ast}\otimes(\theta^{\ast}\diamondsuit \psi^{(1)\ast}),
    \end{aligned}
\end{equation}
where $\psi$  $\in$ $\Hor^\bullet \,Y$, $\phi$  $\in$ $\Hor^k \,Y$, $\vartheta$ $\in$ $\mathfrak{h}^{\#\wedge l}_\C$, $\theta$ $\in$ $\mathfrak{h}^{\#\wedge}_\C$, $\Delta_\Hor(\psi)=\psi^{(0)}\otimes \psi^{(1)}$, $\Delta_\Hor(\phi)=\phi^{(0)}\otimes \phi^{(1)}$ and $\theta \diamondsuit g$  is defined in equations (\ref{2.f16}), (\ref{2.f12.4}). Here, $\Delta_\Hor$ is the pull--back of the right $H$--action on the horizontal forms. In addition, under this isomorphism, $r^\#$ on differential forms is given by $$r^\#(\psi\otimes \theta):=\Delta_\Hor(\psi)\cdot \Delta(\theta):=(\psi^{(0)}\otimes \mathbbm{1})\cdot(\mathbbm{1}\otimes \theta^{(1)}) \otimes \psi^{(1)}\theta^{(2)},$$ with $\Delta_\Hor(\psi)=\psi^{(0)}\otimes \psi^{(1)}$, $\Delta(\theta)=\theta^{(1)}\otimes \theta^{(2)}$. By applying the deformation functor, we obtain  
\begin{equation}
    \label{teodio}
    \Omega^\bullet(Y_\Xi)\cong \Hor^\bullet\,Y_\Xi\otimes \mathfrak{h}^{\#\wedge}_\C \quad \mbox{ and }\quad \Delta_{\Omega^\bullet(Y_\Xi)}=\Delta_\Hor\cdot \Delta.
\end{equation}
\begin{Proposition}
    \label{teodio1}
    The qpb $\zeta_\Xi$ admits a regular qpc.
\end{Proposition}
\begin{proof}
    Consider the linear map $$\omega^c:\mathfrak{h}^\#_\C\longrightarrow \Omega^1(Y_\Xi)$$ such that under the isomorphism of equation (\ref{teodio}), it satisfies $$\omega^c(\theta)=\mathbbm{1}\otimes \theta.$$ Clearly $\omega^c(\theta^\ast)=\omega(\theta)^\ast$ and by equation (\ref{coproduc.1}) $$\Delta_{\Omega^\bullet(Y_\Xi)}(\omega^c(\theta))=\mathbbm{1}\otimes \theta^{(0)}\otimes \theta^{(1)}+\mathbbm{1}\otimes \mathbbm{1}\otimes \theta=(\omega^c\otimes \id_G)\circ \ad(\theta)+\mathbbm{1}\otimes \mathbbm{1}\otimes \theta,$$ with $\ad(\theta)=\theta^{(0)}\otimes \theta^{(1)}$. This shows that $\omega^c$ is a qpb. In addition since in our case, $\theta \diamondsuit g=\epsilon(g)\theta$ (see equation (\ref{2.2.f4})), it follows that a qpc is regular if and only if it graded--commutes with all elements of $\Hor^\bullet\,Y_\Xi$. In the same way, equation (\ref{2.2.f4}) guarantees that $\omega^c$ graded--commutes with all elements of $\Hor^\bullet\,Y_\Xi$ and hence, it is a regular qpc.
\end{proof}

\section{Associated Quantum Vector Bundles and Gauge Quantum Linear Connections}

This section is entirely based on reference \cite{sald2}. The reader is encourage to consult this reference for more details. Let $\zeta=(P,B,\Delta_P)$ be a quantum $\G$--bundle with a differential calculus and a $\G$--corepresentation $\delta^V$ $\in$ $\T$.  According to Section IIC of reference \cite{saldym}, the space 
\begin{equation}
    \label{3.2.f1}
    \Mor(\delta^V,\Delta_\Hor)
\end{equation}
is a finitely generated projective left/right $\Omega^\bullet(B)$--module, where the structure of left/right $\Omega^\bullet(B)$--module is given by left/right multiplication with elements of  $\Omega^\bullet(B)$. Furthermore, the following relation holds (\cite{saldym})
\begin{equation}
\label{3.f5}
\tau=\sum^{d_{V}}_{k=1}\mu^{\tau}_k\,T^\l_k \quad \mbox{ with } \quad \mu^{\tau}_k=\sum^{n_{V}}_{i=1}\tau(e_i)\,x^{V\,\ast}_{ki}\;\in\;\Omega^\bullet(B), \quad T^{\l}_k(e_i)=x^{V}_{ki}.
\end{equation}
In accordance with reference \cite{sald2}, we have a left $\Omega^\bullet(B)$--module isomorphism
$$\Upsilon_V: \Mor(\delta^V,\Delta_\Hor)\longrightarrow \Omega^\bullet(B)\otimes_B \Mor(\delta^V,\Delta_P) $$ given by
\begin{equation}
\label{3.f7}
\Upsilon_{V}(\tau)=\sum^{d_{V}}_{k=1}\mu^\tau_k\otimes_B T^\l_k,
\end{equation}
with inverse $$\Upsilon^{-1}_{V}:\Omega^\bullet(B)\otimes_{B}\Mor(\delta^V,\Delta_P)\longrightarrow\Mor(\delta^V,\Delta_\Hor)$$ given by 
\begin{equation}
\label{3.f5.1.2.3}
\Upsilon^{-1}_{V}(\mu\otimes_{B} T)=\mu\, T.
\end{equation}

Since $$\Mor(\delta^V\bigoplus \delta^W, \Delta_\Hor)=\Mor(\delta^V,\Delta_\Hor)\bigoplus \Mor(\delta^W,\Delta_\Hor)$$ for any $\delta^V$, $\delta^W$ $\in$ $\T$ and every finite--dimensional $\G$--corepresentation $\delta^{V}$ is the direct sum of a finite number of elements of $\T$ (\cite{woro1}),  it follows that $\Mor(\delta^V,\Delta_\Hor)$ is a finitely generated projective left $\Omega^\bullet(B)$--module for every $\delta^{V}$ $\in$ $\Obj(\Rep_{\G})$. Moreover, by applying equation (\ref{3.f5}) to each summand $\delta^{V_j}$ $\in$ $\T$ in the decomposition $\delta^V\cong \bigoplus^m_{j=1} \delta^{V_j}$, this equation naturally extends to every  $\tau$ $\in$ $\Mor(\delta^V,\Delta_\Hor)$.  Similarly, the map $\Upsilon_V$ can be extended to every $\delta^{V}$ $\in$ $\Obj(\Rep_{\G})$. Apparently, $\Upsilon_{V}$ depends on the set of generators $\{ T^\l_k\}^{d_V}_{k=1}$ of $\Mor(\delta^{V},\Delta_P)$. However,  the uniqueness of the inverse map ensures that $\Upsilon_{V}$ is independent of the choice of this set.

Let $\delta^{V}$ $\in$ $\Obj(\Rep_{\G})$. In light of reference \cite{saldym} we have that the space
\begin{equation}
\label{3.f5.1}
\Mor(\delta^{V},\Delta_\Hor)
\end{equation} 
is also a finitely projective right $\Omega^\bullet(B)$--module, where the structure of right $\Omega^\bullet(B)$--module is given by right multiplication with elements of  $\Omega^\bullet(B)$. Furthermore, the following relation holds (\cite{saldym})
\begin{equation}
\label{3.f5.2}
\tau=\displaystyle \sum_k T^\r_k (\mu^{\tau^\ast}_k)^\ast
\end{equation} 
with $T^\r_k:=T^{\l\,\ast}_k$. Here, the maps $\{T^{\l}_k\}$ are corresponding left generators of $\Mor(\delta^{\overline{V}},\Delta_P)$, where $\delta^{\overline{V}}$ is the complex conjugate corepresentation of $\delta^V$ and hence $$\mu^{\tau^\ast}_k=\sum_i \tau(e_i)^\ast\,x^{\overline{V}\,\ast}_{ki} \;\in\;\Omega^\bullet(B),\qquad T^{\l}_k(e_i)=x^{\overline{V}}_{ki}.$$

As before, there exists a right $\Omega^\bullet(B)$--module isomorphism (\cite{sald2})
$$\widehat{\Upsilon}_{V}: \Mor(\delta,\Delta_\Hor)\longrightarrow \Mor(\delta^V,\Delta_P)\otimes_B \Omega^\bullet(B)$$ given by
\begin{equation}
\label{3.f7.5}
  \widehat{\Upsilon}_{V}(\tau)=\sum_k T^\r_k\otimes_B (\mu^{\tau\ast}_k)^\ast 
\end{equation}
with inverse 
$$\widehat{\Upsilon}^{-1}_{V}:  \Mor(\delta^V,\Delta_P)\otimes_B \Omega^\bullet(B)\longrightarrow \Mor(\delta,\Delta_\Hor) $$ given by
\begin{equation}
    \label{3.f5.1.2.3.1}
    \widehat{\Upsilon}^{-1}_{V}(T\otimes_B\mu)=T\,\mu.
\end{equation}
Apparently, the map $\widehat{\Upsilon}_{V}$ depends on the set of generators $\{ T^\l_k\}^{d_V}_{k=1}$ of $\Mor(\delta^{\overline{V}},\Delta_P)$. However,  the uniqueness of the inverse map ensures that $\widehat{\Upsilon}_{V}$ is independent of the choice of this set.

Let $\delta^V$ $\in$ $\Obj(\Rep_\G)$. We define $$E^V_\l:=\Mor(\delta^V,\Delta_P)$$  as left $B$--module,  and we define $$E^V_\r:=\Mor(\delta^V,\Delta_P)$$ as right $B$--module. By our previous constructions and the fact that $\Delta_\Hor=\Delta_P$ in zero--degree, it follows that $E^V_\l$ is a finitely generated projective left $B$--module; while $E^V_\r$ is a finitely generated projective right $B$--module. 

In the {\it classical} case, given a principal $H$--bundle  $\pi_X:Y\longrightarrow X$ and a linear representation of $H$ on $V$, the space of global smooth sections $\Gamma(E^V)$ of the associated vector bundle is isomorphic to the space of $G$--equivariant maps $C^\infty_\C(Y,V)^G$ as $C^\infty_\C(X)$--bimodules (\cite{nodg}).  The {\it non--commutative geometrical} counterpart of $C^\infty_\C(Y,V)^G$ is the space $\Mor(\delta^V,\Delta_P)$ and hence, in  light of reference \cite{con, dv} and the Serre--Swan theorem we have the following definition

\begin{Definition}
\label{qvb's}
    Let $\zeta=(P,B,\Delta_P)$ be a quantum principal $\G$--bundle and let $\delta^V$ $\in$ $\Obj(\Rep_\G)$. We define the {\it associated left quantum vector bundle} (abbreviated ``associated left qvb") as the space $E^V_\l$, and we define the {\it associated right quantum vector bundle} (abbreviated ``associated right qvb") as the space $E^V_\r$.
\end{Definition}

We will write $\Mor(\delta^V, \Delta_P)$ to refer to the $B$--bimodule structure or the $\C$--vector space structure of this set, while the notation $E^V_\l/E^V_\r$ is to refer to the left/right structure of this set. Similarly, we will write $\Mor(\delta^V, \Delta_\H)$ to refer to the $\Omega^\bullet(B)$--bimodule structure or the $\C$--vector space structure of this set; while the notation 
\begin{equation}
    \label{algo}
    \E^V_\l,\qquad \E^V_\r
\end{equation}
is to refer to the left/right  $\Omega^\bullet(B)$--module structure of this set, respectively. Notice that elements of $$\E^V_\l \cong \Omega^\bullet(B)\otimes_{B}E^V_\l \;\; (\mbox{as left } \Omega^\bullet(B)\mbox{--modules}) $$ can be interpreted as {\it left qvb--valued differential forms} of $B$; while elements of $$\E^V_\r \cong E^V_\r\otimes_{B}\Omega^\bullet(B) \;\; (\mbox{as right } \Omega^\bullet(B)\mbox{--modules}) $$ can be interpreted as {\it right qvb--valued differential forms} of $B$.

\begin{Definition}
    \label{iqlcs}
    Let $\omega$ a qpc and let $\delta^V$ $\in$ $\Obj(\Rep_\G)$.  By equation (\ref{new1}) we have that the linear map
\begin{equation}
\label{3.f8}
\nabla^{\omega}_{V}:E^V_\l \longrightarrow \Omega^{1}(B)\otimes_B E^V_\l,\qquad
T  \longmapsto \Upsilon_{V}(D^{\omega}(T)),
\end{equation}
is a {\it quantum linear connection} on $E^V_\l$, in the sense of reference \cite{dv}, i.e., $\nabla^{\omega}_{V}$ satisfies the left Leibniz rule: for every $b$ $\in$ $B$ and every $T$ $\in$ $E^V_\l$ we have $$\nabla^{\omega}_V(b\,T)=db\otimes_B T+b\,\nabla^\omega_V(T).$$  Similarly, by equation (\ref{new2}) we have that the linear map
\begin{equation}
\label{3.f8.1}
\widehat{\nabla}^{\omega}_{V}:E^V_\r \longrightarrow E^V_\r\otimes_{B}\Omega^{1}(B),\qquad
T  \longmapsto \widehat{\Upsilon}_{V}(\widehat{D}^{\omega}(T)),
\end{equation}
is a {\it quantum linear connection} on $E^V_\r$, i.e., $\widehat{\nabla}^{\omega}_{V}$ satisfies the right Leibniz rule: for every $b$ $\in$ $B$ and every $T$ $\in$ $E^V_\r$ we have $$\widehat{\nabla}^{\omega}_{V}(T\,b)=\widehat{\nabla}^{\omega}_{V}(T)\,b+T\otimes_B db .$$ The maps  $\nabla^{\omega}_{V}$ and $\widehat{\nabla}^{\omega}_{V}$ receive the name of {\it gauge quantum linear connections} of $\omega$ (abbreviated ``gauge qlc's").
\end{Definition}

Using the graded Leibniz rule, we can extend $\nabla^{\omega}_{V}$ to the exterior covariant derivative $$d^{\nabla^{\omega}_{V}}: \Omega^\bullet(B)\otimes_B E^V_\l\longrightarrow \Omega^\bullet(B)\otimes_B E^V_\l$$ such that for all $\mu$ $\in$ $\Omega^k(B)$ 
\begin{equation}
\label{3.f10.1}
d^{\nabla^{\omega}_{V}}(\mu\otimes_B T)=d\mu\otimes_B T +(-1)^k\mu \nabla^{\omega}_{V}(T), 
\end{equation}
the curvature of $\nabla^{\omega}_{V}$ is defined as
\begin{equation}
\label{3.f10.2}
R^{\nabla^{\omega}_{V}}:=d^{\nabla^{\omega}_{V}}\circ \nabla^{\omega}_{V}:E^V_\l\longrightarrow \Omega^2(B)\otimes_B E^V_\l,
\end{equation}
and the following formula holds (\cite{sald2}):
\begin{equation}
\label{3.f10.3}
d^{\nabla^{\omega}_{V}}= \Upsilon_{V}\circ D^{\omega}\circ \Upsilon^{-1}_{V}.
\end{equation}
Similarly, using the graded Leibniz rule, we can extend $\widehat{\nabla}^{\omega}_{V}$ to the exterior covariant derivative  $$d^{\widehat{\nabla}^{\omega}_{V}}: E^V_\r \otimes_B \Omega^\bullet(B)\longrightarrow E^V_\r \otimes_B \Omega^\bullet(B),$$ which is given by  
\begin{equation}
\label{3.f10.4}
d^{\widehat{\nabla}^{\omega}_{V}}(T\otimes_B \mu)=\widehat{\nabla}^{\omega}_{V}(T) \mu +T \otimes_B d\mu , 
\end{equation}
the curvature is defined as
\begin{equation}
\label{3.f10.5}
R^{\widehat{\nabla}^{\omega}_{V}}:=d^{\widehat{\nabla}^{\omega}_{V}}\circ \widehat{\nabla}^{\omega}_{V}:E^V_\r\longrightarrow E^V_\r\otimes_B \Omega^2(B),
\end{equation}
and the following formula holds (\cite{sald2}):
\begin{equation}
\label{3.f10.6}
d^{\widehat{\nabla}^{\omega}_{V}}= \widehat{\Upsilon}_{V}\circ \widehat{D}^{\omega}\circ \widehat{\Upsilon}^{-1}_{V}.
\end{equation}

In accordance with Section 3.2 of reference \cite{sald2}, for every $\delta^V$ $\in$ $\Obj(\Rep_\G)$ we have that
\begin{equation}
    \label{modulestrucute}
    \begin{aligned}
        E^V_\l \cong B^{m_{V}} \cdot \rho^V \;\; &(\mbox{as left B--modules}), \\ &E^V_\r\cong \varrho^V\cdot B^{s_V} \;\; (\mbox{as right B--modules})
    \end{aligned}
\end{equation}
for some $m_{V}$, $s_V$ $\in$ $\N$, where $B^n=\underbrace{B\times \cdots \times B}_{n-times}$. Here, $\rho^V$ and $\varrho^V$ are idempotent square matrices with entries in $B$ of dimensions $m_V$ and $s_V$, respectively, such that $$\rho^{V\,\dagger}=\rho^V,\qquad \varrho^{V\,\dagger}=\varrho^V,$$ where $\dagger$ denotes the composition of
the $\ast$ operation with the usual matrix transposition.  For example, for $\delta^V$ $\in$ $\T$, we have (see Remark \ref{rema}) $$\rho^V=(\rho^V_{kl})=\left(\sum^{n_V}_{i=1}x^V_{ki}\,x^{V\,\ast}_{li}\right),\qquad \varrho^V=(\varrho^V_{kl})=\left(\sum^{n_V}_{i=1}x^{\overline{V}}_{ki}\,x^{\overline{V}\,\ast}_{li}\right).$$   

Therefore, there are non--degenerate \emph{canonical} Hermitian structures on $E^V_\l$, $E^V_\r$ (see Section 8.2 of \cite{lan} and Section 3.2 of \cite{sald2}).  Specifically, they are given by
\begin{equation}
    \label{canoleft}
        \langle -,-\rangle_\l: E^V_\l\times E^V_\l\longrightarrow B,\qquad 
        (T_1,T_2) \longmapsto \sum^n_{k=1}T_1(e_k)\,T_2(e_k)^\ast
\end{equation}
\begin{equation}
    \label{canoright}
        \langle -,-\rangle_\r:E^V_\r\times E^V_\r\longrightarrow B,\qquad
        (T_1,T_2) \longmapsto \sum^n_{k=1}T_1(e_k)^\ast\,T_2(e_k),
\end{equation}
where $\{e_k\}^n_{k=1}$ is any orthonormal basis (with respect to the inner product that makes $\delta^V$ unitary) of $V$. It is worth mentioning that  $\langle -,-\rangle_\l$, $\langle -,-\rangle_\r$ are  actually $B$--valued inner products (\cite{sald2}).

The canonical Hermitian structure on $E^V_\l$ can be extended to 
\begin{equation}
\label{algo11}
\begin{aligned}
\langle-,-\rangle_\l:(\Omega^\bullet(B)\otimes_B E^V_\l)\times (\Omega^\bullet(B)\otimes_B E^V_\l)\longrightarrow  \Omega^\bullet(B)
\end{aligned}
\end{equation}
by means of $\langle \mu_1\otimes_B T_1,\mu_2\otimes_B T_2\rangle_\l=\mu_1 \,\langle T_1,T_2\rangle_\l\,\mu^\ast_2.$ (\cite{lan}). Similarly, the canonical Hermitian structure on $E^V_\r$ can be extended to  
\begin{equation}
\label{algo12}
\begin{aligned}
\langle-,-\rangle_\r:(E^V_\r\otimes_B \Omega^\bullet(B))\times (E^V_\r\otimes_B \Omega^\bullet(B))\longrightarrow  \Omega^\bullet(B)
\end{aligned}
\end{equation}
by means of  $ \langle  T_1\otimes_B \mu_1, T_2\otimes_B \mu_2 \rangle_\r=\mu^\ast_1 \,\langle T_1,T_2\rangle_\r\,\mu_2$ (\cite{lan}). 

\begin{Definition}
    \label{definne}
    Let $\nabla: E^V_\l\longrightarrow \Omega^1(B)\otimes_B E^V_\l$ be a qlc, i.e., $\nabla$ is a linear map that satisfies the left Leibniz rule. We say that $\nabla$ is {\it compatible} with the metric $\langle-,-\rangle_\l$ if
    \begin{equation}
\label{3.f18}
\langle \nabla(T_1),T_2 \rangle_\l-\langle T_1,\nabla(T_2)\rangle_\l=d\langle T_1,T_2\rangle_\l 
\end{equation}
for all $T_1$, $T_2$ $\in$ $E^V_\l$. 

Similarly, let $\nabla: E^V_\r\longrightarrow E^V_\r\otimes_B\Omega^1(B) $ be a qlc, i.e., $\nabla$ is a linear map that satisfies the right Leibniz rule. We say that $\nabla$ is {\it compatible} with the metric $\langle-,-\rangle_\r$  if
    \begin{equation}
\label{3.f18.1}
\langle T_1,\nabla(T_2)\rangle_\r-\langle\nabla(T_1),T_2 \rangle_\r=d\langle T_1,T_2\rangle_\r
\end{equation}
for all $T_1$, $T_2$ $\in$ $E^V_\r$.
\end{Definition}

\begin{Definition}
\label{compatiblemetric}
    We define the space of all quantum linear connections compatible with the canonical Hermitian structure of $E^V_\l$ as $\mathfrak{qlc}_C(E^V_\l).$   In addition, we define the space of all quantum linear connections compatible with the canonical Hermitian structure of $E^V_\r$ as  $\mathfrak{qlc}_C(E^V_\r).$ 
\end{Definition}

In Theorem 3.10 of reference \cite{sald2}, the reader can find a proof of the following result, which is formulated under the convention that the differential is Hermitian. The reader should have no difficulty verifying that the theorem also holds under the anti--Hermitian convention.

\begin{Theorem}
\label{fgs}
Let $(\zeta,\omega)$ be a qpb with a qpc and $\delta^V$ $\in$ $\Obj(\Rep_{\G})$. Then $\nabla^\omega_\V$ $\in$ $\mathfrak{qlc}_C(E^V_\l)$ and $\widehat{\nabla}^\omega_\V$  $\in$ $\mathfrak{qlc}_C(E^V_\r)$.
\end{Theorem}

In this paper, we have defined a qpc as a linear map $\omega:\mathfrak{qg}^\#\longrightarrow \Omega^1(P)$ that satisfies 
$\Delta_{\Omega^\bullet(P)}(\omega(\theta))=(\omega\otimes \id_G)\ad(\theta)+\mathbbm{1}\otimes\theta$ and $\omega(\theta^\ast)=\omega(\theta)^\ast$; while in reference \cite{sald2}, a qpc is a linear map $\omega:\mathfrak{qg}^\#\longrightarrow \Omega^1(P)$ that only satisfies $\Delta_{\Omega^\bullet(P)}(\omega(\theta))=(\omega\otimes \id_G)\ad(\theta)+\mathbbm{1}\otimes\theta$, i.e., the condition $\omega(\theta^\ast)=\omega(\theta)^\ast$ is not necessary. In reference \cite{sald2}, when a qpc fulfills $\omega(\theta^\ast)=\omega(\theta)^\ast$ is called {\it real}; so the last theorem in reference \cite{sald2} is written in terms of real qpc's. We will discuss about the condition $\omega(\theta^\ast)=\omega(\theta)^\ast$ in the final section.

\begin{Proposition}
    \label{grassmann}
    Let $(\zeta,\omega)$ be a qpb with a qpc and consider $\delta^V$ a trivial $\G$--corepresentation on $V$. Then $\nabla^\omega_V$ is the Grassmann connection on $E^V_\l$, and $\widehat{\nabla}^\omega_V$ is the Grassmann connection on $E^V_\r$.
\end{Proposition}

\begin{proof}
    Since $\delta^V\cong \oplus_k \,\delta^\C_\triv$ for some $k$ $\in$ $\N$, it is sufficient to prove the statement for $\delta^\C_\triv$, the trivial $\G$--corepresentation on $\C$. Thus, the linear map $$T^\triv: \C \longrightarrow P $$ given by $T^\triv(1)=\mathbbm{1}$ is a $B$--bimodule basis of $\Mor(\delta^\C_\triv,\Delta_P)$ and for every $T$ $\in$ $\Mor(\delta^\C_\triv,\Delta_P)$ we have $$T=b^{_T}\,T^\triv=T^\triv\,b^{_T},$$ where $b^{_T}=T(1)$ $\in$ $B$. Actually, $\{T^\triv\}$ is the set of generators of $\Mor(\delta^\C_\triv,\Delta_P)$ given in Remark \ref{rema}. Thus, for every $w$ $\in$ $\C$, in light of equation (2.22) of reference \cite{micho3} we get 
    \begin{eqnarray*}
        D^\omega(T(w))=D^\omega(b^{_T} T^\triv(w))
        =
        db^{_T}\,T^\triv(w)+b^{_T} d(T^\triv(w))
        =
        db^{_T}\,T^\triv(w); 
    \end{eqnarray*}
     so  $D^\omega(T)=db^{_T}\,T^\triv$ and it follows that $$\nabla^\omega_\C(T)=\Upsilon_{\C}(D^\omega(T))= \Upsilon_{\C}(db^{_T}\,T^\triv )=db^{_T}\otimes_B T^\triv.$$
     
      By a similar calculation using equation (2.23) of reference \cite{micho3}, we obtain $$ \widehat{\nabla}^\omega_\C(T)=\widehat{\Upsilon}_{\C}(\widehat{D}^\omega(T))= \widehat{\Upsilon}_{\C}(T^\triv\, db^{_T})=T^\triv \otimes_B db^{_T}.$$   
\end{proof}

The last proposition shows that the maps 
\begin{equation}
    \label{conecleft}
    \mathfrak{qpc}(\zeta)\longrightarrow \mathfrak{qlc}_C(E^V_\l), \qquad \omega\longmapsto \nabla^\omega_V
\end{equation}
\begin{equation}
    \label{conecright}
    \mathfrak{qpc}(\zeta)\longrightarrow \mathfrak{qlc}_C(E^V_\r), \qquad \omega\longmapsto \widehat{\nabla}^\omega_V
\end{equation}
are, in general, neither surjective nor injective, as in the {\it classical} case (\cite{nodg}).

\section{Yang--Mills--Connes Theory}

In this section we are going to build a {\it link} between Connes' theory of spinors, Yang--Mills fields and the theory of quantum principal bundles. Clearly, Theorem \ref{fgs} gives us a first {\it link}  between both theories since Connes' formulation works for compatible quantum linear connections. However, under certain conditions, we can use all the {\it geometrical} structure of a qpb in order to define a new Yang--Mills functional.

\subsection{Generalities}
Since Connes' theory is  more well--known than Durdevich's formulation of quantum principal bundles, and to maintain a reasonable length for this paper, we will not present the foundational aspects of this framework here. Nevertheless, the reader can consult \cite{con} for further details.

Let 
\begin{equation}
    \label{3.3.f1}
    (B,\mathfrak{H},\D)
\end{equation}
be a $n$--dimensional spectral triple in the sense of \cite{lan}, where $B$ is a $\ast$--(sub)algebra stable under functional calculus \cite{con}. We will identify $B$ with its representation in $\mathcal{B}(\mathfrak{H})$, the space of bounded linear operators in the Hilbert space $(\mathfrak{H},\langle-|-\rangle_\mathfrak{H})$.

Take  $(\Omega^\bullet_U(B),d_U,\ast)$  the universal graded differential $\ast$--algebra over $B$  (\cite{lan}). We define
\begin{equation}
    \label{3.3.f2}
    \pi_\r: \Omega^\bullet_U(B)\longrightarrow \mathcal{B}(\mathfrak{H})
\end{equation}
given by $$\pi_\r(b_0\,d_Ub_1\,\cdots d_Ub_n)=b_0 [\D,b_1]\cdots [\D,b_n].$$ 
It is well--known that $\pi_\r$ is a $\ast$--representation  and considering 
 $J^\bullet=J_0+d_U\,J_0,$  where $J_0=\bigoplus_n J^n_0$ and $$J^n_0:=\{ \mu\in \Omega^n_u(B) \mid \pi_\r(\mu)=0 \},$$ we obtain the Connes' differential forms space (\cite{lan})
\begin{equation}
    \label{3.3.f3}
    \Omega^\bullet_\D(B):=\Omega^\bullet_U(B)/J^\bullet\cong \pi_\r(\Omega^\bullet_U(B))/\pi_\r(d_U J_0),
\end{equation}
\begin{equation}
    \label{3.3.f4}
        d_\D: \Omega^n_\D(B)\longrightarrow \Omega^{n+1}_\D(B),\qquad
        [\mu]\longmapsto d_\D[\mu]:=[d_U\mu]\cong [\pi_\r(d_U \mu)].
\end{equation}

For any quantum group $\G$, let us take any quantum principal $\G$--bundle over $B$
\begin{equation}
    \label{3.3.f5}
    \zeta=(P,B,\Delta_P).
\end{equation}

\begin{Remark}
  \label{remadirac}
From now on, we will consider that $\zeta$ is equipped with a differential calculus such that $(\Omega^\bullet_\D(B),d_\D,\ast)$ is the space of quantum base forms.  
\end{Remark}

It is worth mentioning that such a structure always exists. Indeed,
let $(B,\mathfrak{H},\mathcal{D})$ be a $n$--dimensional spectral triple and $\G$ a quantum group. If there does not exist a \emph{natural} qpb structure such that Remark \ref{remadirac} holds, one can always take the trivial qpb 
\begin{equation}
    \label{3.3.f6}
    \zeta=(P=B\otimes G,B,\Delta_P=\id_B\otimes \Delta).
\end{equation}
with the trivial differential calculus (taking the corresponding tensor products)
\begin{equation}
    \label{3.3.f6.1}
    \Omega^\bullet(P)=\Omega^\bullet(B)\otimes \Gamma^\wedge,\qquad \Delta_{\Omega^\bullet(P)}=\id_{\Omega^\bullet(B)}\otimes \Delta,
\end{equation}
where $(\Gamma^\wedge,d,\ast)$ is the universal differential envelope $\ast$--calculus of some bicovariant $\ast$--FODC $(\Gamma,d)$ over $\G$, and $\Delta:\Gamma^\wedge\longrightarrow \Gamma^\wedge\otimes \Gamma^\wedge$ is the extension of the coproduct $\Delta:G\longrightarrow G\otimes G$ on $\Gamma^\wedge$ (\cite{micho1}).

When the regularity degree of $B$ is $B^2=B$ (i.e., more than Lipschitz), in light of \cite{con}, there is a canonical {\it quantum integration} on $\pi_\r(\Omega^\bullet_U(B))$ given by
\begin{equation}
\label{3.3.f17.1}
    \int_B \mu=\mathrm{Tr}_w(\mu\,|\D|^{-n}),
\end{equation}
where $\mathrm{Tr}_w$ denotes the Dixmier trace. This also defines an inner product on $\pi_\r(\Omega^\bullet_U(B))$ given by 
\begin{equation}
\label{3.3.f17.2}
    \langle \mu_1 \mid \mu_2\rangle_\r=\int_B \mu^\ast_1\,\mu_2= \mathrm{Tr}_w(\mu^\ast_1\,\mu_2\,|\D|^{-n})
\end{equation}
for $\mu_1$, $\mu_2$ $\in$ $\pi_\r(\Omega^m_U(B))$ and $0$ for differential forms of different degree. Let $\mathfrak{H}'$ be the Hilbert space completion of $\pi_\r(\Omega^m_U(B))$ with respect to the previous inner product. If $\mathcal{P}$ denotes orthogonal projection of $\mathfrak{H}'$ on the orthogonal complement of $\pi_\r(d_U J^{m-1}_0)$, we define a quantum integration on $\Omega^\bullet_\D(B)$ given by
\begin{equation}
\label{3.3.f17}
    \int_B \,[\mu]=\mathrm{Tr}_w(\mathcal{P}(\mu)\,|\D|^{-n}),
\end{equation}
and an inner product on $\Omega^\bullet_\D(B)$ given by
\begin{equation}
\label{3.3.f18}
    \langle [\mu_1] \mid [\mu_2]\rangle_\r=\int_B \,[\mu_1]^\ast\,[\mu_2]=\mathrm{Tr}_w(\mathcal{P}(\mu_1)^\ast\,\mathcal{P}(\mu_2)\,|\D|^{-n}),
\end{equation}
for $[\mu_1]$, $[\mu_2]$ $\in$ $\Omega^m_\D(B)$ and $0$ for differential forms of different degree. For more details, see reference \cite{con}.

\subsection{The Gauge Dirac Operator}

With all the generalities established, we proceed to show the Connes' theory of spinors {\it tight} to the theory of quantum principal bundles.

\begin{Definition}
    \label{gaugespinors}
    Let $\delta^V$ $\in$ $\Obj(\Rep_\G)$. We define the space of non--commutative geometrical (right) gauge spinor space as  $\qgs_\r:= E^V_\r\otimes_B \mathfrak{H}.$
\end{Definition}

According to \cite{con,lan}, the map $\pi_R$ induces a projection of the space of qlc's 
$$^{U}\widehat{\nabla}:E^V_\r\longrightarrow E^V_\r\otimes_B \Omega^1_U(B)$$  build over $(\Omega^\bullet_U(B),d_U,\ast)$
compatible with the metric, to the space of qlc's
$$\widehat{\nabla}:E^V_\r\longrightarrow E^V_\r\otimes_B \Omega^1_\D(B)$$ built over $(\Omega^\bullet_\D(B),d_\D,\ast)$
compatible with the metric. In light of references \cite{con,lan} and Theorem \ref{fgs}, we define

\begin{Definition}
    \label{diracoperator}
    Let $\delta^V$ $\in$ $\Obj(\Rep_\G)$ and let $\omega$ be a qpc. We define the non--commutative (right) gauge Dirac operator induced by $\omega$ (or by $\widehat{\nabla}^\omega_\V$)  $$\D_{\widehat{\nabla}^\omega_V}: \qgs_\r \longrightarrow \qgs_\r  $$ as the operator $$\D_{\widehat{\nabla}^\omega_V}(T\otimes_B \psi)=T\otimes_B \D\psi+\left((\id_{E^V_\r}\otimes \pi_\r)\circ\, {^{U}\widehat{\nabla}^\omega_V}T \right)\psi,$$ where ${^{U}\widehat{\nabla}^\omega_V}:E^V_\r\longrightarrow E^V_\r\otimes_B \Omega^1_U(B)$ is a qlc for $(\Omega^\bullet_U(B),d_U,\ast)$ compatible with the metric that projects to the gauge qlc $\widehat{\nabla}^\omega_V$, and the action $((\id_{E^V_\r}\otimes \pi_\r)\circ\, {^{U}\widehat{\nabla}^\omega_V}T)\,\psi$  is component--wise.
\end{Definition}

According to Lemma 1 Section 6.3 of reference \cite{con}, we have

\begin{Theorem}
    \label{diracright}
    Let $\delta^V$ $\in$ $\Obj(\Rep_\G)$ and let $\omega$ be a qpc.
    \begin{enumerate}
        \item The space of non--commutative geometrical (right) gauge spinor space can be promoted to a Hilbert space by defining 
 \begin{equation*}
    \begin{aligned}
        \langle-|-\rangle^\bullet_{\mathfrak{H}}: \qgs_\r \times \qgs_\r \longrightarrow \C
    \end{aligned}
\end{equation*}
given by $$\langle  T_1\otimes_B \psi_1 \mid  T_2\otimes_B \psi_2\rangle^\bullet_{\mathfrak{H}}:=\langle\psi_1\mid \langle T_1,T_2\rangle_\r\, \psi_2\rangle_\mathfrak{H}$$ for all $T_1$, $T_2$ $\in$ $E^V_\r$, $\psi_1$, $\psi_2$ $\in$ $\mathfrak{H}$. 
        \item The operator $\D_{\widehat{\nabla}^\omega_V}$ is self--adjoint with domain  $E^V_\r\otimes_B \mathrm{Dom}(\D)$.
    \end{enumerate}  
\end{Theorem}

Now consider the following two definitions

\begin{Definition}
    \label{diracaction}
     Let $\delta^V$ $\in$ $\Obj(\Rep_\G)$. We define the non--commutative geometrical Dirac action as the functional 
     \begin{equation*}
     \begin{aligned}
         \qS_\dirac: \mathfrak{qpc}(\zeta)\times \qgs_\r\longrightarrow \C,\qquad
         (\omega,\Psi)\longmapsto \langle \Psi\mid\D_{\widehat{\nabla}^\omega_V}\Psi\rangle^\bullet_\mathfrak{H}.
     \end{aligned}
     \end{equation*}
\end{Definition}

\begin{Definition}
\label{qgg}
    We define the quantum gauge group of the qpb $\zeta$ as
    \begin{equation*}
        \begin{aligned}
            \qGG=\{\F:&\Omega^\bullet(P)\longrightarrow \Omega^\bullet(P)\mid \F  \mbox{ is a } \mbox{graded left } \Omega^\bullet(B)-\mbox{module isomorphism  }\\   &\;\mbox{such that}\;\;\;\F(\mathbbm{1})=\mathbbm{1}, \; \Delta_{\Omega^\bullet(P)}\circ \F=(\F\otimes \id_{\Gamma^\wedge})\circ \Delta_{\Omega^\bullet(P)},
            \\
            & \qquad\qquad\quad\F(p^\ast)=\F(p)^\ast \;\mbox{ for all }\; p \in P\; \mbox{ and }
            \\ 
     &\qquad\quad\F(\Im(\omega)^\ast)=\F(\Im(\omega))^\ast \;\mbox{ for all }\; \omega \in \mathfrak{qpc}(\zeta) \}.
        \end{aligned}
    \end{equation*}
\end{Definition}
Elements of $\qGG$ are called {\it quantum gauge transformations}. It is worth mentioning that in general, a quantum gauge transformation {\bf does not} commute with the corresponding differential. 

In accordance with Theorem 4.7 points 1 and 2 of reference \cite{sald2}, $\qGG$ has a well--defined group action on the space $\mathfrak{qpc}(\zeta)$ by 
\begin{equation}
    \label{actionqpc}
    \F^\circledast \omega:=\F\circ \omega,
\end{equation}
and this formula induces a well--defined action on the curvature (see Proposition 4.8 of reference \cite{sald2}). The reader is encouraged to consult the reference \cite{sald2} for more details. Notice that $\F^\circledast \omega$ is only the {\it dualization} of the action of the gauge group on principal connections via the pull--back in differential geometry.

On the other hand, for any $\delta^V$ $\in$ $\Obj(\Rep_\G)$, the group $\qGG$ acts on the associated right qvb by means on (see Section 4.3 of reference \cite{sald2})
\begin{equation}
    \label{actionsection}
    \mathbf{A}_\F: E^{V}_\r\longrightarrow E^{V}_\r,\qquad \mathbf{A}_\F(T)(v)=\F(T(v)).
\end{equation}
Since $\F\circ \ast|_P=\ast \circ \F|_P$, it follows that $\F(p\,b)=\F(p)\,b$ for all $p$ $\in$ $P$, $b$ $\in$ $B$ and hence, the action of $\qGG$ on $E^V_\r$ extends to $\qgs_\r=E^V_\r\otimes_B \mathfrak{H}$ component--wise. 

In general, the quantum gauge group is quite large and is very challenging to work with. Fortunately, in each specific situation, one can work with suitable \emph{ad hoc} subgroups.

\begin{Definition}
\label{diracgauge}
We define the quantum gauge group of the non--commutative geometrical Dirac action as the group
\begin{equation*}
    \begin{aligned}
        \qGG_\dirac:=\{\F \in \qGG\mid \qS_\dirac(\omega,\Psi)&=\qS_\dirac(\F^{\circledast} \omega,\mathbf{A}_\F(\Psi))\\ 
        &\mbox{ for all }\, \omega \in \mathfrak{qpc}(\zeta),\;\Psi\, \in\, \qgs_\r  \} \subseteq \qGG.
    \end{aligned}
\end{equation*}
\end{Definition}

\begin{Proposition}
    \label{gaugesclar}
    Let $\F:\Omega^\bullet(P)\longrightarrow \Omega^\bullet(P)$ be a quantum gauge transformation such that $\F$ is a graded differential $\ast$--algebra morphism. Then $\F$ $\in$ $\qGG_\dirac$.
\end{Proposition}

\begin{proof}
  According to Section 4.3 of reference \cite{sald2}, under the hypothesis of this proposition, we have that $$\mathbf{A}_\F \;\in\; U(E^V_\r),$$ where $U(E^V_\r)$ is the group of all unitary operators of $E^V_\r$ with respect to $\langle-,-\rangle_\r$. In light of Section 4.3 of reference \cite{sald2}, we get $$\widehat{\nabla}^{\F^{\circledast}\omega}_{V}=(\mathbf{A}_\F\otimes_B\id_{\Omega^\bullet_\D(B)} )\circ \,{\widehat{\nabla}^{\omega}_{V}}\circ \mathbf{A}^{-1}_{\F}.$$ 
  
  Let $^U\widehat{\nabla}^\omega_{V}$ be a qlc for $(\Omega^\bullet_U(B),d_U,\ast)$ compatible with the metric that projects to $\widehat{\nabla}^{\omega}_{V}$, i.e., $$\widehat{\nabla}^{\omega}_{V}=(\id_{E^V_\r}\otimes_B \pi_\r)\;{^U}\widehat{\nabla}^\omega_{V}.$$ So 
  \begin{eqnarray*}
      \widehat{\nabla}^{\F^{\circledast}\omega}_{V}&=&(\mathbf{A}_\F\otimes_B\id_{\Omega^\bullet_\D(B)} )\circ \,{\widehat{\nabla}^{\omega}_{V}}\circ \mathbf{A}^{-1}_{\F}
      \\
      &=&(\mathbf{A}_\F\otimes_B\id_{\Omega^\bullet_\D(B)} )\circ \,(\id_{E^V_\r}\otimes_B \pi_\r)\,{^U}\widehat{\nabla}^\omega_{V}\circ \mathbf{A}^{-1}_{\F}
      \\
      &=&
      (\mathbf{A}_\F\otimes_B \pi_\r)\circ \,{^U}\widehat{\nabla}^\omega_{V}\circ \mathbf{A}^{-1}_{\F}
      \\
      &=&
      (\id_{E^V_\r}\otimes_B \pi_\r)\circ (\mathbf{A}_\F\otimes_B\id_{\Omega^\bullet_U(B)})\circ \,{^U}\widehat{\nabla}^\omega_{V}\circ \mathbf{A}^{-1}_{\F}
      \\
      &=&
      (\id_{E^V_\r}\otimes_B \pi_\r)\circ {^U}\widehat{\nabla}^{\F^{\circledast}\omega}_{V},
  \end{eqnarray*}
  where  $${^U}\widehat{\nabla}^{\F^{\circledast}\omega}_{V}=(\mathbf{A}_\F\otimes_B\id_{\Omega^\bullet_U(B)} )\circ \,^U{\widehat{\nabla}^{\omega}_{V}}\circ \mathbf{A}^{-1}_{\F}. $$ Notice that  ${^U}\widehat{\nabla}^{\F^{\circledast}\omega}_{V}$ is given by the {\it standard} action of $U(E^V_\r)$ on $^U{\widehat{\nabla}^{\omega}_{V}}$ and it is well--known that this action sends qlc's compatible with the metric to qlc's compatible with the metric (\cite{lan}); therefore, ${^U}\widehat{\nabla}^{\F^{\circledast}\omega}_{V}$ is a qlc  with respect to $(\Omega^\bullet_U(B),d_U,\ast)$ compatible with the metric that projects to $\widehat{\nabla}^{\F^{\circledast}\omega}_{V}$.  
  
  If  $^U\widehat{\nabla}^\omega_V(T)=\displaystyle\sum_{k}T_k\otimes_B \mu_k$, we have $$^U\widehat{\nabla}^{\F^{\circledast}\omega}_{V}(\mathbf{A}_{\F}(T))= (\mathbf{A}_\F\otimes_B\id_{\Omega^\bullet_U(B)} )\,{^U}\widehat{\nabla}^{\omega}_{V}T=\displaystyle\sum_{k} \mathbf{A}_\F(T_k)\otimes_B \mu_k.$$ Thus, for $\Psi=T\otimes_B \psi$ $\in$ $\qgs_\r$ we obtain
  
\begin{eqnarray*}
 \langle \mathbf{A}_\F(\Psi)\mid\D_{\widehat{\nabla}^{\F^{\circledast}\omega}_V}\mathbf{A}_\F(\Psi)\rangle^\bullet_\mathfrak{H}&=&\langle \mathbf{A}_\F(T)\otimes_B \psi \mid \mathbf{A}_\F(T)\otimes_B \D\psi
 \\
 &+&
 \left((\id_{E^V_\r}\otimes \pi_\r)\circ\, {^{U}\widehat{\nabla}^{\F^{\circledast}\omega}_V}\mathbf{A}_\F(T)\right)\psi\rangle^\bullet_\mathfrak{H}
  \\
  &= &
 \langle\psi \mid \langle \mathbf{A}_\F(T),\mathbf{A}_\F(T) \rangle_\r\,\D\psi\rangle_{\mathfrak{H}}
 \\
 &+&
 \sum_k \langle\psi \mid \langle \mathbf{A}_\F(T),\mathbf{A}_\F(T_k) \rangle_\r \,\pi_\r(\mu_k)\psi \rangle_{\mathfrak{H}}
  \\
  &= &
  \langle\psi \mid \langle T,T \rangle_\r\,\D\psi\rangle_{\mathfrak{H}}
  \\
 &+&
  \sum_k \langle\psi \mid \langle T,T_k \rangle_\r \,\pi_\r(\mu_k)\psi \rangle_{\mathfrak{H}}
  \\
  &= &
  \langle T \otimes_B\psi \mid T \otimes_B \D\psi
  \\
 &+&
  \left((\id_{E^V_\r}\otimes \pi_\r)\circ\, {^{U}\widehat{\nabla}^{\omega}_V}T\right)\psi \rangle^\bullet_{\mathfrak{H}}
  \\
  &= &
  \langle \Psi\mid\D_{\widehat{\nabla}^{\omega}_V}\Psi\rangle^\bullet_\mathfrak{H}
\end{eqnarray*}
and proposition follows by linearity. 
\end{proof}

Following the classical case, we have

\begin{Definition}
\label{6.1.14}
A stationary point of $\qS_\dirac$ is an element $\Psi$ $\in$ $\qgs_\r$ such that for all $\Phi$ $\in$ $\qgs_\r$ we have  $$\left.\dfrac{\partial}{\partial z}\right|_{z=0} \qS_\dirac(\omega,\Psi+z\,\Phi)=0$$ for  $z$ $\in$ $\C$. Stationary points are also called non--commutative geometrical Dirac spinor fields and in terms of a traditional physical interpretation, they can be interpreted as {\it Dirac spinor fields possessing the symmetry $\qGG_\dirac$}.
\end{Definition}

Notice that $\qS_\dirac(\omega,\Psi+ z\, \Phi)$ depends on a polynomical way in $z$, so we have considered the derivative of $\qS_\dirac$ as a formal derivative in the obviously way. The proof of the following theorem is a straightforward calculation, so it will be omitted. 

\begin{Theorem}
\label{6.1.15}
An element $\Psi$ $\in$ $\qgs_\r$ is a non--commutative geometrical Dirac spinor if and only if 
\begin{equation}
\label{3.3.f19}
    \D_{\widehat{\nabla}^{\omega}_V}\Psi=0
\end{equation}
We will refer to equation (\ref{3.3.f19}) as the {\it non--commutative geometrical Dirac equation}.
\end{Theorem}

According to \cite{lan3,dabro}, in order to {\it play} with left structure, it is necessary to consider a real $n$--dimensional spectral triple 
$$(B,\mathfrak{H},\D,\mathcal{J})$$ (or a twisted $n$--dimensional spectral triple),  the opposite $\ast$--algebra 
\begin{equation}
    \label{3.3.f20}
    B^\circ \quad \mbox{ with }\quad b^\circ:= \mathcal{J}b^\ast\mathcal{J}^{-1}
\end{equation}
 and the $\ast$--representation $\pi_\l$
\begin{equation}
    \label{3.3.f22}
   \pi_\l: \Omega^\bullet_U(B)\longrightarrow \mathcal{B}(\mathfrak{H}) 
\end{equation}
given by $$\pi_\l(b_0\,d_Ub_1\cdots d_Ub_n)=b^\circ_0[\D,b^\circ_1]\cdots [\D,b^\circ_n].$$ In this way, {\it mutatis mutandis}, we can recreate all the theory presented in this subsection for $E^V_\l$, the corresponding quotient graded differential $\ast$--algebra, now denoting by $\Omega^\bullet_\D(B^\circ)$ (\cite{lan3,dabro}), and the gauge qlc $$\nabla^\omega_\V:E^V_\l\longrightarrow \Omega^1_\D(B^\circ)\otimes_B E^V_\l.$$

\subsection{Yang--Mills Theory}

The Connes' theory also involves Yang--Mills theory, and in this section we will focus on that aspect.  We will still work in the context of Remark \ref{remadirac} with the regularity degree of $B$ more than Lipschitz (i.e., $B^2=B$). To simplify the notation, in this section, elements of $\Omega^\bullet_\D(B)$ will be denoted without the bracket $[-]$. This section constitutes the core of the paper, and Theorem \ref{soperatoradjunto} is its most important result.

Let $\delta^V$ $\in$ $\Obj(\Rep_\G)$ and take $$\widehat{\nabla}: E^V_\r\longrightarrow E^V_\r\otimes_B \Omega^1_\D(B)$$ any qlc compatible with the metric. In addition, consider $$d^{\widehat{\nabla}}:E^V_\r\otimes_B\Omega^k_\D(B)\longrightarrow E^V_\r\otimes_B\Omega^{k+1}_\D(B) $$ its extension by the graded Leibniz rule  (as in equation (\ref{3.f10.4})), and $$R^{\widehat{\nabla}}:=d^{\widehat{\nabla}}\circ \widehat{\nabla}: E^V_\r\longrightarrow E^V_\r\otimes_B \Omega^2_\D(B)$$ its curvature. Every right $B$--linear map $$\widehat{\Lambda}: E^V_\r\longrightarrow E^V_\r\otimes_B \Omega^k_\D(B)$$ can be extended to a right $\Omega^\bullet(B)$--linear map $$\widehat{\Lambda}^{\mathrm{ext}}: E^V_\r\otimes_B \Omega^\bullet_\D(B) \longrightarrow E^V_\r\otimes_B \Omega^\bullet_\D(B),\qquad  T\otimes_B \mu\longmapsto \widehat{\Lambda}(T)\,\mu$$ for all elementary elements $T\otimes_B \mu$ $\in$ $E^V_\r\otimes_B\Omega^\bullet_\D(B)$. In particular, the curvature map is right $B$--linear and  (\cite{lan}) $$R^{\widehat{\nabla}\,\mathrm{ext}}=d^{\widehat{\nabla}}\circ d^{\widehat{\nabla}}.$$ 

We define the right $B$--linear map 
\begin{equation}
    \label{4.3.f25}
    [\widehat{\nabla},\widehat{\Lambda}]:=d^{\widehat{\nabla}} \circ \widehat{\Lambda} -(-1)^k\,\widehat{\Lambda}^{\mathrm{ext}}\circ \widehat{\nabla}: E^V_\r \longrightarrow E^V_\r\otimes_B \Omega^{k+1}_\D(B) 
\end{equation}
for a right $B$--linear map $\widehat{\Lambda}:E^V_\r\longrightarrow E^V_\r\otimes_B \Omega^k_\D(B).$ Notice that
\begin{equation}
    \label{4.3.f25.1}
[\widehat{\nabla},R^{\widehat{\nabla}}]=d^{\widehat{\nabla}}\circ d^{\widehat{\nabla}}\circ \widehat{\nabla}-d^{\widehat{\nabla}}\circ d^{\widehat{\nabla}}\circ \widehat{\nabla}=0 .
\end{equation}

Consider the  inner product
\begin{equation}
    \label{4.3.f26}
    \langle-|-\rangle^\bullet_\r: (E^V_\r\otimes_B \Omega^2_\D(B))\times (E^V_\r\otimes_B \Omega^2_\D(B))\longrightarrow \C
\end{equation}
given by (see equation (\ref{3.3.f18})) $$ \langle T_1\otimes_B \mu_1\mid T_2\otimes_B \mu_2\rangle^\bullet_\r=\langle \,\mu_1\mid \langle T_1,T_2\rangle_\r\,\mu_2\,\rangle_\r,$$ and the inner product 
\begin{equation}
    \label{4.3.f27}
    \langle\langle-|-\rangle\rangle_\r:(\mathrm{End}(E^V_\r)\otimes_B \Omega^2_\D(B))\times (\mathrm{End}(E^V_\r)\otimes_B \Omega^2_\D(B))\longrightarrow \C
\end{equation}
given by  $$\langle\langle \widehat{\Lambda}_1\mid \widehat{\Lambda}_2\rangle\rangle_\r=\sum_k \langle \widehat{\Lambda}_1(T^\r_k)\mid \widehat{\Lambda}_2(T^\r_k)\rangle^\bullet_\r.$$
In accordance with Section 3 of reference \cite{sald2}, we have $\varrho^V(\overline{e}_k)=T^\r_k$ (see equation (\ref{modulestrucute})), where $\overline{e}_k$ is the canonical basis of $B^{s_V}$. In this way, the inner product of equation (\ref{4.3.f27}) agrees with the one proposed on references \cite{guin,guin2}.  
\begin{Definition}
    \label{yangmills}
    We define the non--commutative analytical Yang--Mills functional as (see Definition \ref{compatiblemetric})
\begin{equation*}
    \begin{aligned}
        \YM: \mathfrak{qlc}_{C}(E^V_\r)\longrightarrow \R,\qquad
        \widehat{\nabla} \longmapsto || R^{\widehat{\nabla}}||^2:=\langle\langle R^{\widehat{\nabla}}\mid R^{\widehat{\nabla}}\rangle\rangle_\r.
    \end{aligned}
\end{equation*}
\end{Definition}
This functional is invariant under the standard action of the group $U(E^V_\r)$ of unitary operators of $E^V_\r$ \cite{con,lan,guin,guin2}. When the operator $$[\widehat{\nabla},-]$$ {\it admits a formally adjoint operator} 
\begin{equation}
    \label{4.3.f28}
    [\widehat{\nabla},-]^\star
\end{equation}
with respect to $\langle\langle-|-\rangle\rangle_\r$ for all qlc $\widehat{\nabla}$ compatible with the metric, we can perform a variation on the functional $\YM$ of the form $$\widehat{\nabla}\longmapsto \widehat{\nabla} +t \,\widehat{\Lambda}$$ for $t$ $\in$ $\R$ and $\widehat{\Lambda}:E^V_\r \longrightarrow E^V_\r\otimes_B \Omega^1_\D(B)$ being a Hermitian right $B$--linear map,  and one obtains the {\it non--commutative analytical Yang--Mills equation} 
\begin{equation}
    \label{4.3.f29}
        [\widehat{\nabla},R^{\widehat{\nabla}}]^\star=0.
\end{equation}
\noindent This happens, for example, for toric non--commutative manifolds as the reader can verify in reference \cite{lan4,lan2}. We will refer to equation (\ref{4.3.f29}) as the {\it non--commutative analytical Yang--Mills equation} and solutions of this equation will be called {\it Yang--Mills quantum linear connections} (Yang--Mills qlc's).

\begin{Remark}
\label{remayanmills}
     In literature, for example,  reference \cite{lan2}, there are developments in which the authors used simultaneously the theory of quantum principal bundles and Connes' theory, as us in this paper. However, all these works are developed for some concrete cases; while this paper provides a general framework.
     
     In principle, our approach does not provide additional information about the quantum spaces themselves; rather, it shows that Connes' theory can be applied to \textbf{any $n$--dimensional spectral triple  as the base space of a quantum principal $\G$--bundle for any quantum group $\G$}, if Remark \ref{remadirac} holds and $B^2=B$. Moreover, by equations (\ref{3.3.f6}) and~(\ref{3.3.f6.1}), one can always construct quantum principal bundles satisfying Remark \ref{remadirac}. Notice that, due to the structure of these equations, any bicovariant $\ast$--FODC over $\G$ can be used.
     
 In other words, the significance of this part of the manuscript lies, for instance, in the fact that we have proven that the gauge Dirac operator associated with $\omega$ can be defined and studied for every qpc $\omega$, precisely as is done in differential geometry.
\end{Remark}

In the rest of this section, we will present another {\it non--commutative   Yang--Mills equation} than the one of equation (\ref{4.3.f29}), for another {\it non--commutative} Yang--Mills functional.

Let $\pi_X: Y\longrightarrow X$ be a principal $H$--bundle. Then, the Yang--Mills functional can be defined as (\cite{gtvp,diff1})
\begin{equation}
\label{classicalyangmills}
        {\qS_\YM}_{\mathrm{class}}: \mathfrak{pc}(\pi_X)\longrightarrow \R,\qquad
        \omega \longmapsto ||\Omega^\omega||^2:=\langle GP(\Omega^\omega)\mid GP(\Omega^\omega)\rangle, 
\end{equation}
 where $\mathfrak{pc}(\pi_X)$ denotes the set of all principal connection of $\pi_X: Y\longrightarrow X$ and $\Omega^\omega$ is the curvature of the principal connection $\omega$  (\cite{diff1}). Here, the map $GP$ is the canonical isomorphism between basic differential $m$--forms of $Y$ of type $\ad$ (with $\ad$ the adjoint representation of $H$ on its Lie algebra) and differential $m$--forms of $X$  with values in the associated vector bundle $E^{\mathfrak{h}}$ of $\ad$ (\cite{nodg,gtvp,diff1}).

By performing a variation of the form $$\omega\longmapsto \omega+t\,\lambda,$$ on the functional ${\qS_\YM}_{\mathrm{class}}$, where $\lambda$ is a basic $1$--form of $Y$ of type $\ad$, one obtains the  {\it classical} Yang--Mills equation 
\begin{equation}
    \label{classicalymequation}
    d^{\nabla^\omega_{\mathfrak{h}}\star}(GP(\Omega^\omega))=0,
\end{equation}
where $d^{\nabla^\omega_{\mathfrak{h}}\star}$ is the formal adjoint operator of the exterior covariant derivative of the gauge linear connection $\nabla^\omega_{\mathfrak{h}}$ \cite{gtvp,diff1}.

Let $$\zeta=(P,B,\Delta_P)$$ be a qpb with a differential calculus as in the previous subsection. Assume that the $\G$--corepresentation $\ad$ of equation (\ref{2.f15}) is finite--dimensional, in other words, $\ad$ $\in$ $\Obj(\Rep_\G)$. Hence, we can apply all the theory of Section 4 to $\ad$. In particular, we can consider the associated right qvb $E^{\mathfrak{qg}^\#}_\r$ and $\E^{\mathfrak{qg}^\#}_\r:=\Mor(\ad,\Delta_\Hor)$. It is worth remembering that the {\it non--commutative geometrical counterpart}  of the map $GP$ is the map (see \cite{sald2}) $$\widehat{\Upsilon}_{\mathfrak{qg}^\#}: \E^{\mathfrak{qg}^\#}_\r\longrightarrow E^{\mathfrak{qg}^\#}_\r\otimes_B \Omega^\bullet_\D(B);$$ while the {\it non--commutative geometrical counterpart} of the inner product us in equation (\ref{classicalyangmills}) is the inner product of equation (\ref{4.3.f26}) for $E^{\mathfrak{qg}^\#}_\r$. Moreover, we have defined the curvature of a qpc $\omega$ in such a way that $R^\omega$ $\in$ $\E^{\mathfrak{qg}^\#}_\r$ with $\Im(R^\omega)$ $\subseteq$ $\Hor^2 P$  through the dualization of the {\it classical} case (see Section 3). Therefore,  it should be natural to consider the following functional \cite{saldym}.

\begin{Definition}
    \label{yangmills1}
    We define the non--commutative geometrical Yang--Mills functional as 
\begin{equation*}
        \qS_\YM: \mathfrak{qpc}(\zeta)\longrightarrow \R,\qquad
        \omega \longmapsto || R^{\omega}||^2:=\langle \widehat{\Upsilon}_{\mathfrak{qg}^\#}(R^{\omega})\mid \widehat{\Upsilon}_{\mathfrak{qg}^\#}(R^{\omega})\rangle^\bullet_\r. 
\end{equation*}
\end{Definition}

It is worth mentioning that in Definition \ref{yangmills}, we measure the square norm of the curvature of a qlc compatible with the metric of $E^V_\r$ for some $\delta^V$ $\in$ $\Obj(\Rep_\G)$; while in Definition \ref{yangmills1}, we measure the square norm of the curvature of a qpc of $\zeta$, as in the equation (\ref{classicalyangmills}). As in the last subsection, we can define

\begin{Definition}
\label{6.1.3}
We define the quantum gauge group of the non--commutative Yang-Mills functional as the group (see Definition \ref{qgg} and equation (\ref{actionqpc})) $$\qGG_\YM:=\{\F \in \qGG\mid \qS_\YM(\omega)=\qS_\YM(\F^{\circledast} \omega)\; \mbox{ for all }\; \omega \,\in\, \mathfrak{qpc}(\zeta) \}\, \subseteq\, \qGG.$$  
\end{Definition}
\noindent In addition, the proof of the following statement can be found in  Proposition 31 of reference \cite{saldym}

\begin{Proposition}
    \label{gaugeym}
    Let $\F:\Omega^\bullet(P)\longrightarrow \Omega^\bullet(P)$ be a quantum gauge transformation such that $\F$ is a graded differential $\ast$--algebra morphism. Then $\F$ $\in$ $\qGG_\YM$.
\end{Proposition}

Consider the dual twisted covariant derivative $\widehat{DS}^\omega=\widehat{D}^\omega-\widehat{S}^\omega$ (see Definition \ref{twisted}) and the operator (see equations (\ref{3.f10.4}), (\ref{3.f10.6}))
\begin{eqnarray}
        \label{ec.000}
d^{\widehat{DS}^\omega} &:=& \widehat{\Upsilon}_{\mathfrak{qg}^\#} \circ \widehat{DS}^\omega \circ \widehat{\Upsilon}^{-1}_{\mathfrak{qg}^\#}\nonumber
\\
&=&\widehat{\Upsilon}_{\mathfrak{qg}^\#} \circ \widehat{D}^\omega \circ \widehat{\Upsilon}^{-1}_{\mathfrak{qg}^\#}-\widehat{\Upsilon}_{\mathfrak{qg}^\#} \circ \widehat{S}^\omega \circ \widehat{\Upsilon}^{-1}_{\mathfrak{qg}^\#}
\\
&=&
d^{\widehat{\nabla}^\omega_{\mathfrak{qg}^\#}}-d^{\widehat{S}^\omega},\nonumber
\end{eqnarray}
with
\begin{equation}
    \label{ec.001}
d^{\widehat{S}^\omega}:=\widehat{\Upsilon}_{\mathfrak{qg}^\#} \circ \widehat{S}^\omega \circ \widehat{\Upsilon}^{-1}_{\mathfrak{qg}^\#}.
\end{equation}

\begin{Theorem}
\label{6.1.5}
Assume that the operator $$d^{\widehat{DS}^\omega}$$ admits a formally adjoint operator $$d^{\widehat{DS}^\omega\star} $$
with respect to $\langle-|-\rangle^\bullet_\r$ for every qpc $\omega$. Then, a qpc $\omega$ is a critical point of the non--commutative geometrical Yang--Mills functional if and only if
\begin{equation}
    \label{4.3.f29.1}
        d^{\widehat{DS}^\omega\star}\,(\widehat{\Upsilon}_{\mathfrak{qg}^\#}(R^{\omega}))=0.
\end{equation}
\end{Theorem}
 
\begin{proof}
According to the proof of Theorem 4.7 of reference \cite{saldym}, a qpc $\omega$ is a critical point of $\qS_\YM$ if and only if\footnote{It is worth mentioning that, in Theorem~4.7 of reference \cite{saldym}, we treat the left and right module structures of $\Mor(\ad,\Delta_\Hor)$ simultaneously, whereas in the present paper we focus exclusively on the right module structure.}
\begin{equation*}
    \mathrm{Re}\left( \langle d^{\widehat{DS}^\omega\star}\,( \widehat{\Upsilon}_{\mathfrak{qg}^\#}(R^{\omega}))\mid \widehat{\Upsilon}_{\mathfrak{qg}^\#}(\lambda)\rangle^\bullet_\r\right)=0
\end{equation*}
for every $\lambda$ $\in$ $\overrightarrow{\mathfrak{qpc}(\zeta)}$ (see equation  (\ref{2.f24.2})), where $\mathrm{Re}(w)$ denotes the real part of the complex number $w$ $\in$ $\C$.

Using equation (\ref{inersperado}), it is easy to check that the twisted covariant derivative $\widehat{DS}^\omega$ is a first--order operator on the direct sum of the spaces $$\{\tau:\mathfrak{qg}^\# \longrightarrow \Hor^k P\mid \tau \mbox{ is linear},\;\,(\tau\otimes \id_G)\circ \ad=\Delta_\Hor \circ \tau \;\mbox{ and }\; \widehat{\tau}=(-1)^{k+1}\tau\}$$ and it follows from (\ref{twcoder1}) that $\lambda$ and $R^\omega$ are elements of this direct sum, for every $\lambda$ $\in$ $\overrightarrow{\mathfrak{qpc}(\zeta)}$. Hence $$\widehat{\Upsilon}^{-1}_{\mathfrak{qg}^\#}(d^{\widehat{DS}^\omega\star}\,( \widehat{\Upsilon}_{\mathfrak{qg}^\#}(R^{\omega})))\;\in\; \overrightarrow{\mathfrak{qpc}(\zeta)}$$ and due to the fact that $\langle-|-\rangle^\bullet_\r$ is an inner product, it follows that $\omega$ is a critical point of $\qS_\YM$ if and only if
equation (\ref{4.3.f29.1}) holds.
\end{proof}


We will refer to equation (\ref{4.3.f29.1}) as the {\it non--commutative geometrical Yang--Mills equation} and solutions of this equation will be called {\it Yang--Mills quantum principal connections}.

\begin{Remark}
    \label{assuptions}
   From now on, we shall restrict our study exclusively to  quantum principal bundles $\zeta_\Xi$ as in Section 3.2 such that $H=SU(N)$.
\end{Remark}

The motivation for the condition on $H$ is that we aim to develop a {\it physical} Yang--Mills theory with qpb's; so considering $SU(N)$ is appropriate. There are reasons to believe that space--time at the Planck scale could be a non--commutative space (\cite{con}); however, as far as the author know, there are no reasons to assume that the structure groups of physical theories must be a purely non--commutative geometrical object. Nevertheless, from a purely mathematical point of view,  Definition \ref{yangmills1} and equation (\ref{4.3.f29.1}) (when $d^{\widehat{DS}^\omega}$ is formally adjointable) hold even if the structure group is a purely non--commutative geometrical object. 

According to reference \cite{lan2}, the qpb's of Section 3.2 satisfy Remark \ref{remadirac} and the regularity condition ($X^2_\Xi=X_\Xi$). This is because the Connes’ differential form spaces is isomorphic to
the graded differential $\ast$--algebras of equation (\ref{ec.012}). Furthermore, there exists $n$ $\in$ $\N$ such that $\Omega^{n+l}(X_\Xi)=0$ for all $l$ $\in$ $\N$; and there exists an element $\dvol$ $\in$ $\Omega^n(X_\Xi)$ that satisfies 
\begin{equation}
    \label{ec.new.1}
    \Omega^{n}(X_\Xi)=X_\Xi\,\dvol,\qquad \dvol\,x=x\,\dvol,\qquad x\,\dvol=0 \,\Longleftrightarrow x=0
\end{equation}
for all $x$ $\in$ $X_\Xi$ (\cite{lan2}).  Thus, every element of $\mu$ $\in$ $\Omega^n(X_\Xi)$ is of the form $L_\Xi(f)\,\dvol$ for some $f$ $\in$ $C^\infty_\C(X)$. So   (\cite{lan2,lan4})  
\begin{equation}
    \label{teodio4}
    \int_{X_\Xi} \mu:=\mathrm{Tr}_w(L_\Xi(f)\, |\D|^{-1})=\mathrm{Tr}_w(f\, |\D|^{-1})=\int_X f\,\dvol,
\end{equation}
where the integral of the right--hand side of the last equation is the usual (normalized) integral on the manifold $X$. Moreover, this implies that
\begin{equation}
    \label{ec.new.2}
    \int_{X_\Xi}  d\mu=0
\end{equation}
for all $\mu$ $\in$ $\Omega^{n-1}(X_\Xi)$. Finally, there exists a quantum star Hodge operator (\cite{lan2,lan4}) 
\begin{equation}
    \label{ec.new.3}
    \star_q: \Omega^n(X_\Xi)\longrightarrow \Omega^{n-k}(X_\Xi)
\end{equation}
such that 
$$\star_q(\mathbbm{1})=\dvol,\,\qquad \star_q(\dvol)=\mathbbm{1}$$ and 
\begin{equation}
    \label{ec.new.4}
    \langle \mu_1\mid \mu_2\rangle_\r=\int_{X_\Xi} \mu^\ast_1\, (\star_q\mu_2).
\end{equation}




In our case, according to Section 2.2, the quantum dual Lie algebra $\mathfrak{qg}^\#$ is exactly $$\mathfrak{h}^\#_\C=\mathfrak{su}^\#_\C(N),$$ the complexification of the dual space of the Lie algebra  $\mathfrak{h}=\mathfrak{su}(N)$ of $SU(N)$. Moreover, in accordance with equation (\ref{adstructure1}), the $\ad$ $G$--corepresentation of equation (\ref{2.f15}) is the pull--back of the complexification of the (right) adjoint action of $H$ on $\mathfrak{h}$; and hence, $\ad$ $\in$ $\Obj(\Rep_\G)$. This implies that we can consider the associated right qvb $$E^{\mathfrak{h}^\#_\C}_\r=\Mor(\ad,\Delta_P) $$ and the inner product 
\begin{equation}
    \label{lonecesito}
    \langle-|-\rangle^\bullet_\r: (E^{\mathfrak{h}^\#_\C}_\r\otimes_{X_\Xi}\Omega^2(X_\Xi))\times (E^{\mathfrak{h}^\#_\C}_\r\otimes_{X_\Xi}\Omega^2(X_\Xi))\longrightarrow \C 
\end{equation}
of the equation (\ref{4.3.f26}).

According to Theorem 27 of reference \cite{saldym}, since there exists a quantum Hodge operator $\star_q$, for $\delta^V$ $\in$ $\Obj(\Rep_\G)$, the exterior covariant derivative $d^{\widehat{\nabla}^\omega_V}$ of the gauge qlc $\widehat{\nabla}^\omega_V$ is formally adjointable with respect to the inner product of equation (\ref{4.3.f26}) for every qpc $\omega$. In concrete, the operator
\begin{equation}
    \label{formaldiff}  
    d^{\widehat{\nabla}^{\omega}_{V}\star} :E^V_\r\otimes_B \Omega^{a+1}(B)\longrightarrow E^V_\r\otimes_B \Omega^a(B)
\end{equation}
given by
\begin{equation*}
d^{\widehat{\nabla}^{\omega}_{V}\star}:=(-1)^{a} \,(\id_{E^V_\r}\otimes_B \star^{-1}_q) \circ\; d^{\widehat{\nabla}^{\omega}_V}\circ (\id_{E^V_\r}\otimes_B \star_q)
\end{equation*}
is the formally adjoint operator of $d^{\widehat{\nabla}^\omega_V}$. This holds, in particular, for the corepresentation $\ad$ $\in$ $\Obj(\Rep_\G)$. It is worth mentioning that the sing $(-)^a$ in the definition of $d^{\widehat{\nabla}^{\omega}_{V}\star}$ is because the differential is anti--Hermitian. When the differential is Hermitian, the sign is $(-1)^{a+1}$ (\cite{saldym}). 

Equation (\ref{formaldiff}) implies that the operator $d^{\widehat{DS}^\omega}$ admits a formally adjoint operator if and only the operator $d^{\widehat{S}^\omega}$ admits a formally adjoint operator. In this way, to conclude this section, we shall prove that this formal adjoint operator always exists.

Let $N\geq 1$  and let $\langle-|-\rangle_{\mathfrak{h}_\C}$ be the negative of the canonical sesquilinear extension in $\mathfrak{h}_\C=\mathfrak{h}\otimes \C$ of the Killing form of $\mathfrak{h}$. Since the $SU(N)$--action $\ad^{\mathrm{class}}$ (see equation (\ref{adjointclassical})) is irreducible and unitary with respect to $\langle-|-\rangle_{\mathfrak{h}_\C}$ and $\ad$ is its pull--back, it follows that $\ad$ is irreducible and unitary with respect to the canonical inner product $\langle-|-\rangle_{\mathfrak{h}^\#_\C}$ in $\mathfrak{h}^\#_\C$ induced by $\langle-|-\rangle_{\mathfrak{h}_\C}$. In addition, if $\{ v_k\}^m_{k=1}$ is an orthonormal basis of $\mathfrak{h}$ with respect to the negative of the Killing form ($m=N^2-1$), then, for its dual basis $\{\widehat{\theta}_k \}^m_{k=1}$ on $\mathfrak{h}^\#_\C$, we have
\begin{equation}
    \label{eq.0.1}
    \langle \theta_l\mid \theta_j\rangle_{\mathfrak{h}^\#_\C}=\delta_{lj},\qquad \theta^\ast_k=-\theta_k\qquad \mbox{ with } \qquad \theta_k=i\,\widehat{\theta}_k. 
\end{equation}
From now on, we will work with the inner product $\langle-|-\rangle_{\mathfrak{h}^\#_\C}$ and the basis $\beta_N=\{\theta_k \}^m_{k=1}$.  We have decided to work with the basis $\beta_N$ instead of the basis $\{\widehat{\theta}_k \}^m_{k=1}$ because in physics, the basis used is always anti--Hermitian.

As we commented on Section 3, when a \emph{classical} Lie group equipped with its \emph{canonical} differential calculus is used in a qpb, there is a canonical embedded differential (see Definition \ref{embbededdifferential}), which is given by  $$\Theta=-{1\over 2}\mathrm{c}^T.$$ In addition, in this case, $\mathrm{c}^T$ is the pull--back of the {\it classical} adjoint Lie algebra representation $\mathrm{c}$ (see equations (\ref{adstructure2}), (\ref{adstructure4})).

For $N=1$, we have $H=SU(N)=U(1)$ and $\mathrm{c}=0$; so $\mathrm{c}^T=0$. This implies that the operator $\widehat{S}^\omega$ is identically zero  (see Definition \ref{Soperator}) and therefore $$\widehat{DS}^\omega=\widehat{D}^\omega.$$  By equation (\ref{3.f10.6}) we obtain $$d^{\widehat{DS}^\omega}=d^{\widehat{\nabla}^\omega_{\mathfrak{h}^\#_\C}},$$ which is formally adjointable. This proves that equation (\ref{4.3.f29.1}) makes sense in the context of Remark \ref{assuptions} for $N=1$, and it takes the form 
\begin{equation}
    \label{4.3.f29.2}
    d^{\widehat{\nabla}^\omega_{\mathfrak{h}^\#_\C}\star}\,(\widehat{\Upsilon}_{\mathfrak{h}^\#_\C}(R^{\omega}))=0.
\end{equation}

In this particular case, the non--commutative geometrical Yang--Mills equation is completely similar to the classical one (see equation (\ref{classicalymequation})). Furthermore, in this particular case, the only difference between them is the fact that one operates in the {\it classical space--time} modeled by $X$, and the other one operates in the {\it quantum space--time} modeled by $X_\Xi$.

For $N\geq 2,$ notice that (see equation (\ref{adstructure2})) $$\mathrm{c}(v_l,v_j)=[v_l,v_j]_\C=[v_l,v_j]=\sum^m_{k=1}f_{ljk}\,v_k $$ for some $f_{ljk}$ $\in$ $\R$. The numbers $f_{ljk}$ are usually called {\it structure constants}. According to \cite{diff1},  the structure constants satisfy
\begin{equation}
    \label{lieprop2}
    f_{ljk}=-f_{jlk},\qquad f_{ljk}=f_{klj}.
\end{equation}
In this way, by equation  (\ref{adstructure4}) we get 
$$\mathrm{c}^T(\theta_k)(v_l,v_j)=i\,(\theta_k \circ \mathrm{c})(v_l,v_j)=-\widehat{\theta}_k([v_l,v_j])=-\f_{ljk};$$ and we conclude that
\begin{equation}
    \label{classicalcommutator}
\mathrm{c}^T(\theta_k)=\sum^m_{l,j}-f_{ljk}\,\widehat{\theta}_l\otimes \widehat{\theta}_j=\sum^m_{l,j}\,f_{ljk}\,\theta_l\otimes \theta_j.
\end{equation} 

On the other hand, under the isomorphism $\widehat{\Upsilon}_{\mathfrak{h}^\#_\C}$ of equation (\ref{3.f7.5}), the inner product of equation  (\ref{lonecesito}) induces an inner product on $\Mor(\ad,\Delta_\Hor)$  
\begin{equation}
    \label{lonecesito2.1}
    \langle-|-\rangle^\bullet_\Hor:  \Mor(\ad,\Delta_\Hor)\times \Mor(\ad,\Delta_\Hor)\longrightarrow \C
\end{equation}
given by $$\langle \tau_1\mid \tau_2\rangle^\bullet_\Hor:=\langle \widehat{\Upsilon}_{\mathfrak{h}^\#_\C}(\tau_1)\mid \widehat{\Upsilon}_{\mathfrak{h}^\#_\C}(\tau_2)\rangle^\bullet_\r =\sum^m_{k=1} \int_{X_\Xi} \tau_1(\theta_k)^\ast\,\tau_2(\theta_k)=\sum^m_{k=1} \int_{X_\Xi} \tau_1(\theta_k)^\ast\,\star'_q(\tau_2(\theta_k)),$$ where $\star'_q$ is the extension on $\Mor(\ad,\Delta_\Hor)$ of the quantum Hodge of the operator $\star_q$ given by $$\star'_q:=\widehat{\Upsilon}^{-1}_{\mathfrak{h}^\#_\C}\circ (\id\otimes_{X_\Xi}\star_q) \circ \widehat{\Upsilon}_{\mathfrak{h}^\#_\C}.$$ For example, by equations (\ref{generators}),  (\ref{3.f5.2}), we obtain that
\begin{equation}
      \label{porqueno}
      \star'_q(T\,\mu)=T\,\star_q(\mu)
\end{equation}
for $T$ $\in$ $E^{\mathfrak{h}^\#_\C}_\r$, $\mu$ $\in$ $\Omega^\bullet(X_\Xi)$.

In this way,  $d^{\widehat{S}^\omega}$ admits a formally adjoint operator with respect to $\langle-|-\rangle^\bullet_\r$ if and only if $\widehat{S}^\omega$ admits a formally adjoint operator with respect to $\langle-|-\rangle^\bullet_\Hor$.

\begin{Theorem}
    \label{soperatoradjunto}
    For $N\geq 2$, the operator $\widehat{S}^{\omega} $ is formally adjointable with respect to $\langle-|-\rangle^\bullet_\Hor$ for every qpc $\omega$.
\end{Theorem}
\begin{proof}
Let $\omega$ be a qpc and let $\tau$ $\in$ $\Mor(\ad, \Delta_\Hor)$ with $\Im(\tau)\subseteq \Hor^a \,P$. Then 
\begin{eqnarray*}
    S^\omega(\tau)&=&\langle\omega,\tau\rangle-(-1)^a\langle\tau,\omega\rangle-(-1)^a[\tau,\omega]
    \\
    &=&\langle\omega,\tau\rangle-(-1)^a\langle\tau,\omega\rangle+2\,(-1)^a  \langle\tau,\omega\rangle 
    \\
    &=&
    \langle\omega,\tau\rangle+(-1)^a\langle\tau,\omega\rangle.
\end{eqnarray*}
According to Section 3, $\mathfrak{qpc}(\zeta)$ is an affine space modeled by $\overrightarrow{\mathfrak{qpc}(\zeta)}$. Hence (see Proposition \ref{teodio1}) $$\omega=\omega^c+\lambda$$ for some $\lambda$ $\in$ $\overrightarrow{\mathfrak{qpc}(\zeta)}$. This implies that $$S^\omega=S^{\omega^c}+S^\lambda,$$ where $$S^\lambda:\Mor(\ad,\Delta_\Hor)\longrightarrow \Mor(\ad,\Delta_\Hor) $$ is given by $$S^\lambda(\tau)=\langle\lambda,\tau\rangle+(-1)^a\langle\tau,\lambda\rangle.$$ However, $\omega^c$ is a regular qpc; so $S^{\omega^c}=0$ (\cite{micho2}); which implies that $$S^\omega=S^\lambda.$$ Furthermore, using equation (\ref{inersperado}), we get 
\begin{eqnarray*}
\widehat{S}^\lambda(\tau)=\widehat{\langle\lambda,\widehat{\tau}\rangle}+(-1)^a\widehat{\langle\widehat{\tau},\lambda\rangle}=(-1)^a\langle \tau,\widehat{\lambda}\rangle+\langle \widehat{\lambda},\tau\rangle
=S^\lambda(\tau)
\end{eqnarray*}
and we conclude that $S^\lambda=\widehat{S}^\lambda$.

Notice that (\cite{stheve}) $$ \tau\longmapsto \langle\lambda,\tau\rangle \qquad \mbox{ and }\qquad \tau\longmapsto \langle\tau,\lambda\rangle$$ are operators on $\Mor(\ad,\Delta_\Hor)$ and  in this way, for  $\theta_k$ $\in$ $\beta_N$  we get  
\begin{eqnarray*}
    \langle \lambda,\tau\rangle^\ast(\theta_k)=-{1\over 2}\sum^m_{l,j=1}(f_{ljk}\,\lambda(\theta_l)\,\tau(\theta_j))^\ast&=&-{1\over 2}(-1)^a\sum^m_{l,j=1} f_{ljk}\, \tau(\theta_j)^\ast\,\lambda(\theta_l)^\ast
    \\
    &=&
    -{1\over 2}(-1)^a\sum^m_{l,j=1} f_{ljk}\, \tau(\theta_j)^\ast\,\lambda(\theta^\ast_l)
    \\
    &=&
    {1\over 2}(-1)^a\sum^m_{l,j=1} f_{ljk}\, \tau(\theta_j)^\ast\,\lambda(\theta_l).
\end{eqnarray*}
Let $\sigma$ $\in$ $\Mor(\ad,\Delta_\Hor)$ with  $\Im(\sigma)\subseteq \Hor^{a+1} \,P$. Then, by equation (\ref{lieprop2}) we obtain
\begin{eqnarray*}
    \langle \langle\lambda,\tau\rangle\mid \sigma\rangle^\bullet_\Hor&=&\sum^m_{k=1}\int_{X_\Xi} \langle\lambda,\tau\rangle^\ast(\theta_k)\,\star'_q(\sigma(\theta_k)) 
    \\
    &=&
    {1\over 2}(-1)^a\sum^m_{l,j,k=}f_{ljk}\int_{X_\Xi}\tau(\theta_j)^\ast\,\lambda(\theta_l)\,\star'_q(\sigma(\theta_k))
    \\
    &=&
    {1\over 2}(-1)^a\sum^m_{l,j,k=}f_{lkj}\int_{X_\Xi}\tau(\theta_k)^\ast\,\lambda(\theta_l)\,\star'_q(\sigma(\theta_j))
    \\
    &=&
    -{1\over 2}(-1)^a\sum^m_{l,j,k=}f_{ljk}\int_{X_\Xi}\tau(\theta_k)^\ast\,\lambda(\theta_l)\,\star'_q(\sigma(\theta_j))
    \\
    &=&
\sum^m_{k=1}\int_{X_\Xi}\tau(\theta_k)^\ast\,\star'_q\left(\star'^{-1}_q\left(-{1\over 2}(-1)^a\sum^m_{l,j=1}f_{ljk}\, \lambda(\theta_l)\,\star'_q(\sigma(\theta_j))\right)\right)
    \\
    &=&
\sum^m_{k=1}\int_{X_\Xi}\tau(\theta_k)^\ast\,\star'_q\left((-1)^a\star'^{-1}_q\left( \langle \lambda,\star'_q(\sigma)\rangle(\theta_k) \right)\right)
    \\
    &=&
    \langle \tau\mid (-1)^a \star'^{-1}_q\left( \langle \lambda,\star'_q(\sigma)\rangle\right)\rangle^\bullet_\Hor.
\end{eqnarray*}
Therefore, the formal adjoint operator of $$\tau\longmapsto \langle\lambda,\tau\rangle$$ with respect to $\langle-|-\rangle^\bullet_\Hor$ is
\begin{equation}
    \label{teodio2}
    \sigma \longmapsto (-1)^a  \star'^{-1}_q ( \langle \lambda,\star'_q(\sigma)\rangle).
\end{equation}
On the other hand, for $\theta_k$ $\in$ $\beta_N$  we get  
\begin{eqnarray*}
    \langle \tau,\lambda\rangle^\ast(\theta_k)=-{1\over 2}\sum^m_{l,j=1}(f_{ljk}\,\tau(\theta_l)\,\lambda(\theta_j))^\ast&=&-{1\over 2}(-1)^a\sum^m_{l,j=1} f_{ljk}\, \lambda(\theta_j)^\ast\,\tau(\theta_l)^\ast
    \\
    &=&
    -{1\over 2}(-1)^a\sum^m_{l,j=1} f_{ljk}\, \lambda(\theta^\ast_j)\,\tau(\theta_l)^\ast
    \\
    &=&
    {1\over 2}(-1)^a\sum^m_{l,j=1} f_{ljk}\, \lambda(\theta_j)\,\tau(\theta_l)^\ast.
\end{eqnarray*}
Thus, by equation (\ref{lieprop2}) and the fact that the quantum integral is the usual integral on $X$ (see equation (\ref{teodio4})), we obtain
\begin{eqnarray*}
    (-1)^a\langle \langle\tau,\lambda\rangle\mid \sigma\rangle^\bullet_\Hor&=& (-1)^a\sum^m_{k=1}\int_{X_\Xi}\langle \tau,\lambda\rangle^\ast(\theta_k)\,\star'_q(\sigma(\theta_k))
    \\
    &=&
    {1\over 2}\sum^m_{l,j,k=1}f_{ljk}\int_{X_\Xi}\lambda(\theta_j)\,\tau(\theta_l)^\ast\,\star'_q(\sigma(\theta_k))
    \\
    &=&
    {1\over 2}(-1)^{n-1}\sum^m_{l,j,k=1}f_{ljk}\int_{X_\Xi}\tau(\theta_l)^\ast\,\star'_q(\sigma(\theta_k))\,\lambda(\theta_j)
    \\
    &=&
    {1\over 2}(-1)^{n-1}\sum^m_{l,j,k=1}f_{kjl}\int_{X_\Xi}\tau(\theta_k)^\ast\,\star'_q(\sigma(\theta_l))\,\lambda(\theta_j)
    \\
    &=&
    -{1\over 2}(-1)^{n-1}\sum^m_{l,j,k=1}f_{ljk}\int_{X_\Xi}\tau(\theta_k)^\ast\,\star'_q(\sigma(\theta_l))\,\lambda(\theta_j)
    \\
    &=&
    \sum^m_{k=1} \int_{X_\Xi} \tau(\theta_k)^\ast\,\star'_q\left(\star'^{-1}_q\left(-{1\over 2}(-1)^{n-1}\sum^m_{l,j=1}f_{ljk}\star'_q(\sigma(\theta_l))\,\lambda(\theta_j) \right)\right)
    \\
    &=&
    \sum^m_{k=1} \int_{X_\Xi} \tau(\theta_k)^\ast\,\star'_q\left((-1)^{n-1}\star'^{-1}_q\left(\langle \star'_q(\sigma),\lambda\rangle(\theta_k) \right)\right)
    \\
    &=&
    \langle \tau\mid (-1)^{n-1}\star'^{-1}_q\left(\langle \star'_q(\sigma),\lambda \rangle\right)\rangle^\bullet_\Hor.
\end{eqnarray*}
Hence, the formal adjoint operator of $$\tau\longmapsto \langle \tau,\lambda\rangle $$ with respect to $\langle-|-\rangle^\bullet_\Hor$ is $$\sigma\longmapsto (-1)^{n-1}\star'_q\left(\langle\star'_q(\sigma),\lambda\rangle\right). $$ It follows from this result and equation (\ref{teodio2}) that the operator
\begin{equation}
    \label{teodio5}
\widehat{S}^{\omega\star}:\Mor(\ad,\Delta_\Hor)\longrightarrow \Mor(\ad,\Delta_\Hor)
\end{equation}
given by $$\widehat{S}^{\omega}(\sigma)=(-1)^{a}\star'^{-1}_q(\widehat{S}^{\omega}(\star'_q(\sigma))) $$ for all $\sigma$ $\in$ $\Mor(\ad,\Delta_\Hor)$ with $\Im(\sigma) \subseteq \Hor^{a+1}\,Y_\Xi$
is the formal adjoint operator of $\widehat{S}^{\omega}$ with respect to $\langle-|-\rangle^\bullet_\Hor$.
\end{proof}

In this way,  we have that the operator $$d^{\widehat{DS}^\omega}:=\widehat{\Upsilon}_{\mathfrak{h}^\#_\C} \circ \widehat{DS}^\omega \circ \widehat{\Upsilon}^{-1}_{\mathfrak{h}^\#_\C}=\widehat{\Upsilon}_{\mathfrak{h}^\#_\C} \circ (\widehat{D}^\omega-\widehat{S}^\omega) \circ \widehat{\Upsilon}^{-1}_{\mathfrak{h}^\#_\C}=d^{\widehat{\nabla}^\omega_{\mathfrak{h}^\#_\C}}-d^{\widehat{S}^{\omega}}$$ admits a formally adjoint operator with respect to $\langle-|-\rangle^\bullet_\r$ and therefore, the non--commutative geometrical Yang--Mills equation (equation (\ref{4.3.f29.1})) makes sense in the context of Remark \ref{assuptions} for $N\geq 1$:
\begin{equation}
    \label{ymgeometric}
    (d^{\widehat{\nabla}^\omega_{\mathfrak{h}^\#_\C}\star}-d^{\widehat{S}^{\omega}\star})\,(\widehat{\Upsilon}_{\mathfrak{h}^\#_\C}(R^{\omega}))=0.
\end{equation}
Recall that, for $N=1$, we have $d^{S^{\omega}\star}=0$ and therefore the last equation and equation (\ref{4.3.f29.2}) are the same. 

\begin{Remark}
\label{comparar}
    In the {\it classical} case, equation (\ref{ymgeometric}) matches with dualization (and complexification) of equation (\ref{classicalymequation}). In fact, in the {\it classical} case, every quantum principal connection is regular (see Section 3), hence $S^\omega=0$; which implies that $d^{\widehat{S}^{\omega}}=d^{\widehat{S}^{\omega}\star}=0$ and thus, equation (\ref{ymgeometric}) takes the form 
$$d^{\widehat{\nabla}^\omega_{\mathfrak{h}^\#_\C}\star}(\widehat{\Upsilon}_{\mathfrak{h}^\#_\C}(R^{\omega}))=0.$$ 
On the other hand, the definition of $R^\omega$ has been formulated exactly to be the pull--back (and complexification) of the the curvature $\Omega^\omega$ of a principal connection in differential geometry with $-\displaystyle{1\over 2}c^T$ as the embedded differential used (\cite{nodg}). Since the non--commutative geometrical counterpart (by dualizing with the pull--back) of the space of basic differential forms of type $\ad^{\mathrm{class}}$ is exactly the space $\Mor(\ad,\Delta_\Hor)$ and $\widehat{\Upsilon}^{-1}_{\mathfrak{h}^\#_\C}$ is only the product of elements of $\Mor(\ad,\Delta_P)$ with $\Omega^\bullet(B)$, it follows that $\widehat{\Upsilon}_{\mathfrak{h}^\#_\C}$ coincides with the pull--back (and complexification) of the map $GP$, the canonical isomorphism between basic differential forms of type $\ad^{\mathrm{class}}$ and vector bundle--valued differential forms (\cite{nodg}). 

Finally, the definition of  $d^{\widehat{\nabla}^\omega_{\mathfrak{h}^\#_\C}\star}$ coincides with the definition of $d^{\nabla^\omega_{\mathfrak{h}}\star}$ (\cite{nodg}); therefore, equation (\ref{ymgeometric}) matches with dualization of equation (\ref{classicalymequation}). A more detail explanation is given in Appendix C of reference \cite{saldym}.

On the other hand, in \cite{con} A. Connes proves that his functional also reproduces the classical case. Thus, we conclude that, in the context of differential geometry, both the non--commutative analytical Yang--Mills equation and the non--commutative geometrical Yang--Mills equation coincide with the classical Yang--Mills equation, as expected. 
\end{Remark}

Clearly, the novelty of the theory of Section 5.3 lies in the formulation of the non--commutative geometrical Yang--Mills equation. Since equations (\ref{4.3.f29}), (\ref{ymgeometric}) are proper generalizations of the {\it classical} Yang--Mills equation, the following questions naturally arise in the context of non--commutative geometry:

\begin{enumerate}
\item Do both equations lead to the same result?
\item If not, do their solution spaces share any common features?
\end{enumerate}

These questions will be addressed in the next section through a concrete example. From a physical perspective, the Lie groups $$U(1),\qquad SU(2),\qquad  SU(3)$$ are particularly important, as they correspond to the structure groups of the principal bundles that describe three of the four fundamental interactions: electromagnetism, the weak nuclear interaction, and the strong nuclear interaction.
\begin{Remark}
    \label{gaugeboson}
    In differential geometry, given $\pi_X: Y\longrightarrow X $  a principal $G$--bundle over a manifold $X$, solutions of the Yang--Mills equation represent gauge fields in the vacuum moving through the space--time $X$, for the Lie groups $G=U(1)$, $SU(2)$, $SU(3)$.
    
    Following this interpretations, given  $\zeta=(P,B,\Delta_P) $  a quantum principal $\G$--bundle over a quantum space $B$ with a differential calculus, solutions of the non--commutative geometrical Yang--Mills equation represent non--commutative gauge boson fields in the vacuum  moving through the quantum space--time $B$ for $\G$ the canonical quantum group associated to the Lie group $G$, where $G=U(1)$, $SU(2)$, $SU(3)$ and with their canonical differential structures.
\end{Remark}

To develop a concrete example, we will focus on the electromagnetic interaction, i.e., a Yang--Mills theory with the Lie group $U(1)$. It is worth mentioning that the study of quantum principal $U(1)$--bundles is by no means a minor topic. Indeed, it continues to be an active area of research, as the reader can verify from~\cite{braming}.

Examples of the theory for the remaining Lie groups of the standard model will be addressed in forthcoming publications.

\section{Example: Electromagnetism in the Quantum $n$--Torus}

Throughout the paper, we have made four non--trivial assumptions. The only general assumption appears in Remark \ref{rema}, where we assumed that the base space $B$ is a $\ast$--(sub)algebra stable under functional calculus. The other tree assumptions are: Remark \ref{remadirac}, the regularity condition $B=B^2$ and Remark \ref{assuptions}. The regularity condition is necessary in order to define the inner product of equation (\ref{3.3.f18}), Remark \ref{remadirac} is required to obtain the non--commutative geometrical Dirac equation (equation (\ref{3.3.f19})); while Remark \ref{assuptions} is a particular case of Remark \ref{remadirac} and it is necessary to derive the non--commutative geometrical Yang--Mills equation (equation (\ref{ymgeometric})).

 It is worth mentioning that for toric non--commutative manifolds (for example, the quantum base space of a qpb under Remark \ref{assuptions}), the non--commutative analytical Yang--Mills equation (equation (\ref{4.3.f29})) takes the form (\cite{lan2,lan4})
\begin{equation}
    \label{ymanalitic}
    [\widehat{\nabla},(\id_{E^V_\r}\otimes_{B} \star_q)\circ R^{\widehat{\nabla}}]=0.
\end{equation}
Therefore, quantum principal $SU(N)$--bundles over toric non--commutative manifolds provide a natural class of qpb's for studying the non--commutative analytical Yang--Mills equation and the non--commutative geometric Yang--Mills equation simultaneously, with the aim of addressing the two questions posed above.

In this section we are going to use a non--commutative $n$--torus as the quantum base space of a qpb to illustrate the theory presented in the last section; putting particular emphasis in the theory showed in Section 5.3. We want to use a non--commutative $n$--torus because this space has been studied for a long time in detail, for example in \cite{con,dabro1,dabro2,lap} and probably the reader is familiar with some results about this quantum space. So, it helps us to highlight the differences and similarities of both non--commutative Yang--Mills functionals. 

Additionally, in this section, we will see that it is possible to change the {\it classical} $\ast$--FODC of the Lie group $H=SU(N)$ to a non--commutative one in our framework. However, in this situation, in general, there is no a {\it canonical} embedded differential and therefore, it is quite challenging to prove Theorem  \ref{soperatoradjunto} in this context. 

\subsection{A quantum principal $\U(1)$--bundle over the Quantum $n$--Torus.}

Consider $$\torus^n=U(1)\times \cdots \times U(1)$$ an $n$--torus for $n\geq 2$. Then $\torus^n$ is a closed, oriented, spin Riemannian manifold of dimension $n$. We will consider  $\torus^n$  with the Riemannian metric defined by the product Riemannian metric of each $U(1)$ (the flat $n$--torus). Then, the group product of $\torus^n$ is an isometrical action of $\torus^n$ on itself and hence, $\torus^n$ is in the domain of the deformation map $L_\Xi$. In this way, we define  the quantum $n$--torus as the non--commutative (unital) $\ast$--algebra (\cite{dabro1}) $$\torus^n_\Xi:=L_\Xi(C^\infty_\C(\torus^n)).$$ A generic element of $\torus^n_\Xi$ is  a formal sum $$\sum a_{m_1\cdots m_n} \u^{m_1}_1\,\u^{m_2}_2\cdots \u^{m_n}_n,$$ where $\{ a_{m_1\cdots m_n}\}_{m_i\in \Z}$ form a rapid decay sequence in $\C$ (decay faster than the inverse of any polynomial in $(m_1,...,m_n)$), and the operators  $\u_1$,..., $\u_n$ satisfy
\begin{equation}
    \label{eq.6.1}
    \u^\ast_k=\u^{-1}_k,\qquad \u_k\,\u_j=\mathrm{e}^{2\pi i\,\Xi_{kj}}\,\u_j\,\u_k
\end{equation}
for all $k$, $j$ $\in$ $\{1,...,n\}$, with $\Xi=(\Xi_{kj})$ an antisymmetric $n$--dimensional square matrix \cite{con}. By taking the space $C(\torus^n_\Xi)$ defined as the universal $C^\ast$--algebra of $\torus^n_\Xi$, it is well--known that $\torus^n_\Xi$  is stable under functional calculus \cite{con}. 

The next step is to consider the canonical spectral triple of $\torus^n_\Xi$ and with that, the Connes' differential form space. In fact, let $$\tau_0:\torus^n_\Xi\longrightarrow \C$$ be the unique normalized tracial state of $\torus^n_\Xi$  
\begin{equation}
    \label{eq.6.1.1}
    \tau_0\left(\sum a_{m_1\cdots m_n} \u^{m_1}_1\,\u^{m_2}_2\cdots \u^{m_n}_n\right)=a_{0,...,0}.
\end{equation}
and consider the GNS Hilbert space $ L^2(\torus^n_\Xi,\tau_0)$. In this way, 
\begin{equation}
    \label{eq.6.1.2}
(\torus^n_\Xi,\,\mathfrak{H}:=L^2(\torus^n_\Xi,\tau_0)\otimes \C^{2^{[n/2]}},\D),
\end{equation}
where $[n/2]$ denotes the integer part of $n/2$, is the canonical $n$--dimensional spectral triple. Here (\cite{con}),  
\begin{equation}
    \label{eq.6.1.3}
     \D:=\sum^n_{j=1} \gamma^j\,\delta_j,
\end{equation}
where $\gamma^j$ are the square $2^{[n/2]}$--dimensional gamma matrices (\cite{diff1}), and 
\begin{equation}
    \label{eq.6.1.4}
      \delta_j: \torus^n_\Xi\longrightarrow \torus^n_\Xi
\end{equation}
is the $\ast$--preserving derivation in $\torus^n_\Xi$ defined by $$\delta_j(\mathbbm{u}_k)=\left\{
\begin{array}{ll}
 2\pi i\, \mathbbm{u}_j & \text{ if } j=k \\
  \,\,0 & \text{ if } j \not= 0.
\end{array}
\right.$$

Associated to the canonical $n$--dimensional spectral triple, we can take the Connes' differential form space  $(\Omega^\bullet_\D(\torus^n_\Xi),d_\D,\ast)$  and according to references \cite{lan4,lan2}, the previous graded differential $\ast$--algebra is isomorphic to (see Section 3.2) $$(\Omega^\bullet(\torus^n_\Xi):=L_\Xi(\Omega^\bullet_\C(\torus^n)),d=i\,d',\ast).$$ Since the derivation $\delta_j$ is the non--commutative counterpart of the operator $\displaystyle {\partial \over \partial x_j }$ in the {\it classical} case (\cite{con,dabro1,dabro2}), the previous graded differential $\ast$--algebra can be seen as follows: consider the space $\C^n$ and its exterior algebra $$\bigwedge^\bullet \C^{n}=\bigoplus^n_{k=0}\wedge^k \C^{n}\qquad \mbox{ and } \qquad \Omega^\bullet(\torus^n_\Xi)\cong \torus^n_\Xi\otimes \bigwedge^\bullet \C^{n}$$ with the $\ast$--algebra structure given by $$(x_1\otimes v_1)\cdot (x_2\otimes v_2)=x_1\cdot x_2\otimes v_1\wedge v_2,\qquad (x\otimes v)^\ast=x^\ast\otimes v^\ast $$ and the differential 
\begin{equation}
    \label{lonecesito1}
   d: \Omega^\bullet(\torus^n_\Xi)\longrightarrow \Omega^\bullet(\torus^n_\Xi),\qquad x\longmapsto i\,\delta_1(x)\,d\mathbbm{U}_1+\cdots i\,\delta_n(x)\,d\mathbbm{U}_{n}
\end{equation}
for $x$ $\in$ $\Omega^0(\torus^n_\Xi)=\torus^n_\Xi\cong \torus^n_\Xi\otimes \C,$ where
\begin{equation}
    \label{eq.6.2.0}
   d\mathbbm{U}_j:=-{1\over 2\pi}\u^\ast_jd\u_j=(0,...,0,\mathbbm{1},0...,0).
\end{equation}
Notice that the differential can also be viewed as
\begin{equation}
    \label{lonecesito2}
    d(x)= i(\delta_1(x),\delta_2(x),\cdots \delta_n(x)).
\end{equation}
For $x\,d\mathbbm{U}_1\cdots d\mathbbm{U}_k$ $\in$ $\Omega^k(\torus^n_\Xi)$ with $x$ $\in$ $\torus^n_\Xi$, we have $$d(x\,d\mathbbm{U}_1\cdots d\mathbbm{U}_k)=dx\cdot d\mathbbm{U}_1\cdots d\mathbbm{U}_k.$$

This space has a canonical $n$--volume form given by 
\begin{equation}
    \label{eq.6.2.0.1}
    \dvol:=d\mathbbm{U}_1\cdots d\mathbbm{U}_n=\mathbbm{1}\otimes \mathrm{vol}=\mathbbm{1}\otimes \overline{e}_1\wedge\cdots \wedge \overline{e}_n,
\end{equation}
with $\{\overline{e}_i \}^n_{i=1}$ the canonical linear basis of $\C^n$. Notice that  
\begin{equation}
    \label{eq.6.2.1}
    \Omega^n(\torus^n_\Xi)=\torus^n_\Xi\,\dvol,\qquad \,\dvol\,x=x\,\dvol,\qquad x\,\dvol=0\;\Longleftrightarrow \;x=0,
\end{equation}
for all $x$ $\in$ $\torus^n_\Xi$.

On the other hand, taking the inner product  $\langle-|-\rangle_{\wedge}$ on $\bigwedge^\bullet \C^n$ induced by the canonical inner product of $\C^n$, there is an associated star Hodge operator $\star$ and we can define 
\begin{equation}
    \label{eq.6.2.2}
    \star_q: \Omega^k(\torus^n_\Xi)\longrightarrow \Omega^{n-k}(\torus^n_\Xi),\qquad \mu=x\otimes v\longmapsto x\otimes \star v.
\end{equation}
In the same way, taking into account 
\begin{equation}
     \label{eq.6.2.3}
     \int_{\torus^n_\Xi}: \Omega^n(\torus^n_\Xi)\longrightarrow \C,\qquad x\,\dvol\longmapsto \mathrm{Tr}_w(x\,|\D|^{-n}),
\end{equation}
where $\mathrm{Tr}_w$ denotes Dixmier trace (\cite{con}),  we define the inner product of $\Omega^\bullet(\torus^n_\Xi)$ such that differential forms of different degrees are orthogonal, and
\begin{equation}
    \label{eq.6.2.4}
    \langle x_1\otimes v_1\mid x_2\otimes v_2\rangle_\r:=\int_{\torus^n_\Xi} (x_1\otimes v_1)^\ast\,\star_q(x_2\otimes v_2)=\int_{\torus^n_\Xi} x^\ast_1x_2\,(v^\ast_1\star v_2)
\end{equation}
for $x_1\otimes v_1$, $x_2\otimes v_2$ $\in$ $\Omega^k(\torus^n_\Xi)$. According references \cite{lan4,lan2}, the previous inner product matches with the one of equation (\ref{3.3.f18}), and also we have 
\begin{equation}
    \label{eq.6.2.5}
    \int_{\torus^n_\Xi}d\mu=0
\end{equation}
for all $\mu$ $\in$ $\Omega^{n-1}(\torus^n_\Xi)$ \cite{lan4,lan2}.

Now, let us consider the canonical quantum group $\mathcal{U}(1)$ associated to $U(1)$ (see Section 2.2 and the beginning of Section 2.3). In order to develop a concrete example, we will focus in the qpb over $\torus^n_\Xi$ used in Section 6 of reference \cite{dabro2}. This is a quantum principal $\mathcal{U}(1)$--bundle 
\begin{equation}
    \label{eq.6.2}
    \zeta=(P:=\torus^{n+1}_\Xi,\torus^n_\Xi,\Delta_P)
\end{equation}
where  $\Delta_P:P\longrightarrow P\otimes G $ is the $\ast$--algebra morphism defined by
\begin{equation}
    \label{eq.6.4}
    \Delta_P(\u_i)=\u_i\otimes \mathbbm{1},\qquad \Delta_P(\u_{n+1})=\u_{n+1}\otimes z
\end{equation}
for $i=1,...,n$. For a differential calculus on the qpb, let us take 
\begin{equation}
    \label{eq.6.7.1}
    (\Omega^\bullet(P),d,\ast):=(\Omega^\bullet_\D(\torus^{n+1}_\Xi),d_\D,\ast)\cong(\Omega^\bullet(\torus^{n+1}_\Xi),d=id',\ast)
\end{equation}
and
\begin{equation}
    \label{eq.6.7}
\Delta_{\Omega^\bullet(P)}:\Omega^\bullet(P)\longrightarrow \Omega^\bullet(P)\otimes \Gamma^\wedge
\end{equation}
is given by 
\begin{equation*}
    \Delta_{\Omega^\bullet(P)}|_{P}=\Delta_P,  \qquad\qquad \Delta_{\Omega^\bullet(P)}(d\u_i)=d\u_i\otimes \mathbbm{1},
\end{equation*}
$$ \Delta_{\Omega^\bullet(P)}(d\u_{n+1})=d\u_{n+1}\otimes z+\u_{n+1}\otimes dz. $$
The space of horizontal forms is 
\begin{equation}
    \label{eq.6.8}
\Hor^\bullet\,P=P\,\Omega^\bullet(\torus^n_\Xi)=\Omega^\bullet(\torus^n_\Xi)\,P;
\end{equation}
the $\U(1)$--corepresentation on $\Hor^\bullet\,P$ is given by  
\begin{equation}
    \label{eq.6.9}
\Delta_\Hor:\Hor^\bullet\,P\longrightarrow \Hor^\bullet\,P\otimes G,\qquad p\,\mu\longrightarrow p^{(0)}\mu\otimes p^{(1)}
\end{equation}
and the space of base forms is 
\begin{equation}
    \label{eq.6.10}
\Omega^\bullet(\torus^n_\Xi).
\end{equation}

The qpb $\zeta$ with this differential calculus is an example of the class of qpb's presented in Section 3.2 and, for example, equation (\ref{eq.6.8}) illustrates the result presented in equation (\ref{ec.3.31}). In addition, in degree $1$, we have
$$\Omega^1(P)\cong(\Hor^1 P\otimes\mathbbm{1}\C)\oplus (P\otimes \mathfrak{u}(1)^\#_\C)$$ under the identification $$ \vartheta:=\pi(z)=z^\ast dz\;\;\;\longleftrightarrow \;\;\; \u^\ast_{n+1}d\u_{n+1}.$$ Therefore, the regular qpc $\omega^c$ proposed in Proposition (\ref{teodio1}) is given by 
\begin{equation}
    \label{teodio6}
    \omega^c : \mathfrak{u}(1)^\#_\C \longrightarrow \Omega^1(P),\qquad w\,\vartheta\longmapsto w\,\u^\ast_{n+1}d\u_{n+1}
\end{equation}
for all $w$ $\in$ $\C$ (recall that $\mathfrak{u}^\#_\C=\mathrm{span}_\C(\vartheta)$). Let us explicitly show that $\omega^c$ is a regular qpc for our qpb $\zeta$. First of all, according to equation (\ref{properties}) we get  $\vartheta^\ast=\vartheta$  and $$(\u^\ast_{n+1}d\u_{n+1})^\ast=(0,...,0,-2\,\pi\,\mathbbm{1})^\ast=(0,...,0,-2\,\pi\,\mathbbm{1})=\u^\ast_{n+1}d\u_{n+1}.$$ It follows that   $\omega^c(\theta^\ast)=\omega^c(\theta)^\ast$  for all $\theta$ $\in$ $\mathfrak{u}(1)^\#_\C=\mathrm{span}_\C\{ \pi(z)\}$. 

    On the other hand (see equation (\ref{2.f15}))
    \begin{equation}
        \label{adtrivial}
        \ad(\vartheta)=\ad(\pi(z))=(\pi\otimes \id_{\mathfrak{u}(1)^\#_\C})\Ad(z)=\pi(z)\otimes \mathbbm{1}=\vartheta\otimes \mathbbm{1}
    \end{equation}
    and then 
    \begin{eqnarray*}
        \Delta_{\Omega^\bullet(P)}(\omega^c(\vartheta))=\Delta_{\Omega^\bullet(P)}(\mathbbm{u}^\ast_{n+1}d\mathbbm{u}_{n+1})
        &=&
        \Delta_{P}(\mathbbm{u}^\ast_{n+1})\cdot\Delta_{\Omega^\bullet(P)}(d\mathbbm{u}_{n+1})
        \\
        &=&
         \mathbbm{u}^\ast_{n+1}d\mathbbm{u}_{n+1}\otimes z^\ast z
         \\
         &+&
         \mathbbm{u}^\ast_{n+1}\mathbbm{u}_{n+1}\otimes z^\ast dz
        \\
        &=&
         \mathbbm{u}^\ast_{n+1}d\mathbbm{u}_{n+1}\otimes \mathbbm{1}+\mathbbm{1}\otimes \pi(z)
         \\
        &=&
         \mathbbm{u}^\ast_{n+1}d\mathbbm{u}_{n+1}\otimes \mathbbm{1}+\mathbbm{1}\otimes \vartheta
         \\
        &=&
        (\omega^c\otimes \id_{\mathfrak{u}(1)^\#_\C})\ad(\vartheta)+\mathbbm{1}\otimes \vartheta.
    \end{eqnarray*}
    Hence, $\omega^c$ is qpc. Finally, let $g$ $\in$ $G$ and $\varphi$ $\in$ $\Hor^k P$. Then, by equation (\ref{2.2.f4}), we have 
    \begin{eqnarray*}
        \varphi^{(0)}\omega(\vartheta\diamondsuit \varphi^{(1)})=\varphi^{(0)}\omega(\pi(z)\diamondsuit \varphi^{(1)})&=&\varphi^{(0)}\epsilon(\varphi^{(1)})\omega(\pi(z))
        \\
        &=&
        \varphi \,\omega(\pi(z))
        \\
        &=&
        \varphi\, \omega(\vartheta)=(-1)^k \omega(\vartheta)\,\varphi,
    \end{eqnarray*}
   where the last equality is because $\omega(\vartheta)=(0,...,0,-2\pi\mathbbm{1})$ is a 1--degree element of the graded center of $\Omega^\bullet(P)$. We conclude that $\omega$ is a regular qpc.

\begin{Proposition}
\label{flatqpc}
    The qpc $\omega^c$ is flat, i.e., $R^{\omega^c}=0$.
\end{Proposition}
\begin{proof}
    For the Lie group $U(1)$, we have $c=0$; thus $\mathrm{c}^T=0$. This implies that $$ R^{\omega^c}(\vartheta)=d\omega^c(\vartheta)-\langle\omega^c,\omega^c\rangle(\vartheta)=d\omega^c(\vartheta)=d(0,...,0,-2\,\pi\,\mathbbm{1})=0.$$ 
    We conclude that $R^{\omega^c}=0$.
\end{proof}

According to Section 3.1, every qpb $\omega$ is of the form 
\begin{equation}
    \label{eq.6.11}
    \omega=\omega^c+\lambda,
\end{equation}
where $\lambda$ $\in$ $\overrightarrow{\mathfrak{qpc}(\zeta)}$ (see equation  (\ref{2.f24.2})). Since $\vartheta^\ast=\vartheta$ and $\ad(\vartheta)=\vartheta\otimes \mathbbm{1}$, it follows that every $\lambda$ is of the form 
\begin{equation}
    \label{eq.6.12}
    \lambda(\vartheta)=\mu \,\in\, \Omega^1(\torus^n_\Xi) \qquad \mbox{ such that }\qquad \mu^\ast=\mu.
\end{equation}
Moreover, the curvature of every qpc is  
\begin{equation}
    \label{eq.6.13.1}
    R^\omega=d\omega^c+d\lambda=R^{\omega^c}+d\lambda=d\lambda \;\;\Longrightarrow \;\;R^\omega(\vartheta)=d\mu.
\end{equation}

\begin{Proposition}
    \label{zerocurvature1}
    Let $\omega$ be a qpc different from $\omega^c$. Then $\omega$ has zero curvature if and only $$\mu=\sum^n_{j=1}t_j\,d\mathbbm{U}_j+db\quad \mbox{ where }\quad t_j\,\in\,\R\quad \mbox{ and }\quad b\,\in\,\torus^n_\Xi \quad \mbox{ with }\quad b=-b^\ast.$$
\end{Proposition}

\begin{proof}
    Let $\omega\not=\omega^c$ be a qpc with zero curvature. Then $0\not=\mu$ $\in$ $\Ker(d|_{\Omega^1(\torus^n_\Xi)})$ with $\mu=\mu^\ast$. According to \cite{lan4}, the cohomology of the quantum $n$--torus with respect to  
 $(\Omega^\bullet_\D(\torus^n_\Xi), d_\D, \ast) \cong (\Omega^\bullet(\torus^n_\Xi), d = \id', \ast)
$ coincides with the cohomology of the {\it classical} $n$--torus in differential geometry.  In this way, $$\Ker(d|_{\Omega^1(\torus^n_\Xi)})=\left\{\sum^n_{j=1}w_j \,d\mathbbm{U}_j\mid w_j \,\in\,\C\right\}+\{db\mid b\, \in\, \torus^n_\Xi \}$$ and hence, the Proposition follows taking into account that $\mu$ has to be Hermitian and $d$ is anti--Hermitian.
\end{proof}

It is well--known that a complete set of non--equivalent irreducible $\mathcal{U}(1)$--corepresentations is given by  $\T'=\{ \delta^m\mid m\,\in\,\Z \},$ where 
\begin{equation}
    \label{corepu(1)}
    \delta^m: \C\longrightarrow \C\otimes G,\qquad w\longmapsto w\otimes z^m. 
\end{equation}
 
\begin{Proposition}
    \label{generetostorus}
    Let $\delta^m$ $\in$ $\T$ for $m\not=0$. Then the set of maps $\{ T^\l_k\}^{d_m}_{k=1}$ $\subseteq$ $\Mor(\delta^m,\Delta_P)$ (for some $d_m$ $\in$ $\N$) of Remark \ref{rema} are given by $$T^\l_1(1)=\u^m_{n+1} \qquad (d_m=1).$$ 
\end{Proposition}

\begin{proof}
     First of all, notice that $\{1\}$ is an orthonormal basis of the inner product that makes $\delta^n$ unitary. Thus, the linear map $$T^\l_1: \C \longrightarrow P \quad \mbox{ such that }\quad T^\l_1(1)=\u^m_{n+1}=:x^{\C}_{11}$$ satisfies 
    all the requirements.
\end{proof}

Recall that case $m=0$ has already been treated in Proposition \ref{grassmann}, for an arbitrary qpb with an arbitrary differential calculus. In fact, we have that the map $$T^\l_1=T^\triv :\C\longrightarrow P$$ with $T^\triv(1)=\mathbbm{1}$ is the corresponding map for $\Mor(\delta^0,\Delta_P)$ given in Remark \ref{rema}. 

As we have checked in equation (\ref{3.f5}), the map $\{T^\l_1 \}$ form a set of $\Omega^\bullet(\torus^n_\Xi)$--left generators of $\Mor(\delta^m,\Delta_\Hor)$. In the particular case of this section, $\{T^\l_1 \}$ is actually a $\Omega^\bullet(\torus^n_\Xi)$--bimodule basis and a direct calculation shows that $T^\l_1=T^\r_1$.

Since every finite--dimensional $\U(1)$--corepresentation is the finite direct sum of elements of $\T$ and $\Delta_\Hor|_P=\Delta_P$, it follows that the associated left/right qvb of every finite--dimensional $\U(1)$--corepresentation is a free left/right $\torus^n_\Xi$--module, i.e., a trivial left/right qvb, in light of the Serre--Swan theorem.

\subsubsection{The Gauge Dirac Operator}
To exemplify this part of the theory, we will focus in the gauge Dirac theory associated to $\omega^c$ for $\delta^1$ $\in$ $\T$. For the differential calculus on $\zeta$ based on the universal differential $\ast$--algebra of $P$, the linear map $\omega^c$ is also a qpc (and the proof is completely similar); however, it is not regular and $R^{\omega^c}\not=0$. For this differential calculus, the dual covariant derivative of $T^\r_1(1)$ is  (see equation (\ref{2.f30.5.1})) 
\begin{eqnarray*}
\widehat{D}^{\omega^c}_U(T^\r_1(1))=\widehat{D}^{\omega^c}_U\u_{n+1}&=&d_U\u_{n+1}+\omega^c(\pi(S^{-1}(z)))\u_{n+1}
\\
&=&
d_U\u_{n+1}-\u^\ast_{n+1}(d_U\u_{n+1})\u_{n+1}
\end{eqnarray*}
 and the gauge qlc is (see equation (\ref{3.f8.1})) $${^U}\widehat{\nabla}^{\omega^c}_{\C} (T^\r_1)=T^\r_1\otimes_{\torus^n_\Xi} (\mu^{_{\widehat{D}^{\omega^c}_U(T^\r_1)^\ast}})^\ast,$$ where $$ \mu^{_{\widehat{D}^{\omega^c}_U(T^\r_1)^\ast}}=\widehat{D}^{\omega^c}_U(T^\r_1(1))^\ast\,x^{\C\,\ast}_{11}=\widehat{D}^{\omega^c}_U(T^\r_1(1))^\ast\,\u_{n+1}.$$ Since every element of $E^\C_\r=\Mor(\delta^1,\Delta_P)$ is of the form (see equation (\ref{3.f5.2})) $$T=T^\r_1\,(b^{_{T^\ast}})^\ast \qquad \mbox{ for }\qquad b^{_{T^\ast}}=T(1)^\ast\,x^{\C\,\ast}_{11}=T(1)^\ast\,\u_{n+1}\in\, \torus^n_\Xi,$$ we obtain by the right Leibniz rule
\begin{eqnarray*}
    {^U}\widehat{\nabla}^{\omega^c}_{\C} (T)={^U}\widehat{\nabla}^{\omega^c}_{\C}(T^\r_1\,(b^{_{T^\ast}})^\ast)&=&{^U}\widehat{\nabla}^{\omega^c}_{\C}(T^\r_1)\,(b^{_{T^\ast}})^\ast+T^\r_1\otimes_{\torus^n_\Xi} d_U((b^{_{T^\ast}})^\ast)
    \\
    &=&
    T^\r_1\otimes_{\torus^n_\Xi} (\mu^{_{\widehat{D}^{\omega^c}_U(T^\r_1)^\ast}})^\ast\,(b^{_{T^\ast}})^\ast-T^\r_1\otimes_{\torus^n_\Xi} d_Ub^{_{T^\ast}} 
    \\
    &=&
    T^\r_1\otimes_{\torus^n_\Xi} (b^{_{T^\ast}}\mu^{_{\widehat{D}^{\omega^c}_U(T^\r_1)^\ast}})^\ast-T^\r_1\otimes_{\torus^n_\Xi} d_Ub^{_{T^\ast}} 
    \\
    &=&
    T^\r_1\otimes_{\torus^n_\Xi}[(b^{_{T^\ast}}\mu^{_{\widehat{D}^{\omega^c}_U(T^\r_1)^\ast}})^\ast-d_Ub^{_{T^\ast}}].
\end{eqnarray*}

On the other hand, for the differential calculus on $\zeta$ based on  $(\Omega^\bullet(P),d=id',\ast)$, we have that 
\begin{eqnarray*}
    \widehat{D}^{\omega^c}(T^\r_1(1))=\widehat{D}^{\omega^c}\u_{n+1}&=&d\u_{n+1}+\omega^c(\pi(S^{-1}(z)))\u_{n+1}
    \\
    &=&
    d\u_{n+1}-\u^\ast_{n+1}(d\u_{n+1})\u_{n+1}=0.
\end{eqnarray*}
Thus $$\widehat{\nabla}^{\omega^c}_{\C} (T^\r_1)=0.$$  Since every element of $E^\C_\r$ is of the form  $$T=T^\r_1\,(b^{_{T^\ast}})^\ast \qquad \mbox{ for }\qquad b^{_{T^\ast}}=T(1)^\ast\,x^{\C\,\ast}_{11}=T(1)^\ast\,\u_{n+1}\in\, \torus^n_\Xi,$$ we obtain by the right Leibniz rule
$$ \widehat{\nabla}^{\omega^c}_{\C} (T)= \widehat{\nabla}^{\omega^c}_{\C}(T^\r_1\,(b^{_{T^\ast}})^\ast)=T^\r_1\otimes_{\torus^n_\Xi} d((b^{_{T^\ast}})^\ast)=-T^\r_1\otimes_{\torus^n_\Xi} d(b^{_{T^\ast}}).$$ Under the projection $\pi_\r$ (see equation (\ref{3.3.f2})) we have that 
\begin{eqnarray*}
    \pi_\r(\widehat{D}_U^{\omega^c}(T^\r_1(1)))=
    \widehat{D}^{\omega^c}(T^\r_1(1))=0
\end{eqnarray*}
and it follows that $$(\id_{E^\C_\r}\otimes_{\torus^n_\Xi} \pi_\r)\; {^U}\widehat{\nabla}^{\omega^c}_{\C}= \widehat{\nabla}^{\omega^c}_{\C}.$$ Hence, ${^U}\widehat{\nabla}^{\omega^c}_{\C}$ is a qlc with respect to $(\Omega^\bullet_U(\torus^n),d_\U,\ast)$ compatible with the metric (by Theorem \ref{fgs}) that projects to $\widehat{\nabla}^{\omega^c}_{\C}$ and therefore, we can take
$$\D_{\widehat{\nabla}^{\omega^c}_\C}(T\otimes_{\torus^n_\Xi} \psi)=T\otimes_{\torus^n_\Xi} \D\psi+\left((\id_{E^V_\r}\otimes_{\torus^n_\Xi} \pi_\r)\circ\, {^{U}\widehat{\nabla}^{\omega^c}_\C}T \right)\psi.$$ This is the gauge Dirac operator induced by $\omega^c$.

Let $\Psi=\displaystyle\sum_{_k}T_k\otimes_{\torus^n_\Xi} \psi_k=\sum_{_k}T^\r_1\,(b^{_{T^\ast}}_k)^\ast\otimes_{\torus^n_\Xi} \psi_k$ $\in$ $\qgs_\r$. Then 
\begin{equation}
\label{finaldirac}
  \D_{\widehat{\nabla}^{\omega^c}_\C}\Psi=0 \;\;\Longleftrightarrow\;\; \sum_k (b^{{_{T^\ast}}}_k\,\D - d(b^{_{T^\ast}}_k))\,\psi_k=0.
\end{equation}
This characterize all non--commutative geometrical Dirac spinor fields. For example, if $\psi$ is an eigenvector of $\D$ with eigenvalue $0$ (\cite{con}), then $\Psi=T^\r_1\otimes \psi$ is a non--commutative geometrical Dirac spinor field.

As we commented in Remark \ref{remayanmills}, the importance of this part of the theory lies on the fact that we can always take the gauge Dirac operator induced by any qpc for any qpb (in particular, for any quantum group and any bicovariant $\ast$--FODC on the quantum group) with the base space $B$ having structure of a $n$--dimensional spectral triple.

\subsubsection{Yang--Mills Theory}

As we mentioned at the end of Section 5, the importance of this part of the theory clearly lies in the non--commutative geometrical Yang--Mills equation. In particular, in this subsection, we will compare this equation with the well--known non--commutative analytic Yang--Mills equation.

\begin{Theorem}
\label{qpc-qlc1}
    Under the assignation $$\mathfrak{qpc}(\zeta)\longrightarrow \mathfrak{qlc}_C(E^{\mathfrak{u}(1)^\#_\C}_\r), \qquad \omega\longmapsto \widehat{\nabla}^\omega_{{\mathfrak{u}(1)^\#_\C}}$$ every qpc is a solution of the non--commutative analytical Yang--Mills equation on the associated right qvb $E^{\mathfrak{u}(1)^\#_\C}_\r$ of the $\U(1)$--corepresentation $\ad$. In other words, the gauge qlc of every qpc is a Yang--Mills qlc.
\end{Theorem}

\begin{proof}
    Let $\omega$ be a qpc. By equation (\ref{adtrivial}), we known that $\ad(\theta)=\theta\otimes \mathbbm{1}$  for every $\theta$ $\in$ $\mathfrak{u}(1)^\#_\C=\mathrm{span}_\C\{\vartheta \}$. Therefore, $\ad\cong\delta^0=\delta^\C_\triv$ and hence, the associated right qvb $E^{\mathfrak{u}(1)^\#_\C}_\r$ for $\ad$ is isomorphic to the associated qvb $E^{\C}_\r$ for $\delta^\C_\triv$ (and in this case, $E^{\C}_\r\cong \torus^n_\Xi$). In light of Proposition \ref{grassmann}, we have that the gauge qlc $$\widehat{\nabla}^\omega_{{\mathfrak{u}(1)^\#_\C}}=:\widehat{\nabla}^\omega$$ induced by $\omega$ is the Grassmann connection: $$\widehat{\nabla}^\omega(T)=T^\triv\otimes_B db^{_T},$$ with $T=T^\triv\,b^{_T}$. Since  $\widehat{\nabla}^\omega(T^\triv)= T^\triv\otimes_{\torus^n_\Xi}db^{{T^\triv}}=0,$  a straightforward calculation shows that 
\begin{eqnarray*}
    R^{\widehat{\nabla}^\omega}(T)=(d^{\widehat{\nabla}^\omega}\circ \widehat{\nabla}^\omega)(T)=
   0. 
\end{eqnarray*}
Hence, $\widehat{\nabla}^\omega:=\widehat{\nabla}^\omega_{{\mathfrak{u}(1)^\#_\C}}$ is trivially a solution of the non--commutative analytical Yang--Mills equation (equation (\ref{ymanalitic})).
\end{proof}

Of course, these are not all the solutions of the non--commutative analytical Yang--Mills equation; however, these are all the solutions that come from the quantum principal bundle structure. 

Let $$\widehat{\nabla}: E^{\mathfrak{u}(1)^\#_\C}_\r\longrightarrow E^{\mathfrak{u}(1)^\#_\C}_\r\otimes_{\torus^n_\Xi} \Omega^1(\torus^n_\Xi) $$ be a qlc of $E^{\mathfrak{u}(1)^\#_\C}_\r$ compatible with the metric different from $\widehat{\nabla}^\omega$. Then, according to \cite{lan}, there exists an Hermitian right $\torus^n_\Xi$--linear map $$\widehat{\Lambda}: E^{\mathfrak{u}(1)^\#_\C}_\r \longrightarrow E^{\mathfrak{u}(1)^\#_\C}_\r\otimes_{\torus^n_\Xi}\Omega^1(\torus^n_\Xi)\qquad \mbox{such that }\qquad \widehat{\nabla}=\widehat{\nabla}^\omega+\widehat{\Lambda}.$$  Since $\{ T^\triv\}$ is a $\torus^n_\Xi$--basis of $E^{\mathfrak{u}(1)^\#_\C}_\r$, the map $\widehat{\Lambda}$ is completely characterized by $$\widehat{\Lambda}(T^\triv)=T^\triv \otimes_{\torus^n_\Xi} \mu$$  for some $0\not=\mu$ $\in$ $\Omega^1(\torus^n_\Xi)$  such that $\mu=\mu^\ast.$ 
    
    In this way, the curvature of  $\widehat{\nabla}$, which is a right $\torus^n_\Xi$--linear map, is completely characterized by 
    \begin{eqnarray*}
        R^{\widehat{\nabla}}(T^\triv)=d^{\widehat{\nabla}}(\widehat{\nabla}(T^\triv)) = d^{\widehat{\nabla}}(\widehat{\nabla}^\omega(T^\triv)+\widehat{\Lambda}(T^\triv))
        =
        T^{\triv}\otimes_{\torus^n_\Xi} (d\mu+\mu\, \mu).
    \end{eqnarray*}

    Assume that $\widehat{\nabla}$ $\in$ $\mathfrak{qlc}_C(E^{\mathfrak{u}(1)^\#_\C}_\r)$ is Yang--Mills qlc. Then, by equation (\ref{ymanalitic}), we have $$0=[\widehat{\nabla},(\id\otimes_{B} \star_q)\circ R^{\widehat{\nabla}}]$$ and it follows that (see equation (\ref{4.3.f25})) 
    \begin{eqnarray*}
        0&=&[\widehat{\nabla},(\id\otimes_{\torus^n_\Xi} \star_q)R^{\widehat{\nabla}}](T^\triv)
        \\
        &=&
        d^{\widehat{\nabla}}(\id\otimes_{\torus^n_\Xi} \star_q)R^{\widehat{\nabla}}(T^\triv)-(-1)^{n-2}\, ((\id\otimes_{\torus^n_\Xi} \star_q)R^{\widehat{\nabla}})^\mathrm{ext}\;  \widehat{\nabla}(T^\triv)
        \\
        &=&
        d^{\widehat{\nabla}}(\id\otimes_{\torus^n_\Xi} \star_q)R^{\widehat{\nabla}}(T^\triv)-(-1)^{n-2}\, ((\id\otimes_{\torus^n_\Xi} \star_q)R^{\widehat{\nabla}})^\mathrm{ext}\; \widehat{\Lambda}(T^\triv)
        \\
        &=&
        d^{\widehat{\nabla}}(T^{\triv}\otimes_{\torus^n_\Xi} \star_q(d\mu+\mu\, \mu))-(-1)^{n-2}\,((\id\otimes_{B} \star_q)R^{\widehat{\nabla}})^\mathrm{ext}\;(T^\triv\otimes_{\torus^n_\Xi} \mu)
        \\
        &=&
        d^{\widehat{\nabla}}(T^{\triv}\otimes_{\torus^n_\Xi} \star_q(d\mu+\mu\, \mu))-(-1)^{n-2}\,((\id\otimes_{B} \star_q)R^{\widehat{\nabla}}(T^\triv))\,\mu
        \\
        &=&
        d^{\widehat{\nabla}}(T^{\triv}\otimes_{\torus^n_\Xi} \star_q(d\mu+\mu\, \mu))-(-1)^{n-2}\, T^{\triv}\otimes_{\torus^n_\Xi}(\star_q(d\mu+\mu\, \mu))\,\mu
        \\
        &=&
        \widehat{\nabla}(T^{\triv})\star_q(d\mu+\mu\, \mu)+T^{\triv}\otimes_{\torus^n_\Xi} d(\star_q(d\mu+\mu\, \mu))
        \\
        &-&
        (-1)^{n-2} \,T^{\triv}\otimes_{\torus^n_\Xi}(\star_q(d\mu+\mu\, \mu))\,\mu
        \\
        &=&
        \widehat{\Lambda}(T^{\triv})\star_q(d\mu+\mu\, \mu)+T^{\triv}\otimes_{\torus^n_\Xi} d(\star_q(d\mu+\mu\, \mu))
        \\
        &+&
        (-1)^{n-1} \,T^{\triv}\otimes_{\torus^n_\Xi}(\star_q(d\mu+\mu\, \mu))\,\mu
        \\
        &=&
        T^{\triv}\otimes_{\torus^n_\Xi} [d(\star_q(d\mu+\mu\, \mu))+\mu \,\star_q(d\mu+\mu\, \mu)+(-1)^{n-1}\star_q(d\mu+\mu\, \mu)\,\mu]
    \end{eqnarray*}

\noindent In this way, we have already proven the following statement 

\begin{Theorem}
\label{qpc-qlc2}
    Let $\widehat{\nabla}$ be a qlc compatible with the metric different from $\widehat{\nabla}^\omega$. Then,  $\widehat{\nabla}$ is a Yang--Mills qlc if and only if
    \begin{equation}
        \label{ahuevo1}
        d(\star_q(d\mu+\mu\, \mu))+\mu \,\star_q(d\mu+\mu\, \mu)+(-1)^{n-1}\star_q(d\mu+\mu\, \mu)\,\mu=0.
    \end{equation}
\end{Theorem}

    According to \cite{con}, if $\widehat{\nabla}$ is a Yang--Mills qlc's, then it must have zero curvature. Hence, equation (\ref{ahuevo1}) turns into 
    \begin{equation}
        \label{yangmillsqlc2}
        d\mu+\mu\, \mu=0.
    \end{equation}
    
     Moreover, in light of \cite{con}, under the standard action on $\mathfrak{qlc}_C(E^{\mathfrak{u}(1)^\#_\C}_\r)$ of the group $U(E^{\mathfrak{u}(1)^\#_\C}_\r)$ of all unitary operators of $E^{\mathfrak{u}(1)^\#_\C}_\r$, there is only one orbit of Yang--Mills qlc's. In other words, for a Yang--Mills qlc $\widehat{\nabla}$ different from $\widehat{\nabla}^\omega$, there exists a unitary (with respect to the canonical Hermitian structure $\langle-,-\rangle_\r$) right $\torus^n_\Xi$--linear map $$\mathbf{A}: E^{\mathfrak{u}(1)^\#_\C}_\r \longrightarrow E^{\mathfrak{u}(1)^\#_\C}_\r \qquad \mbox{ such that }\qquad (\mathbf{A}\otimes_{\torus^n_\Xi}\id_{\Omega^1(\torus^n_\Xi)})\circ \widehat{\nabla}^\omega\circ \mathbf{A}^{-1}=\widehat{\nabla}.$$ The map $\mathbf{A}$ is completely characterized by $$\mathbf{A}(T^\triv)=T^\triv\, x\quad \mbox{ with }\quad \mathbbm{1}\not=x\,\in\, \torus^n_\Xi \quad \mbox{ such that }\quad  x^{-1}=x^\ast.$$ In this way, we have $$\widehat{\nabla}(T^\triv)=(\mathbf{A}\otimes_{\torus^n_\Xi}\id_{\Omega^1(\torus^n_\Xi)})\left( \widehat{\nabla}^\omega\left( \mathbf{A}^{-1}(T^\triv)\right)\right)=T^\triv \otimes_{\torus^n_\Xi} x\,dx^\ast$$ and this implies that $\mu=x\,dx^\ast.$ Since $\mathbbm{1}\not=x$ is unitary, it follows that $$\mu=2\,\pi\,(m_1\,d\mathbbm{U}_1+\cdots m_n\,d\mathbbm{U}_n)$$ for some $m_j$ $\in$ $\Z$ not all zero. In summary, for every $n\geq 2$, every Yang--Mills qlc is characterized by 
\begin{equation}
        \label{yangmillsqlc3}
       \mu=2\,\pi\,(m_1\,d\mathbbm{U}_1+\cdots m_n\,d\mathbbm{U}_n)
    \end{equation}
for $m_j$ $\in$ $\Z$, $j=1,....,n$. 

Now, we proceed to analyze solutions of the non--commutative geometrical Yang--Mills equation. First, recall that for the Lie group $U(1)$, the non--commutative geometrical Yang--Mills equation takes the form (see equation (\ref{4.3.f29.2}))
     $$d^{\widehat{\nabla}^\omega\star}(\widehat{\Upsilon}_{\mathfrak{u}(1)^\#_\C}(R^\omega))=0.$$ 
By Proposition \ref{flatqpc}, it follows that $\omega^c$ is a Yang--Mills qpc;  but also we have

\begin{Theorem}
\label{solqpc1}
    Consider the characterization of every qpc $\omega$ and its curvature presented in equations (\ref{eq.6.11})--(\ref{eq.6.13.1}). Let $\omega$ be a qpc different from $\omega^c$. Then, $\omega$ is Yang--Mills qpc if and only if
    \begin{equation}
        \label{eq.6.15}
        d^{\star_q}d\mu=0.
    \end{equation}
\end{Theorem}

\begin{proof}
    Let $\omega=\omega^c+\lambda$ $\in$ $\mathfrak{qpc}(\zeta)$ with $\lambda(\vartheta)=\mu\not=0$. Then we know that $$R^\omega=d\lambda \quad \mbox{ with }\quad R^\omega(\vartheta)=d\lambda(\vartheta)=d\mu \quad \Longrightarrow \quad R^\omega=T^\triv\,d\mu.$$ In this way, since  $\widehat{\nabla}^\omega(T^\triv)= T^\triv\otimes_{\torus^n_\Xi}db^{{T^\triv}}=0,$ we obtain (see equation (\ref{formaldiff}))
     \begin{eqnarray*}
         0=d^{\widehat{\nabla}^\omega\star}(\widehat{\Upsilon}_{\mathfrak{u}(1)^\#_\C}(R^\omega))  &=&d^{\widehat{\nabla}^\omega\star}(\widehat{\Upsilon}_{\mathfrak{u}(1)^\#_\C}(T^\triv\,d\mu))
         \\
         &=&d^{\widehat{\nabla}^\omega\star}(T^\triv\otimes_{\torus^n_\Xi} d\mu)
         \\
         &=&-(\id\otimes_{\torus^n_\Xi}\star^{-1}_q) d^{\widehat{\nabla}^\omega}(\id\otimes_{\torus^n_\Xi}\star_q)(T^\triv\otimes_{\torus^n_\Xi} d\mu)
         \\
         &=&
         -(\id\otimes_{\torus^n_\Xi}\star^{-1}_q) d^{\widehat{\nabla}^\omega} (T^\triv\otimes_{\torus^n_\Xi} \star_q(d\mu))
         \\
         &=&
         -(\id\otimes_{\torus^n_\Xi}\star^{-1}_q)[\widehat{\nabla}^\omega(T^\triv)\star_q(d\mu) 
         +
         (T^\triv\otimes_{\torus^n_\Xi} d(\star_q(d\mu)))]
         \\
         &=&
         -(\id\otimes_{\torus^n_\Xi}\star^{-1}_q)  (T^\triv\otimes_{\torus^n_\Xi} d(\star_q(d\mu)))
         \\
         &=&
         - T^\triv\otimes_{\torus^n_\Xi} \star^{-1}_q(d(\star_q(d\mu))).
     \end{eqnarray*}
    The last equality implies  $0=-\star^{-1}_q(d(\star_q(d\mu)))$; so $=d^{\star}d\mu=0.$  
\end{proof}

It is worth mentioning that the equation that describes Yang--Mills qlc's (equation (\ref{ahuevo1})) and the equation that describes Yang--Mills qpc's (equation (\ref{eq.6.15})) {\bf are different}.

Let $\omega=\omega^c+\lambda$ be a Yang--Mills qpc. Then $\lambda$ is of the form $$\lambda(\vartheta)=\mu=\sum^n_{j=1}x_j\,d\mathbbm{U}_j$$ for some hermitian elements $x_j$ $\in$ $\torus^n_\Xi$. Thus $$0=d^{\star}d\mu=\sum^n_{j=1}d^{\star}d(x_j)\,d\mathbbm{U}_j=\sum^n_{j=1}\Delta_\l(x_j)\,d\mathbbm{U}_j \;\;\Longrightarrow \;\; \Delta_\l(x_j)=0$$ for all $j=1,...,n$, where $\Delta_\l$ is the (flat) Laplacian \cite{lap}. According to \cite{lap}, the only eigenvectors of $\Delta_\l$ with eigenvalue zero is $w\,\mathbbm{1}$, for $w$ $\in$ $\C$. Since $x_j$ is an hermitian element, we conclude that 
\begin{equation}
    \label{yangmillsqpc1}
    \mu=\sum^n_{j=1}t_j\,d\mathbbm{U}_j\qquad \mbox{ where }\qquad  t_j \,\in\,\R.
\end{equation}
This characterizes all Yang--Mills qpc's for every $n\geq 2$.  

\begin{Proposition}
    \label{gaugetorus2}
    Let $\omega=\omega^c+\lambda$ be a Yang--Mills qpc with $\lambda\not=0$. Then, there exists $\F$ $\in$ $\qGG_\YM$ such that (see Definition (\ref{6.1.3}) and equation (\ref{actionqpc})) $$\F^{\circledast}\omega^\c=\omega.$$ In other words,  up to an element of $\qGG_\YM$, $\omega^c$ is the only Yang--Mills qpc. 
\end{Proposition}

\begin{proof}
    We know that $$\omega(\vartheta)=\omega^c(\vartheta)+\lambda(\vartheta) \quad \mbox{ with }\quad \lambda(\vartheta)=\mu\not=0,  \quad \mu=\mu^\ast,\quad  d\mu=0.$$ Now, consider the graded left $\Omega^\bullet(\torus^n_\Xi)$--module isomorphism $$\F:=m_\Omega\circ (\id_{\Omega^\bullet(P)}\otimes \f)\circ \Delta_{\Omega^\bullet(P)} :\Omega^\bullet(P)\longrightarrow \Omega^\bullet(P),$$ where  $m_\Omega: \Omega^\bullet(P)\otimes \Omega^\bullet(P)\longrightarrow \Omega^\bullet(P)$ is the product map and $$\f: \Gamma^\wedge=G\oplus \Gamma\longrightarrow \Omega^\bullet(\torus^3_\Xi)$$ is the graded linear map defined by $$\f|_G=\epsilon\, \mathbbm{1}\qquad \mbox{ and }\qquad  \f|_{\Gamma}(g\,\vartheta)=\epsilon(g)\mu$$ recalling that  $\Gamma=G\,\mathfrak{u}(1)^\#_\C=G\,\mathrm{span}_\C\{\vartheta \}=\mathrm{span}_G\{\vartheta \}.$ The inverse map of  $\F$ is given by 
    $$\F^{-1}:=m_\Omega\circ (\id_{\Omega^\bullet(P)}\otimes \f^{-1})\circ \Delta_{\Omega^\bullet(P)} :\Omega^\bullet(P)\longrightarrow \Omega^\bullet(P),$$ where $$\f^{-1}: \Gamma^\wedge=G\oplus \Gamma\longrightarrow \Omega^\bullet(\torus^3_\Xi)$$ is the graded linear map defined by  $$\f^{-1}|_G=\epsilon\, \mathbbm{1}\qquad \mbox{ and } \qquad \f^{-1}|_{\Gamma}(g\,\vartheta)=-\epsilon(g)\,\mu\otimes \mathbbm{1}.$$ The notation $\f^{-1}$ is because this map is convolution inverse of $\f$ (see Section 4 of reference \cite{sald2}). Since $\vartheta=\vartheta^\ast$ and $\mu=\mu^\ast$, it follows that  $\f\circ \ast=\ast\circ \f $ and a long but straightforward direct computation shows that $\F \in\qGG$.
    
    According to Proposition 4.8 of reference \cite{sald2}, we obtain 
    \begin{eqnarray*}
        \F^{\circledast}\omega^\c(\vartheta)=m_{\Omega}(\omega^\c\otimes \f )\ad(\vartheta) +\f(\vartheta)
        &=&
        m_{\Omega}(\omega^c\otimes \f )(\vartheta\otimes \mathbbm{1}) +\f(\vartheta)
         \\
        &=&
\omega^c(\vartheta)\f(\mathbbm{1})+f(\vartheta)
        \\
        &=&
        \omega^c(\vartheta) +\mu 
        \\
        &=&
        \omega(\vartheta).
    \end{eqnarray*}
    Hence  $ \F^{\circledast}\omega^\c=\omega.$ 
    
    To finalize the proof, we are going to show that $\F$ $\in$ $\qGG_\YM$. Let $\omega'=\omega^c+\lambda'$ be any other qpc. Then, we get $$\omega'(\vartheta)=\omega^c(\vartheta)+\mu' \qquad \mbox{ with }\qquad \lambda'(\vartheta)=\mu',\quad \mu'=\mu'^\ast.$$  As above, Proposition 4.8 of reference \cite{sald2} implies that  $$\F^{\circledast}\omega'(\vartheta)=\omega'(\vartheta)+\mu=\omega^c(\vartheta)+\mu'+\mu$$ and therefore $$R^{\F^{\circledast}\omega'}(\vartheta)=d\mu'+d\mu=d\mu'=R^{\omega'}(\vartheta)\;\;\Longrightarrow \;\; R^{\F^{\circledast}\omega'}= R^{\omega'}$$ Thus    $||R^{\F^{\circledast}\omega'}||^2=||R^{\omega'}||^2$  and we conclude that $\F$ $\in$ $\qGG_\YM$.
\end{proof}

Let us summarize the results concerning the two non--commutative Yang--Mills functionals and compare them with each other in order to provide a {\bf first} answer to the questions raised at the end of the Section 5 for the qpb $\zeta$ of this subsection. 
\begin{enumerate}
    \item The gauge qlc  $\widehat{\nabla}^\omega:=\widehat{\nabla}^\omega_{{\mathfrak{u}(1)^\#_\C}}$  of every qpc $\omega$ is a Yang--Mills qlc of $E^{\mathfrak{u}(1)^\#_\C}_\r$.
    \item Equation (\ref{ahuevo1}) characterize all Yang--Mills qlc's of $E^{\mathfrak{u}(1)^\#_\C}_\r$ for every $n\geq 2$. This is equivalently to ask for qlc compatible with the metric with zero curvature and in equation (\ref{yangmillsqlc3}) we characterize all these qlc's.
    \item Up to an element of $U(E^{\mathfrak{u}(1)^\#_\C}_\r)$,  $\widehat{\nabla}^\omega$  is the only Yang--Mills qlc.
    \item Equation (\ref{eq.6.15}) characterize all Yang--Mills qpc's for every $n\geq 2$. This is equivalent to equation (\ref{yangmillsqpc1}). 
    \item The non--commutative analytical Yang--Mills equation (equation (\ref{ahuevo1})) and the non--commutative geometrical Yang--Mills (equation (\ref{eq.6.15})) are different.
    \item The solutions presented in equation (\ref{yangmillsqlc3}) are also solutions of equation (\ref{yangmillsqpc1}).
    \item Up to an element of $\qGG_\YM$, $\omega^c$ is the only Yang--Mills qpc.
\end{enumerate}

\begin{Remark}
    \label{gaugeboson1}
    In light of the interpretation of Remark \ref{gaugeboson}, $\omega^c$ (and quantum gauge transformation of this qpc) is a non--commutative photon field moving through the quantum $n$--torus, where $U(1)$ is equipped with its canonical differential structure. As we should expect, module (quantum) gauge transformation, there is only one (non--commutative) photon field.
\end{Remark}

Now, we will focus on showing that our theory accepts some modifications with interesting results.

\subsection{A Modification of the Theory: The $1$--dimensional Non--Standard Differential Calculus with Non--Trivial Higher--Order components}
 
 As we have mentioned throughout the paper, our goal is to develop a {\it physical} Yang--Mills theory based on qpb's. To this end, we have imposed the condition that the quantum group of our qpb's must be (the canonical quantum group associated with) $SU(N)$ for some $N \geq 1$. Moreover, in this case, the quantum group possesses a canonical bicovariant $\ast$--FODC: the space of {\it classical} $\C$--valued differential $1$--forms on $SU(N)$ (Section 2.2).
 
However, in purely mathematical terms, there is no need to restrict ourselves to {\it classical} differential calculus on $SU(N)$; any other finite-dimensional bicovariant $\ast$--FODC of $SU(N)$ may also be used. This section is to show that such changes are possible in the theory. However, in this situation, in general, there
is no a canonical embedded differential and therefore, it is quite challenging to prove Theorem 5.18 in this context; so we will focus solely in a concrete example.

Consider the qpb
\begin{equation}
    \label{new.1}
    \zeta=(P:=\torus^n_\Xi\otimes G,\torus^n_\Xi,\Delta_P:=\id_{\torus^n_\Xi}\otimes \Delta)
\end{equation}
and the differential calculus on $\zeta$ will be given by the tensor products of graded differential $\ast$--algebras
\begin{equation}
    \label{new.2}
    (\Omega^\bullet(P):=\Omega^\bullet(\torus^n_\Xi)\otimes \Gamma^\wedge,d_\otimes,\ast),
\end{equation}
and $$\Delta_{\Omega^\bullet(P)}:=\id_{\Omega^\bullet(\torus^n_\Xi)}\otimes \Delta:\Omega^\bullet(P)\longrightarrow \Omega^\bullet(P)\otimes \Gamma^\wedge,$$  where the graded differential $\ast$--algebra $$ (\Gamma^\wedge,d,\ast)$$ is the one of Section 2.3. We will refer to the $\ast$--FODC of Section 2.3 as the {\it $1$--dimensional non--standard} $\ast$--FODC of $U(1)$. 

It is worth mentioning that we have changed from the qpb of equation (\ref{eq.6.2}) to the one of equation (\ref{new.1}) because the differential calculus of the qpb of equation (\ref{eq.6.2}) is not compatible with the $1$--dimensional non--standard $\ast$--FODC of $U(1)$. 

Notice that the space of horizontal forms is 
\begin{equation}
    \label{hor1}
    \Hor^\bullet\,P=\Omega^\bullet(\torus^n_\Xi)\otimes G 
\end{equation}
and the space of base forms is given by 
\begin{equation}
    \label{base1}
    \Omega^\bullet(\torus^n_\Xi)\otimes \mathbbm{1}\cong \Omega^\bullet(\torus^n_\Xi).
\end{equation}
To simplify the notation, in some calculation we will consider that the space of base forms is exactly $\Omega^\bullet(\torus^n_\Xi)$. Furthermore, notice that the $\U(1)$--corepresentation $\ad$ is still isomorphic to the trivial $\U(1)$--corepresentation on $\C$ (see equation (\ref{2.3.f16})). Hence  $c^T=0$  as in the {\it classical} $\ast$--FODC of $U(1)$. 

Recall that in the case of the {\it classical} $\ast$--FODC of $U(1)$, the embedded differential used is  $$-\displaystyle \displaystyle {1\over 2}c^T=0.$$ In this way, one can think that  $\displaystyle -{1\over 2}c^T=0$  is also an embedded differential for the $1$--dimensional non--standard $\ast$--FODC. The following proposition shows that this is not the case; however, since the proof consists of a straightforward calculation verifying that Definition \ref{embbededdifferential} is satisfied, we will omit it.  

\begin{Proposition}
    \label{embeddeddif}
    For the $1$--dimensional non--standard $\ast$--FODC of $U(1)$, the only embedded differential $$\Theta: \mathfrak{qu}(1)^\#\longrightarrow \mathfrak{qu}(1)^\#\otimes \mathfrak{qu}(1)^\#$$ is given by $$\Theta(\vartheta)=-\vartheta\otimes \vartheta,$$ where $\vartheta=\pi(z)$. 
\end{Proposition}

As in the previous subsection, the qpb $\zeta$ with this differential calculus is an example of the class of qpb's presented in Section 3.2 and, for example, equation (\ref{hor1}) illustrates the result presented in equation (\ref{ec.3.31}). In addition, in degree $1$, we have
$$\Omega^1(P)\cong(\Hor^1 P\otimes\mathbbm{1}\C)\oplus (P\otimes \mathfrak{qu}(1)^\#)$$ under the identification $$\vartheta:=\pi(z)=z^\ast dz\;\;\;\longleftrightarrow \;\;\; \mathbbm{1}\otimes \vartheta.$$ Therefore, the regular qpc $\omega^c$ proposed in Proposition (\ref{teodio1}) is given by 
\begin{equation}
    \label{teodio7}
    \omega^c : \mathfrak{qu}(1)^\#  \longrightarrow \Omega^1(P),\qquad w\,\vartheta\longmapsto w\,\mathbbm{1}\otimes \vartheta
\end{equation}
for all $w$ $\in$ $\C$ (recall that $\mathfrak{qu}^\#=\mathrm{span}_\C(\vartheta)$). As in the previous subsection, one can explicitly show that $\omega^c$ is a regular qpc for our qpb $\zeta$. In addition, as in the previous subsection, we have the following result.
\begin{Proposition}
\label{canoqpc1}
    The qpc $\omega^c$ is flat, i.e., $R^{\omega^c}=0$.
\end{Proposition}

\begin{proof}
    By equation (\ref{2.3.f18}), we get
    \begin{eqnarray*}
        R^{\omega^c}(\vartheta)=d\omega^c(\vartheta)-\langle \omega^c,\omega^c\rangle(\vartheta)&=&d\omega^c(\vartheta)+  \omega^c(\vartheta)\,\omega^c(\vartheta)
        \\
        &=&
        \mathbbm{1}\otimes d\vartheta+(\mathbbm{1}\otimes \vartheta)(\mathbbm{1}\otimes \vartheta)
        \\
        &=&
        \mathbbm{1}\otimes (-\vartheta\,\vartheta+ \vartheta\, \vartheta)
        \\
        &=&
        0.
    \end{eqnarray*}
We conclude that $R^{\omega^c}=0$.
\end{proof}

According to Section 3.1, every qpb $\omega$ is of the form 
\begin{equation}
    \label{eq.6.1.11}
    \omega=\omega^c+\lambda,
\end{equation}
where $\lambda$ $\in$ $\overrightarrow{\mathfrak{qpc}(\zeta)}$ (see equation  (\ref{2.f24.2})). As in the previous subsection, since $\vartheta^\ast=\vartheta$ and $\ad(\vartheta)=\vartheta\otimes \mathbbm{1}$, it follows that every $\lambda$ is of the form 
\begin{equation}
    \label{eq.6.1.12}
    \lambda(\vartheta)=\mu\otimes \mathbbm{1}\quad \mbox{ where }\quad \mu \,\in\, \Omega^1(\torus^n_\Xi) \quad \mbox{ such that }\quad \mu^\ast=\mu.
\end{equation}
Moreover, the curvature of every qpc is  
\begin{eqnarray*}
    R^\omega=d\omega-\langle \omega,\omega\rangle&=& d(\omega^c+\lambda)-\langle \omega^c+\lambda,\omega^c+\lambda\rangle 
    \\
    &=&
    d\lambda-\langle \omega^c,\lambda\rangle-\langle\lambda,\omega^c\rangle-\langle \lambda,\lambda\rangle.  
\end{eqnarray*}
In this way, we get for $\vartheta=\pi(z)$
\begin{eqnarray*}
    R^\omega(\vartheta)&=&d\lambda(\vartheta)-\langle \omega^c,\lambda\rangle(\vartheta)-\langle\lambda,\omega^c\rangle(\vartheta)-\langle \lambda,\lambda\rangle(\vartheta)
\\
    &=&
(d\mu+\mu\mu)\otimes \mathbbm{1}
\\
    &=&
   d\lambda(\vartheta)-\langle \lambda,\lambda\rangle(\vartheta). 
\end{eqnarray*}
Hence 
\begin{equation}
    \label{eq.6.1.14}
    R^\omega=d\lambda-\langle\lambda,\lambda\rangle.
\end{equation}

On the other hand, the associated right qvb for the $\ad$ $\U(1)$--corepresentation $$E^{\mathfrak{qu}(1)^\#}_\r\cong \torus^n_\Xi$$ has a canonical $\torus^n_\Xi$--basis, given by $$T^\triv: \mathfrak{qu}(1)^\#\longrightarrow P, \qquad w\,\vartheta\longmapsto w\,\mathbbm{1} \qquad \mbox{ with }\qquad w\,\in\, \C.$$ 

Since the associated right qvb for $\ad$ with respect to the quantum principal $\U(1)$--bundle of equation (\ref{new.1}) is canonical isomorphic to the assoacited right qvb for $\ad$ with respect to the quantum principal $\U(1)$--bundle of equation (\ref{eq.6.2}), it follows that  Theorems  \ref{qpc-qlc1}, \ref{qpc-qlc2} and equation (\ref{yangmillsqlc3}) hold in the case studied in this subsection. 

Notice that by equation (\ref{eq.6.1.14}), a qpc $\omega$ different from $\omega^c$ has zero curvature if and only if the following equation holds
\begin{equation}
    \label{eq.6.1.15}
    d\mu+\mu\,\mu=0
\end{equation}
with $\mu\not=0$. The previous equation is exactly equation (\ref{yangmillsqlc2}) and, hence, every qpc with zero curvature is given by  equation (\ref{yangmillsqlc3}). It is worth remarking that now, for the $1$--dimensional non--standard $\ast$--FODC of $U(1)$, the equation that describes qlc's on $E^{\mathfrak{qu}(1)^\#}_\r$ compatible with the metric with zero curvature is the same as the equation that describes qpc's with zero curvature.

The proof of the following statement is completely analogous to that of Theorem \ref{soperatoradjunto}, and will therefore be omitted.

\begin{Theorem}
    \label{soperatoradjunto1}
    The formal adjoint operator 
    \begin{equation*}
     d^{\widehat{S}^{\omega}\star} :E^{\mathfrak{qu}(1)^\#}_\r\otimes_B \Omega^{a+1}(\torus^n_\Xi)\longrightarrow E^{\mathfrak{qu}(1)^\#}_\r\otimes_B \Omega^a(\torus^n_\Xi)
\end{equation*}
of $d^{\widehat{S}^{\omega}}:=\widehat{\Upsilon}_{\mathfrak{qu}(1)^\#} \circ \widehat{S}^\omega \circ \widehat{\Upsilon}^{-1}_{\mathfrak{qu}(1)^\#}$ is given by
    \begin{equation*}
d^{\widehat{S}^{\omega}\star}:=(-1)^{a}\, (\id \otimes_B \star^{-1}_q) \circ\; d^{\widehat{S}^{\omega}}\circ (\id \otimes_B \star_q).
\end{equation*}
\end{Theorem}

Putting together equations (\ref{ec.000}), (\ref{4.3.f29.1}),  (\ref{formaldiff}) and the last theorem, it immediately follows that 
\begin{Theorem}
    \label{YANGMILLS}
    For the $1$--dimensional non--standard $\ast$--FODC of $U(1)$, the non--commutative geometrical Yang--Mills equation takes the form  
    \begin{eqnarray*}
        0&=&(d^{\widehat{\nabla}^\omega\star}-d^{\widehat{S}^{\omega}\star})(\widehat{\Upsilon}_{\mathfrak{qu}(1)^\#}(R^\omega))
        \\
        &=&
        (-1)^a (\id\otimes_B \star^{-1}_q)(d^{\widehat{\nabla}^\omega}-d^{\widehat{S}^\omega})(\id\otimes_B \star_q)(\widehat{\Upsilon}_{\mathfrak{qu}(1)^\#}(R^\omega)).
    \end{eqnarray*}
\end{Theorem}

It is worth recalling that in the case of the {\it classical} $\ast$--FODC of $U(1)$, the non--commutative geometrical Yang--Mills equation is  
\begin{eqnarray*}
    0&=&d^{\widehat{\nabla}^\omega\star}(\widehat{\Upsilon}_{\mathfrak{u}(1)^\#_\C}(R^\omega))
    \\
    &=&
    (-1)^a (\id\otimes_B \star^{-1}_q)d^{\widehat{\nabla}^\omega}(\id\otimes_B \star_q)(\widehat{\Upsilon}_{\mathfrak{qu}(1)^\#}(R^\omega))
\end{eqnarray*}
because $d^{\widehat{S}^{\omega}\star}=0$. By Proposition \ref{canoqpc1}, we obtain that $\omega^c$ is a Yang--Mills qpc; but also we have

\begin{Theorem}
\label{solqpc2}
    Consider the characterization of every qpc $\omega$ and its curvature presented in equations (\ref{eq.6.1.11})--(\ref{eq.6.1.14}). Let $\omega$ be a qpc different from $\omega^c$. Then, $\omega$ is a Yang--Mills qpc if and only if
    \begin{equation}
        \label{ahuevo2}
       d(\star_q(d\mu+\mu\,\mu))+\mu\,\star_q(d\mu+\mu\,\mu)+(-1)^{n-1}\,\star_q(d\mu+\mu\,\mu)\,\mu=0.
    \end{equation}
\end{Theorem}

\begin{proof}
    Let $\omega=\omega^c+\lambda$ $\in$ $\mathfrak{qpc}(\zeta)$, where $\lambda(\vartheta)=\mu\otimes \mathbbm{1}$ with $0\not=\mu$ $\in$ $\Omega^1(\torus^n_\Xi)$. Then, we know that $$R^\omega=d\lambda+\langle\lambda,\lambda\rangle \qquad \mbox{ with }\qquad R^\omega(\vartheta)=d\lambda(\vartheta)-\langle \lambda,\lambda\rangle(\vartheta)=(d\mu+\mu\,\mu)\otimes \mathbbm{1}$$ and $$R^\omega=T^\triv\,(d\mu+\mu\,\mu).$$ Thus (see equation (\ref{formaldiff}))
    
     \begin{eqnarray*}
         0&=&(d^{\widehat{\nabla}^\omega\star}-d^{\widehat{S}^{\omega}\star})(\widehat{\Upsilon}_{\mathfrak{qu}(1)^\#}(R^\omega))
         \\
         &=&
         (d^{\widehat{\nabla}^\omega\star}-d^{\widehat{S}^{\omega}\star})(
         \widehat{\Upsilon}_{\mathfrak{qu}(1)^\#}(T^\triv\,(d\mu+\mu\,\mu)))
         \\
         &=&(d^{\widehat{\nabla}^\omega\star}-d^{\widehat{S}^{\omega}\star})(T^\triv\otimes_{\torus^n_\Xi} (d\mu+\mu\,\mu))
         \\
         &=&-(\id\otimes_{\torus^n_\Xi}\star^{-1}_q) d^{\widehat{\nabla}^\omega}(\id\otimes_{\torus^n_\Xi}\star_q)(T^\triv\otimes_{\torus^n_\Xi} (d\mu+\mu\,\mu))
         \\
         &+&
  (\id\otimes_{\torus^n_\Xi}\star^{-1}_q)d^{\widehat{S}^{\omega}}(\id\otimes_{\torus^n_\Xi}\star_q)(T^\triv\otimes_{\torus^n_\Xi} (d\mu+\mu\,\mu))
         \\
         &=&-(\id\otimes_{\torus^n_\Xi}\star^{-1}_q) d^{\widehat{\nabla}^\omega}(T^\triv\otimes_{\torus^n_\Xi} \star_q(d\mu+\mu\,\mu))
         \\
         &+&
(\id\otimes_{\torus^n_\Xi}\star^{-1}_q)d^{\widehat{S}^{\omega}}(T^\triv\otimes_{\torus^n_\Xi} \star_q(d\mu+\mu\,\mu))
         \\
         &=&
         -(\id\otimes_{\torus^n_\Xi}\star^{-1}_q)(\widehat{\nabla}^\omega(T^\triv)\star_q(d\mu+\mu\,\mu)) 
         \\
         &-&
         (\id\otimes_{\torus^n_\Xi}\star^{-1}_q)(T^\triv\otimes_{\torus^n_\Xi} d(\star_q(d\mu+\mu\,\mu)))
         \\
         &-&
          \,T^\triv\otimes_{\torus^n_\Xi} \star^{-1}_q( \mu\,\star_q(d\mu+\mu\,\mu))
         \\
         &-&
         (-1)^{n-1}\,T^\triv\otimes_{\torus^n_\Xi}\star^{-1}_q(\star_q(d\mu+\mu\,\mu)\,\mu)
         \\
         &=&
         T^\triv\otimes_{\torus^n_\Xi} \star^{-1}_q([-d(\star_q(d\mu+\mu\,\mu))
         \\
         &-&
         \mu\,\star_q(d\mu+\mu\,\mu)-(-1)^{n-1}\,\star_q(d\mu+\mu\,\mu)\,\mu])
     \end{eqnarray*}
    The last equality implies $$\star^{-1}_q([-d(\star_q(d\mu+\mu\,\mu))-\mu\,\star_q(d\mu+\mu\,\mu)-(-1)^{n-1}\,\star_q(d\mu+\mu\,\mu)\,\mu])=0 $$ which is equivalent to equation (\ref{ahuevo2}).
\end{proof}

It is worth mentioning that, for the $1$--dimensional non--standard $\ast$--FODC of $U(1)$, the equation that describes Yang--Mills qlc's (equation (\ref{ahuevo1})) and the equation that describes Yang--Mills qpc's (equation (\ref{ahuevo2})) {\bf are exactly the  same for all $n\geq 2$}. Hence, equation (\ref{ahuevo2}) turns into 
\begin{equation}
        \label{yangmillsqpc3}
        d\mu+\mu\, \mu=0,
    \end{equation}
and we conclude that a qpc $\omega$ is a Yang--Mills qpc if and only $\omega$ has zero curvature. This type of qpc are given by equation  (\ref{yangmillsqlc3}).

\begin{Proposition}
    \label{gaugetorus3}
    Let $\omega=\omega^c+\lambda$ be a Yang--Mills qpc with $\lambda\not=0$. Then, there exists $\F$ $\in$ $\qGG_\YM$ such that  $$\F^{\circledast}\omega^\c=\omega.$$ In other words,  up to an element of $\qGG_\YM$, $\omega^c$ is the only Yang--Mills qpc. 
\end{Proposition}

\begin{proof}
    Let $\omega=\omega^c+\lambda$ be another Yang--Mills qpc, where 
    \begin{eqnarray*}
        0\not=\lambda(\vartheta)&=&\mu\otimes \mathbbm{1}
        \\
        &=&2\,\pi(m_1\,d\mathbbm{U}_1,\cdots ,m_1\,d\mathbbm{U}_n)\otimes\mathbbm{1}\cong  2\,\pi(m_1\,d\mathbbm{U}_1,\cdots ,m_1\,d\mathbbm{U}_n)=\mu
    \end{eqnarray*}
    for some $m_j$ $\in$ $\Z$ with $j=1,...,n$, not all zero. Then, there exists a quantum gauge transformation $$\F: \Omega^\bullet(P)\longrightarrow \Omega^\bullet(P)$$ such that  $\F^{\circledast}\omega^c=\omega.$  This is because, \emph{mutatis mutandis}, the map $\F$ defined in the proof of Proposition \ref{gaugetorus2} is still a quantum gauge transformation for the qpb of this section, and the proof of this is completely similar. 

    Let us check that $\F$ is actually an element of $\qGG_\YM$. Let $\omega'=\omega^c+\lambda'$ be any other qpc. Then, we get $$\omega'(\vartheta)=\omega^c(\vartheta)+\mu'\otimes \mathbbm{1} \quad \mbox{ with }\quad \lambda'(\vartheta)=\mu'\otimes \mathbbm{1},\quad \mu'\,\in\,\Omega^1(\torus^n_\Xi),\quad \mu'=\mu'^\ast.$$  Proposition 4.8 of reference \cite{sald2} implies that 
    \begin{eqnarray*}
        \F^{\circledast}\omega'(\vartheta)=m_{\Omega}(\omega'\otimes \f )\ad(\vartheta) +\f(\vartheta)
        &=&
        m_{\Omega}(\omega'\otimes \f )(\vartheta\otimes \mathbbm{1}) +\f(\vartheta)
         \\
        &=&
\omega'(\vartheta)\f(\mathbbm{1})+f(\vartheta)
        \\
        &=&
        \omega'(\vartheta) +\mu.
        \\
        &=&
        \omega^c(\vartheta)+\mu' +\mu
        \\
        &=&
    \omega^c(\vartheta)+\lambda'(\vartheta)+\lambda(\vartheta)
        \\
        &=&
        \omega^c(\vartheta)+\widetilde{\lambda}(\vartheta),
    \end{eqnarray*}
where $\widetilde{\lambda}=\lambda'+\lambda$. According to equation (\ref{eq.6.1.14}), we obtain $$  R^{\F^{\circledast}\omega'}=d\widetilde{\lambda}-\langle\widetilde{\lambda},\widetilde{\lambda}\rangle.$$ Finally, since $\mu$ is a $1$--degree element of the graded center of $\Omega^\bullet(\torus^n_\Xi)$, we get 
\begin{eqnarray*}
    R^{\F^{\circledast}\omega'}(\vartheta)=d\widetilde{\lambda}(\vartheta)-\langle\widetilde{\lambda},\widetilde{\lambda}\rangle(\vartheta)&=&d\widetilde{\lambda}(\vartheta)+ \widetilde{\lambda}(\vartheta)\,\widetilde{\lambda} (\vartheta)
    \\
    &=&
    (d\mu'+d\mu + (\mu'+\mu)\,(\mu'+\mu))\otimes \mathbbm{1}
    \\
    &=&
    (d\mu' +  \mu'\,\mu')\otimes \mathbbm{1}
    \\
    &=&
    d\lambda'(\vartheta)-\langle \lambda',\lambda'\rangle(\vartheta)
    \\
    &=&
    R^{\omega'}(\vartheta).
\end{eqnarray*}
and it follows that $||R^{\F^{\circledast}\omega'}||^2=||R^{\omega'}||^2$. Hence, $\F$ $\in$ $\qGG_\YM$. 
\end{proof}

Let us summarize the results concerning the two non--commutative Yang--Mills functionals for the $1$--dimensional non--standard $\ast$--FODC of $U(1)$ and compare them with each other in order to provide a {\bf second} answer to the questions raised at the end of the Section 5 for the qpb $\zeta$ of this subsection.  
\begin{enumerate}
    \item The gauge qlc  $\widehat{\nabla}^\omega=\widehat{\nabla}^\omega_{{\mathfrak{qu}(1)^\#}}$  of every qpc $\omega$ is a Yang--Mills qlc on $E^{\mathfrak{qu}(1)^\#}_\r$.
    \item Equation (\ref{ahuevo1}) characterize all Yang--Mills qlc's on $E^{\mathfrak{qu}(1)^\#}_\r$ for every $n\geq 2$. This is equivalent to ask for qlc compatible with the metric with zero curvature and in equation (\ref{yangmillsqlc3}) we characterize all qlc's with zero curvature. 
    \item  Up to an element of $U(E^{\mathfrak{qu}(1)^\#}_\r)$,  $\nabla^\omega$  is the only Yang--Mills qlc. 
    \item Equation (\ref{ahuevo2}) characterize all Yang--Mills qpc's for every $n\geq 2$. This is equivalent to ask for qpc with zero curvature and in equation (\ref{yangmillsqlc3}) we characterize all qpc's with zero curvature. 
    \item The non--commutative analytical Yang--Mills equation (equation (\ref{ahuevo1})) and the non--commutative geometrical Yang--Mills (equation (\ref{ahuevo2})) are the same.
    \item Up to an element of $\qGG_\YM$, $\omega^c$ is the only Yang--Mills qpc.
\end{enumerate}

\begin{Remark}
    \label{gaugeboson2}
    In light of the interpretation of Remark \ref{gaugeboson}, $\omega^c$ (and quantum gauge transformation of this qpc) is a non--commutative photon field moving through the quantum $n$--torus, where $U(1)$ is equipped with $1$--dimensional non--standard differential structure. As we should expect, module (quantum) gauge transformation, there is only one (non--commutative) photon field.
\end{Remark}

\section{Concluding Comments}

As mentioned in the introduction, this paper has two main objectives. The first one is to show that, for a given quantum principal $\G$--bundle over an $n$--dimensional spectral triple, one can always study the gauge Dirac operator associated with any quantum principal connection (see Definition \ref{diracoperator}). This was covered in Section 5.2 and is made possible by Theorem \ref{fgs}. Moreover, as noted in Remark \ref{remayanmills}, in principle, this formulation does not provide additional information about the quantum spaces. However, if the studied quantum vector bundle is the associated quantum vector bundle of some quantum principal bundle corresponding to a finite--dimensional corepresentation, then the gauge Dirac operator of a quantum principal connection is always a suitable and interesting object of study.

The second objective constitutes the core of this paper and consists in presenting the two Yang--Mills functionals that can be defined in non--commutative geometry: the non-commutative analytical Yang--Mills functional (see Definition \ref{yangmills}) and the non--commutative geometrical Yang--Mills functional (see Definition \ref{yangmills1}). The first functional has been studied in detail for many quantum spaces over the past decades \cite{con,lan,lan2}. In contrast, the second functional was originally proposed in \cite{saldym}. It is worth mentioning that in \cite{saldym} we showed that the non--commutative geometrical Yang--Mills functional can be applied to two different classes of qpb's. In the present paper, we have worked with another class of qpb's suitable for the theory (see Remark \ref{assuptions}).

The two Yang--Mills functionals are, in principle, mathematically distinct. The non--commutative analytical Yang--Mills functional measures the squared norm of the curvature of a qlc compatible with the metric, whereas the non--commutative geometrical Yang--Mills functional measures the squared norm of the curvature of a qpc. However, by Theorem \ref{fgs}, the gauge qlc associated with every qpc is compatible with the metric, which suggests that their corresponding {\it field equations} (equations (\ref{ymanalitic}), (\ref{ymgeometric}), respectively) may share solutions.

According to Section 6, there are cases in which both non--commutative Yang--Mills functionals lead to the same results, although there are also cases in which the results of the two functionals differ. In the case of quantum principal $\U(1)$--bundles and trivial quantum vector bundles, both defined over the non--commutative $n$--torus, the fact that both functionals are equivalent or not depends on which bicovariant $\ast$--FODC of $U(1)$ is used. One might think that the {\it classical} $\ast$--FODC of $U(1)$ is the correct one to perfectly couple both Yang--Mills functionals; however, this is not the case. The bicovariant $\ast$--FODC of $U(1)$ that perfectly couples the two Yang--Mills functionals is precisely the one presented in Section 2.3.

 It is worth mentioning that in reference \cite{sald5} we appreciate the same phenomenon, that is, only by using the $1$--dimensional non--standard bicovariant $\ast$--FODC of $U(1)$ we can recreate previous results of the electromagnetic theory in the Moyal--Weyl algebra. In addition, in reference \cite{sald3} we applied our theory on the quantum Hopf fibration, also known as the $q$--deformed Dirac monopole bundle. In this situation, the bicovariant  $\ast$--FODC of $U(1)$ used is $1$--dimensional, but again, it is not the {\it classical} one (although, it is not the one presented in Section 2.3 neither).

This may be due to the fact that $U(1)$ is an abelian group, and its {\it classical} $\ast$--FODC is also {\it too commutative} (there are no differential forms of degree higher than $1$, and $0$--forms and $1$--forms commute); so, {\it it does not fit very well} within the non--commutative geometrical framework. Thus,  {\it no--classical} $\ast$--FODCs of $U(1)$ {\it fit better} in the non--commutative geometrical framework.

The study of specific examples for the Lie groups $SU(N)$ with $N\geq 2$ is a strong motivation to continue this line of research, which will be addressed in forthcoming publications.

In differential geometry, the concept of a principal $H$--bundle forms the foundation of Yang--Mills theory and it has the distinctive property of encoding the notion of {\it symmetry} through the action of the Lie group $H$. The non--commutative geometrical Yang--Mills functional has the advantage of being based on the concept of a quantum principal $\G$--bundle, which encodes the notion of {\it non--commutative symmetry} through  the quantum group $\G$.

It is worth mentioning that, in the way we have formulate our theory, the non--commutative geometrical Yang--Mills functional, as well as the non--commutative geometrical Yang--Mills equation, are not limited to the case presented in this paper (see reference reference \cite{saldym}).  In particular, from a {\it completely mathematical} point of view,  we can apply our formulation when the quantum group of the bundle is the quantum $SU(N)$ group with any finite--dimensional bicovariant $\ast$--FODC's. This will be addressed in forthcoming publications. 

For instance, in \cite{sald5} we showed that if the base space is the algebra of $\mathbb{C}$--valued smooth functions on a (compact) manifold, then for any one--dimensional non--standard $\ast$--FODC on $U(1)$ (for example, the one described in Section~2.3), the non--commutative geometrical Yang--Mills equation coincides with the \emph{classical} Yang--Mills equation. However, it is likely that this result holds only for the group $U(1)$.

To finalize this paper, let us talk about our the definition of quantum principal connections. For this paper we have decided to use the condition $\omega(\theta^\ast)=\omega(\theta)^\ast$, first, because it is the {\it standard} definition in Durdevich's formulation, as the reader can verify in references \cite{micho1,micho2,micho3,stheve}. Second, because in this way, Theorem \ref{fgs} holds for every qpc. Finally, because we have related qpc's to gauge boson fields, as in differential geometry, and one of the {\it physical requirements} for gauge boson fields is that they must be real fields. The condition $\omega(\theta^\ast)=\omega(\theta)^\ast$ ensures this.

\appendix

\section{Some Similarities and Differences between Durdevich’s formulation and the Brzezi\'nski–Majid formulation}

Consider a qpc $\zeta=(P,B,\Delta_P)$ as defined in Section 3.1. Since  in the present paper we always work with quantum groups, it can be proven that the quotient map $$\widetilde{\beta}:P\otimes_B P\longrightarrow P\otimes G$$ induced by $\beta$ is always bijective \cite{libro}. In this setting, the triple $(P,B,\Delta_P)$ is commonly referred to as a Hopf--Galois extension. In other words, in this paper, qpb's are always Hopf--Galois extensions.

In the Brzezi\'nski–Majid formulation of qpb's presented in \cite{libro,br1,br2}, a qpb is defined as a Hopf--Galois extension $\zeta=(P,B,\Delta_P)$ with
\begin{enumerate}
    \item A FODC $(\Omega^1(P),d)$ over $P$.
    \item A FODC $(\Omega^1(B),d)$ over $B$, where $\Omega^1(B)\subseteq \Omega^1(P)$.
    \item A bicovariant FODC $(\Gamma,d)$ over $G$
\end{enumerate}
such that the following sequence of $P$--bimodules
\begin{equation}
\label{3.f1.4.5}
0\longrightarrow  P\,\Omega^1(B)\,P \lhook\joinrel\relbar\joinrel\rightarrow  \Omega^1(P) \xlongrightarrow{\pi_\V} \Vert^1 P:=P\otimes \mathfrak{qg}^\# \longrightarrow 0
\end{equation}
is exact, where  $\pi_V: \Omega^1(P)\longrightarrow \Vert^1 P$ is defined as in Proposition \ref{seq} (here, we are imposing that the map $\pi_V$ is well--defined).

Here, we can appreciate some difference with the definition of qpb and differential calculus given is Section 3:
\begin{enumerate}
    \item The involution $\ast$ is not necessary.
    \item The space of horizontal $1$--forms is defined as $P\,\Omega^1(B)\,P$.
    \item The space of base $1$--forms is always generated by $B$.
    \item The exactness of the sequence (\ref{3.f1.4}) is condition.
    \item There are no condition for degrees higher or equal to $2$.
\end{enumerate}

Of course, there are a lot of examples for which both definitions coincide, for instance, the qpb's of Section 3.2. However, there are cases in which they both definitions do not coincide, for instance, Example 2.1 of reference \cite{appendix}.

In the Brzezi\'nski–Majid formulation of qpb's, a qpc  is defined as a linear map $$\omega^\times:G\longrightarrow \Omega^1(P)$$ such that  $$\Delta_R\circ \omega^\times=(\omega\otimes \id_{\mathfrak{qg}^\#})\circ \Ad\qquad \mbox{ and}\qquad (\pi_V\circ \omega^\times)(g)=\mathbbm{1}\otimes  \pi(g)$$ for all $g$ $\in$ $G$, where $$\Delta_R:\Omega^1(P)\longrightarrow\Omega^1(P)\otimes G, \qquad p_1dp_2\longmapsto p^{(0)}_1dp^{(0)}_2\otimes p^{(1)}_1p^{(1)}_2 $$ with $\Delta_P(p)=p^{(0)}\otimes p^{(1)}$. It is worth mentioning that, as in Theorem 12.8 of reference \cite{stheve}, every $\omega^\times$ can be regarded as a left $P$--module monomorphism from $\mathrm{Ver}^1\,P$ to $\Omega^1(P)$ such that sequence (\ref{3.f1.4.5}) splits \cite{libro}.

Quantum principal connections $\omega^\times$ in the previous sense, and qpc's $\omega$ as in Definition \ref{qpc's} are related by (\cite{micho2,libro}) $$\omega^\times=\omega\circ \pi,$$ where $\pi:G\longrightarrow \mathfrak{qg}^\#$ is the quantum germs map (see equation (\ref{2.f14})). 

The covariant derivative of $\omega^\times$ is defined as $$D^{\omega^\times}:P\longrightarrow \Omega^1(P),\qquad p\longmapsto dp-p^{(0)}\,\omega^\times(p^{(1)}).$$ It is worth noticing that the previous definition coincides with the definition of $D^\omega$ in degree $0$ (see equation (\ref{2.f30})). However, $D^{\omega^\times}$ is not defined for higher degrees, and the operator $\widehat{D}^\omega$ does not exists (although, of course, it is easily defined). Moreover,
we say  that a qpc $\omega^\times$ is \emph{strong} if $$\Im(D^{\omega^\times})\subseteq \Omega^1(B)\,P.$$

The curvature of a qpc $\omega^\times$ is defined by (\cite{libro})
$$r^{\omega^\times}:G\longrightarrow \Omega^1(P),\qquad g\longmapsto d\omega^\times(g)+\omega^\times(\pi(g^{(1)}))\,\omega^\times(\pi(g^{(2)})$$ It is easy to check that, if $\omega$ is multiplicative in the sense of Section 3 (for instance, in the dualization of a \emph{classical} principal connection), one obtains $$r^{\omega^\times}=R^\omega\circ \pi.$$
Otherwise, the two curvatures are generally different. This is problematic, since in general there exist qpc's $\omega$ that are not multiplicative \cite{micho2}.

In the Brzezi\'nski–Majid formulation of qpb's, given $\delta^V$ a finite--dimensional $\G$--corepresentation, the associated quantum vector bundle is defined as the finitely generated projective left $B$--module $$E^V:=P\,\square^{G}\, V:=\{s\in P\otimes V\mid (\Delta_P\otimes \id_V)(s)=(\id\otimes \delta^V_L)(s)\},$$ where $$\delta^V_L:V\longrightarrow G\otimes V, \qquad v\longmapsto S^{-1}(v^{(1)})\otimes v^{(0)},$$ with $\delta^V(v)=v^{(0)}\otimes v^{(1)}$, and the $B$--left module structure is given by $$b\cdot (p\otimes v)=b\,p\otimes v.$$

In accordance with Proposition $6.1$ of reference \cite{br2}, for the natural coaction on the dual space $V^{\#}$ of $V$, we have  $$E^V=P\,\square^{G}\, V^{\#} \cong E^V_\l$$ as left $B$--modules.

In addition, if $\omega^\times$ is strong, we can define a gauge quantum linear connection on $E$ by means of (\cite{libro}) $$\nabla^{\omega^\times}_V:E^V\longrightarrow \Omega^1(B)\otimes E^V,\qquad  (p\otimes v)\longmapsto D^{\omega^\times}p\otimes v.$$ It is worth mentioning that, in this setting, the gauge qlc is defined exclusivity for strong qpc's $\omega^\times$; while in the theory showed in Section 4, the gauge qlc is defined for every qpc $\omega$. However, it is worth noticing that, if $\omega$ is a regular qpc, then $\omega^\times=\omega\circ \pi$ is strong qpc. Furthermore, in Durdevich's formulation, there is a general theory of canonical Hermitian structures on associated qvb's and of course, Theorem  \ref{fgs} holds for every qpc $\omega$.

To finalize this presentation of the Brzezi\'nski–Majid formulation of qpb's, let us talk about the quantum gauge group. In this formulation, the quantum gauge group $\qGG^0$ is defined as
$$\qGG^0=\{\F:P\longrightarrow P \mid \F \mbox{ is a left } B\mbox{-- module isomorphim}$$ $$\mbox{ such that }\,\F(\mathbbm{1})=\mathbbm{1},\;\; \Delta_P\circ \F=(\F\otimes \id_{G})\circ \Delta_P.$$ This group is isomorphic to the group of all convolution--invertible maps (\cite{br1}) $$ \f:G\longrightarrow P$$ such that $$\f(\mathbbm{1})=\mathbbm{1} \;\;\mbox{ and }\;\; (\f\otimes \id_G)\circ \Ad=\Delta_P\circ \f.$$
Additionally, $\qGG^0$ acts on the space $\mathfrak{qpc}(\zeta^\times)$ of all qpc's $\omega^\times$ by means of
\begin{equation}
    \label{ap.1}
    \omega^\times \longmapsto \F\omega^\times:=\f\;\widetilde{\ast}\; \omega^\times \;\widetilde{\ast}\; \f^{-1}+ \f\;\widetilde{\ast}\; (d\f^{-1}),
\end{equation}
where $\widetilde{\ast}$ denotes the convolution product, and the following formula holds (\cite{br1}) 
\begin{equation}
    \label{ap.2}
    r^{\F\omega^\times}=\f\;\widetilde{\ast}\; r^{\omega^\times} \;\widetilde{\ast}\; \f^{-1}.
\end{equation}

Notice that the definition of the quantum gauge group $\qGG$ in Durdevich's formulation (see Definition (\ref{qgg})) is a  generalization at the level of differential calculus of  $\qGG^0$ taking into account the $\ast$ operation. In particular, $\qGG$  is isomorphic to a subgroup of the group of all convolution--invertible maps (see reference \cite{appendix}) $$\f:\Gamma^\wedge\longrightarrow \Omega^\bullet(P).$$

In Durdevich's framework, the action of $\qGG$ on qpc's as in equation (\ref{ap.1})  is not well--defined. Furthermore, even if we extend the domain of $\omega$ to $G$ by using the quantum germs map $\pi$ , the induced action on the curvature $$\f\ast R^\omega \ast \f^{-1} $$ remains ill-defined since $R^\omega$ has domain $\mathfrak{qg}^\#$. However, the action of  $\qGG$ on qpc's as in equation (\ref{actionqpc}) is well--defined, and the action on $R^\omega$ too (see Section 4.2 of reference \cite{sald2}).

To conclude this appendix, let us clarify more explicitly our reasons for developing this work within Durdevich's formulation.

\begin{enumerate}
    \item The principal reason is related with the curvature.  In the Brzezi\'nski–Majid formulation of qpb's, we have 
    \begin{enumerate}
    \item Since $r^{\omega^\times}$ has domain $G$, the embedded differential (see equation (\ref{embbededdifferential})) is not required. Consequently, the operator $S^\omega$ does not exist in this formulation. This point is of great importance, since, as shown in Section~6.2, the operator $S^\omega$ is responsible for ensuring that the critical points of the two non--commutative Yang--Mills equations coincide in the qpb used.

    \item Since $r^{\omega^\times}$ has domain $G$ rather than $\mathfrak{qg}^\#$, the $B$--valued inner product introduced in Section~4 cannot be used; unless $G$ is finite--dimensional, which leaves us very limited in examples. Therefore, in general, the expression $\|r^{\omega^\times}\|^{2}$ is not well--defined.

    \item If one dualizes, via pull--back, the definition of the curvature $\Omega^\omega$ of a principal connection in differential geometry, one finds that the curvature has domain equal to the dual of the Lie algebra of the group \cite{nodg}, exactly as in the case of $R^\omega$. In addition, $\Omega^\omega$ is a basic differential $2$--form of type $\ad^{\mathrm{class}}$ \cite{nodg}, and thus, in physical terms, $\Omega^\omega$ can be interpreted as a field strength. This is consistent with the fact that $R^\omega \in \Mor(\ad,\Delta_{\Hor})$, and hence, from a physical perspective, $R^\omega$ can be interpreted as a quantum field strength. Therefore, $R^\omega$, rather than $r^{\omega^\times}$, appears to be the most natural choice for the purposes of this work.
\end{enumerate}

    \item In the Brzezi\'nski–Majid formulation of qpb's, the space $E^V$ does not have, at first sight, canonical generators (as the $\{ T^\l_k\}$ maps) that allow one to perform explicit computations easily. However, since $E^V \cong E^V_{\l}$, such generators can in principle be obtained.

    \item In the Brzezi\'nski–Majid formulation of qpb's, there is no general theory of canonical Hermitian structures on $E^V$ valid for every qpb and every finite--dimensional corepresentation $\delta^V$. Moreover, in this setting, an analogue of Theorem~\ref{fgs} does not exist. As shown in Section~5.2, this theorem is of central importance, since it provides the link between Connes' Dirac theory and the theory of qpb's. Nevertheless, as in the previous point, it is plausible that by following the approach of \cite{sald2} and using the identification $E^V \cong E^V_{\l}$, an analogous result could be obtained.

    \item In the Brzezi\'nski–Majid formulation of qpb's, the gauge quantum linear connection is defined only for strong quantum principal connections. In contrast, in Durdevich's formulation it is defined for every quantum principal connection, without the need to impose additional conditions.

    \item Durdevich's formulation allows one to work simultaneously with $(E^V_{\l},\nabla^\omega_V)$ and $(E^V_{\r},\widehat{\nabla}^\omega_V)$, thanks to the existence of the covariant derivatives $D^\omega$ and $\widehat{D}^\omega$.

    \item In \cite{sald2}, an example of a quantum principal bundle is presented for which the action defined in equation~(\ref{ap.1}) is trivial, whereas the action defined in equation~(\ref{actionqpc}) is transitive. This shows that the two actions are genuinely different and suggests that the action given by equation~(\ref{actionqpc}) is, in some sense, the more general one.\\
\end{enumerate}

\begin{center}
    {\bf{Declarations}:}
\begin{enumerate}
    \item Ethics approval and consent to participate: applicable.
    \item  Consent for publication: applicable.
    \item  Availability of data and materials: applicable.
    \item  Competing interests: Not applicable.
    \item  Funding: Center of Research in Mathematics, CIMAT. 
     \item  Authors' contributions: There is only one author.
     \item Acknowledgments: I would like to thank Professor Stephen Sontz for his comments and for the time he devoted to this paper.
     \item Dedication: In memory of Arturo and his family.
\end{enumerate}
\end{center}

\address{
Center of Research in Mathematics, CIMAT,\\
Jalisco S/N Col.Valenciana, C.P. 36023, Guanjuato, M\'exico\\
\email{gustavo.saldana@cimat.mx}\\
}

\end{document}